\documentclass[reqno,a4paper]{amsart}

\usepackage{mathrsfs,amsmath,amssymb,graphicx}
\usepackage[unicode]{hyperref}
\usepackage{latexsym}
\usepackage{layout}
\usepackage[english]{babel}
\usepackage{color}
\usepackage{fancyvrb}
\usepackage{amsmath,amscd}
\usepackage{txfonts}
\usepackage{enumitem} 
\usepackage{float}
\usepackage{multicol} 

\usepackage{subcaption}

\usepackage{transparent}

\usepackage[normalem]{ulem} 

\usepackage{tikz}
\usetikzlibrary{calc}
\usetikzlibrary{matrix,arrows,decorations,chains,positioning}

\makeatletter
\@namedef{subjclassname@2020}{\textup{2020} Mathematics Subject Classification}
\makeatother

\theoremstyle{plain}
\newtheorem*{theorem*}{Theorem} 
\newtheorem{theorem}{Theorem}[section]
\newtheorem{lemma}[theorem]{Lemma}

\newtheorem{corollary}[theorem]{Corollary}

\theoremstyle{definition}
\newtheorem{definition}{Definition}[section]
\newtheorem{remark}[theorem]{Remark}
\newtheorem{question}{Question}[section]

\newcommand{\tA}{\mathcal{A}}
\newcommand{\tD}{\mathcal{D}}
\newcommand{\tF}{{\mathcal F}}
\newcommand{\slice}{{\mathcal S}}
\newcommand{\annulus}{{\mathsf A}}
\newcommand{\mobius}{{\mathsf M}} 
\newcommand{\disk}{{\mathsf D}}
\newcommand{\torus}{{\mathsf T}}
\newcommand{\surface}{{\mathsf S}}

\newcommand{\HK}{{\rm HK}}
\newcommand{\Mo}[1]{{\mathbf{M}_{#1}}}

\newcommand{\LLtwo}[1]{{\mathbf{L}^-_{#1}}}
\newcommand{\HKA}{{\rm HK_\annulus}}

 \newcommand{\Sbb} {\mathbb S}
\newcommand{\sphere} {\Sbb^3}

\newcommand{\Compl}[1]{E(#1)}
\newcommand{\Complprime}[1]{E(#1')}

\newcommand{\snode}{\bullet}
\newcommand{\hnode}{\circ}

\newcommand{\bm}{\partial_M}
 
\newcommand{\familyLLone}{\left\{ \big(\sphere,\mathbf{L}^\circ_n\big) \right\}_{n\in\mathbb{Z}}}
\newcommand{\familyLLvariant}{\left\{\big(\sphere,\mathbf{L}^\odot_n\big)\right\}_{n\in\mathbb{Z}}}
\newcommand{\familyLLtwo}{\left\{\big(\sphere,\mathbf{L}^-_n\big)\right\}_{n\in\mathbb{Z}}}
\newcommand{\familyMo}{\left\{(\sphere,\mathbf{M}_n)\right\}_{n\in\mathbb{Z}}}
\newcommand{\familyE}{\left\{(\sphere,\mathbf{E}_n)\right\}_{n\in\mathbb{Z}}}
\newcommand{\pairLLone}[1]{(\sphere,\mathbf{L}^\circ_{#1})}
\newcommand{\pairLLvariant}[1]{(\sphere,\mathbf{L}^\odot_{#1})}
\newcommand{\pairLLtwo}[1]{(\sphere,\mathbf{L}^-_{#1})}
\newcommand{\pairMo}[1]{(\sphere,\mathbf{M}_{#1})}
\newcommand{\pairE}[1]{(\sphere,\mathbf{E}_{#1})}
\newcommand{\pairfourone}{(\sphere, 4_1)}
\newcommand{\pairfiveone}{(\sphere, 5_1)}
\newcommand{\pairfivetwo}{(\sphere, 5_2)}
\newcommand{\pairsixone}{(\sphere, 6_1)}
\newcommand{\pairsixten}{(\sphere, 6_{10})}

\newcommand{\pairprime}{(\sphere,\HK')}
\newcommand{\pair}{(\sphere,\HK)}
\newcommand{\pairA}{(\sphere,\HK_\annulus)} 

\newcommand{\hopf}{{\bf h}}
\newcommand{\knot}{{\bf k}}
\newcommand{\link}{{\bf l}}
\newcommand{\EM}{{\bf em}}
\newcommand{\charM}{{\Lambda_M}}
\newcommand{\charE}{{\Lambda_{\textsc{ext}}}}
\newcommand{\anndiag}{{\Lambda_{\textsc{hk}}}}
\newcommand{\anndiagprime}{{\Lambda_{\textsc{hk}'}}}

\newcommand{\cout}[1]   {}

\newcommand{\rnbhd}[1]{\mathfrak N(#1)}

\newcommand{\openrnbhd}[1]{\mathring{\mathfrak N}(#1)}

\definecolor{mygray}{rgb}{0.92,0.92,0.92}

\numberwithin{equation}{section}
\numberwithin{figure}{section}

\title[Annulus Configuration]{Annulus configuration in handlebody-knot exteriors}
 
\author{Yi-Sheng Wang}
\address{National Sun Yat-sen University, Kaohsiung 804, Taiwan}
\email{yisheng@mail.nsysu.edu.tw}

\date{\today}

\begin{document}
 
\thanks{
The author gratefully acknowledges the support 
from MoST (grant no. 110-2115-M-001-004-MY3), Taiwan.}

\begin{abstract}
In contrast to classical knots,
the knot type of a 
genus two handlebody-knot is not determined by its exterior,
and it is often a challenging task to distinguish handlebody-knots with homeomorphic exteriors.
The present paper considers an invariant (the annulus diagram), defined via Johannson's characteristic submanifold
theory and the Koda-Ozawa classification for essential annuli, and demonstrates its capability to distinguish such handlebody-knots; 
particularly, the annulus diagram is 
able to differentiate members in   
the handlebody-knot families given by Motto and Lee-Lee. 
\end{abstract}

\maketitle
 
\section{Introduction}\label{sec:intro}

A genus $g$ \emph{handlebody-knot} $\pair$ is a
genus $g$ handlebody $\HK$ embedded in an oriented $3$-sphere $\sphere$, and  
two handlebody-knots are \emph{equivalent} or of the same \emph{knot type} if they are ambient isotopic.
While the theory of genus \emph{one} handlebody-knot is equivalent to the study of classical knots, 
higher genus handlebody-knots behave quite differently from
classical knots. For instance, in classical knot theory, the Gordon-Luecke theorem \cite{GorLue:89} asserts that the knot type of a knot is determined by
the homeomorphism type of its exterior, yet
the statement does not hold in higher genus case. 
The first genus $g\geq 3$ counterexample is discovered by Suzuki \cite{Suz:75}, 
and later several infinite families of inequivalent
genus two handlebody-knots with homeomorphic exteriors
are constructed by Motto \cite{Mott:90} and Lee-Lee \cite{LeeLee:12}. The knot types of some genus two handlebody-knots though, for example, $\pairfourone$ in the Ishii-Kishimoto-Moriuchi-Suzuki knot table \cite{IshKisMorSuz:12}, are determined by their exteriors. 
The present work focuses on genus two handlebody-knots, abbreviated to handlebody-knots hereinafter.

While an infinite family of 
handlebody-knots with homeomorphic exteriors 
can be generated quite easily with 
a twist construction \cite{Mott:90} (see Sec.\ \ref{subsec:hk_family}), 
the real challenge is to determine whether handlebody-knots so constructed are mutually inequivalent; 
no computational invariant
capable to distinguish \emph{infinitely} 
many such handlebody-knots seems to be known.
\cite{BelPaoWan:20a} develops a 
computational invariant, based on counting homomorphisms on
the knot group\footnote{The fundamental group
of a handlebody-knot exterior.}, 
able to differentiate \emph{finitely} 
many inequivalent handlebody-knots 
with homeomorphic exteriors, yet it cannot cope with
an infinite family of such handlebody-knots. 

\cite{Mott:90} proves the mutual inequivalence of Motto's handlebody-knots by studying the mapping class group of their exteriors, whereas to differentiate 
Lee-Lee's handlebody-knots, \cite{LeeLee:12} carries 
out a detailed analysis on certain essential annuli 
in their exteriors. 
Either case makes essential use of 
annulus configuration of handlebody-knot exteriors. 
\cite{Wan:22p} shows that the configuration of annuli 
in a handlebody-knot exterior can be encoded 
in a labeled diagram, called \emph{the annulus diagram}, and provides a classification for such diagrams. 

One purpose of the present paper is to compute 
the annulus diagrams of Motto's and Lee-Lee's handlebody-knots, and show that their inequivalence 
can be detected by the annulus diagram 
(Theorems \ref{teo:ann_diag_motto}, \ref{teo:ann_diag_leeleeone}, and \ref{teo:ann_diag_leeleetwo}).  
The annulus diagram provides a general framework 
to describe how the annulus configuration of 
a handlebody-knot exterior may differ between
inequivalent handlebody-knots with homeomorphic exteriors,
and can be applied to other infinite families; 
we demonstrate this by constructing a new infinity family of handlebody-knots with homeomorphic exteriors and proving 
their mutual inequivalence by the annulus diagram (Theorem \ref{teo:ann_diag_leeleevariant}). 

The definition of the annulus diagram, built on
Johannson's characteristic submanifold theory
and the Koda-Ozawa classification of essential annuli, 
is reviewed in Sec.\ \ref{sec:annulus},
and the annulus diagrams of Motto's and Lee-Lee's
handlebody-knots are computed in Sec.\ \ref{sec:families} 
after a brief review of the twist operation that produces
them.

While Sec.\ \ref{sec:families} describes
how handlebody-knot exteriors \emph{fail} to determine 
the knot type of a handlebody-knot, 
Sec.\ \ref{sec:classification} focuses on the positive,
and discusses  
\emph{under what condition 
the handlebody-knot exterior does determine the knot type}.
We show that the knot types of handlebody-knots 
with certain annulus diagrams are determined by 
their exteriors (Theorems \ref{teo:circle_stick}, \ref{teo:theta}); 
particularly, as a corollary of 
Theorem \ref{teo:theta} and \cite[Theorems $1.5$]{Wan:22p},
we obtain that, 
if the handlebody-knot exterior 
admits three non-isotopic, 
non-separating essential annuli and no essential tori,
then it determines the knot type of $\pair$. 
Examples of such include $\pairfourone,\pairsixten$ in the 
knot table of \cite{IshKisMorSuz:12}. 
In closing, we construct an infinite family of 
inequivalent handlebody-knots, showing that in general,
even together with the annulus diagram, the handlebody-knot
exterior is not sufficient to determine the
handlebody-knot.  

 


\section{Annulus Diagram}\label{sec:annulus}

Throughout the paper we work in the piecewise linear category.
Given a subpolyhedron $X$ of a $3$-manifold $M$, 
we denote by $\mathring{X}$, $\mathfrak{N}(X)$ 
and $\bm X$ the interior, 
a regular neighborhood and
the frontier of $X$ in $M$.
By the exterior $\Compl X$ of $X$ in $M$, 
we understand the complement of $M-\openrnbhd{X}$ if $X$ has codimension greater than zero, and  
is the closure of $M-X$ otherwise. 
Submanifolds of a manifold $M$ are assumed to be 
proper and in general position except in some obvious cases
where submanifolds are in $\partial M$.  
A surface which is not a disk or sphere
in a $3$-manifold $M$ is essential if 
it is incompressible, $\partial$-incompressible, and non-boundary parallel; an essential disk in $M$ 
is one that does not cut off a $3$-ball from $M$.
When $M$ is a handlebody, an essential disk in $M$ is also called a \emph{meridian} disk. 
By unique, we understand \emph{unique, up to isotopy}. 
An \emph{atoroidal} 
$3$-manifold is one that contains no essential tori. 
A pair $(\sphere,X)$ denotes an embedding 
of $X$ in $\sphere$;
of the greatest interest here is 
the case where $X=\HK$ is a genus two handlebody.  
Unless otherwise specified, all handlebody-knots $\pair$ are assumed to be atoroidal, namely,
their exteriors $\Compl\HK$ being atoroidal. 
By Thurston's hyperbolization theorem, such an exterior either contains an essential annulus or is hyperbolic.

\subsection{Characteristic diagram}
Given an atoroidal, irreducible, $\partial$-irreducible,  compact, oriented $3$-manifold $M$,   
a codimension-zero submanifold $X\subset M$
is admissibly \emph{fibered} if it can either be \emph{Seifert fibered}
with $X\cap \partial M$ consisting of some fibers or be \emph{I-fibered} with 
$X\cap \partial M$ being the two lids of the I-bundle, where a \emph{lid} 
of an I-bundle $\pi:X\rightarrow B$ is a component of the closure of $\partial X-\pi^{-1}(\partial B)$. 
A codimension-zero submanifold $X\subset M$ is \emph{simple} if 
every essential annulus $\annulus\subset X$ not meeting the frontier $\bm X$ is parallel to an annular component $\annulus'$ of $\bm X$---namely, $\annulus$ 
cuts off a submanifold 
of $X$ that admits an I-bundle structure 
with two lids being $\annulus$ and $\annulus'$. 

By Johannson's characteristic submanifold theory \cite{Joh:79},
there exists a unique surface $\surface\subset M$
consisting of essential annuli such that  
\begin{enumerate}
\item the closure of each component of the complement $M-\surface$
is either simple or admissibly fibered, and
\item removing any component of $\surface$ 
causes the first condition to fail.
\end{enumerate}
Components of $\surface$ are called \emph{characteristic} annuli of $M$, and 
the \emph{characteristic diagram} $\charM$ is 
defined to be the graph given by 
assigning to each admissibly fibered (resp.\ simple) component of $M-\surface$ a solid (resp.\ hollow) node, 
and to each component $\annulus$ of $\surface$ 
an edge that connects the node(s) representing component(s) of $M-\surface$ whose closure(s) contains/contain $\annulus$.
For instance, the characteristic diagram of the 
exterior of the handlebody-knot $\pairfourone$ 
is \raisebox{-.3\height}{\includegraphics[scale=.1]{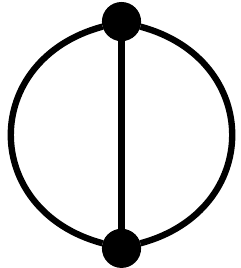}}; 
conversely, any handlebody-knot whose exterior has the characteristic diagram is equivalent to $\pairfourone$ 
\cite[Theorem $1.5$]{Wan:22p}. 
Unlike knots, however, the handlebody-knot exterior 
is in general insufficient 
to distinguish inequivalent handlebody-knots \cite{Mott:90}, \cite{LeeLee:12}, \cite{BelPaoWan:20a}.  
 
\subsection{Types of essential annuli}
Recall that \cite{KodOzaGor:15} and \cite{FunKod:20} classify 
essential annuli in an atoroidal handlebody-knot exterior 
into six types: a type $2$ annulus $\annulus$
is characterized by the property that exactly one component of $\partial \annulus$
bounds a disk $\disk$ in $\HK$; $\annulus$ is said to be of \emph{type $2$-$1$} (resp.\ of \emph{type $2$-$2$})
if $\disk$ is \emph{non-separating} (resp.\ \emph{separating}). 
We use the notation
$\hopf_i$ for a type $2$-$i$ annulus, $i=1,2$, as it is 
also called a Hopf type annulus. 
A \emph{type $3$-$2$} (resp.\ \emph{type $3$-$3$}) annulus $\annulus$ is characterized by the property that 
components of $\partial \annulus$ do not bound disks in $\HK$, and are \emph{parallel} (resp.\ \emph{non-parallel}),
and there exists, up to isotopy, a unique non-separating (resp.\ separating) disk $\disk\subset \HK$ 
disjoint from $\partial \annulus$ \cite[Lemma $2.3$]{LeeLee:12}, \cite[Lemmas $2.1$, $2.3$]{FunKod:20}, \cite[Lemma $2.9$]{Wan:23}. 

A type $3$-$2$ annulus $\annulus$ 
can be further classified into two subtypes: 
if $\annulus$ is essential 
in the exterior of the solid torus $\HK-\openrnbhd{\disk}$,
it is of type $3$-$2$i and otherwise is of type $3$-$2$ii. 
In addition, since $\annulus$ 
cuts off a solid torus $V$ from $\Compl\HK$, 
we define the \emph{slope} of $\annulus$ to be 
the slope of the core of $\annulus$ 
with respect to $(\sphere, V)$,
and denote by $\knot_\ast(r)$ a type $3$-$2\ast$ annulus with a slope of $r$, $\ast=1,2$; the essentiality of $\annulus$ implies $r$ is neither integral nor $\infty$.

Similarly, there is a finer classification for 
a type $3$-$3$ annulus $\annulus$. 
Let $l_1,l_2$ be components of $\partial\annulus$.
Then the unique separating disk $\disk$ cuts $\HK$ into two solid tori $V_1,V_2$ with $l_i\subset \partial V_i$, $i=1,2$. 
The \emph{slope pair} of $\annulus$ is then defined to the unordered pair $(r,s)$
with $r,s$ by the slopes of $l_i$, 
with respect to $(\sphere, V_i)$, $i=1,2$. 
Denote by $\link(r,s)$ a type $3$-$3$ annulus with a slope pair of $(r,s)$. By \cite{Wan:23}, the pair $(r,s)$ is either of the form $(\frac{p}{q},\frac{q}{p})$, $pq\neq 0$ or of the form $(\frac{p}{q},pq)$, $q \neq 0$, where $p,q\in\mathbb{Z}$.

A \emph{type $4$-$1$} annulus is characterized by the property that
components of $\partial \annulus$ are parallel and \emph{no} essential disks in $\HK$ disjoint from $\partial \annulus$ exist. $\annulus$ is necessarily separating 
and cuts off a solid torus whose core in $\sphere$ is an Eudave-Mu\~noz knot (see \cite[Proof of Theorem $3.3$]{KodOzaGor:15} or
Proof of Lemma \ref{lm:typefourone_exclusion}). 

The \emph{annulus diagram} $\anndiag$ of $\pair$ 
is defined to be the characteristic diagram $\charE$ of $\Compl\HK$
together with a labeling that assigns to each edge a label $\hopf_i, \knot_i(r), \link(r,s)$ or $\EM$, depending
on the type of the annulus it represents, where $i=1,2$ and $r,s\in\mathbb{Q}$.

\subsection{Examples}

\begin{figure}[t]
\begin{subfigure}{.45\linewidth}
\center
\includegraphics[scale=.1]{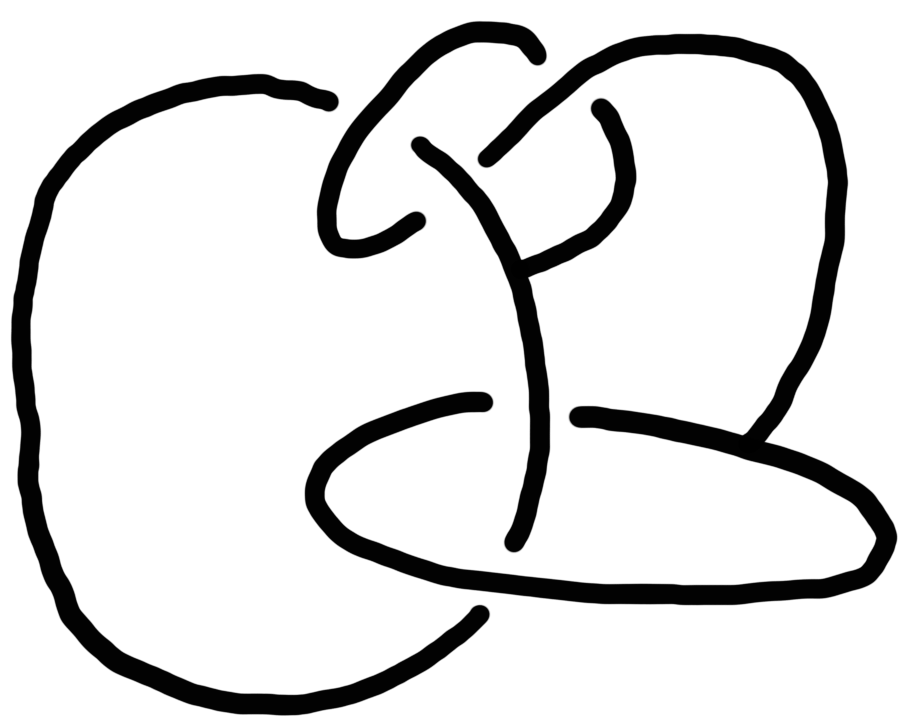}
\caption{$\pairfiveone$.}
\label{fig:hkfiveone}
\end{subfigure}
\begin{subfigure}{.45\linewidth}
\center
\includegraphics[scale=.1]{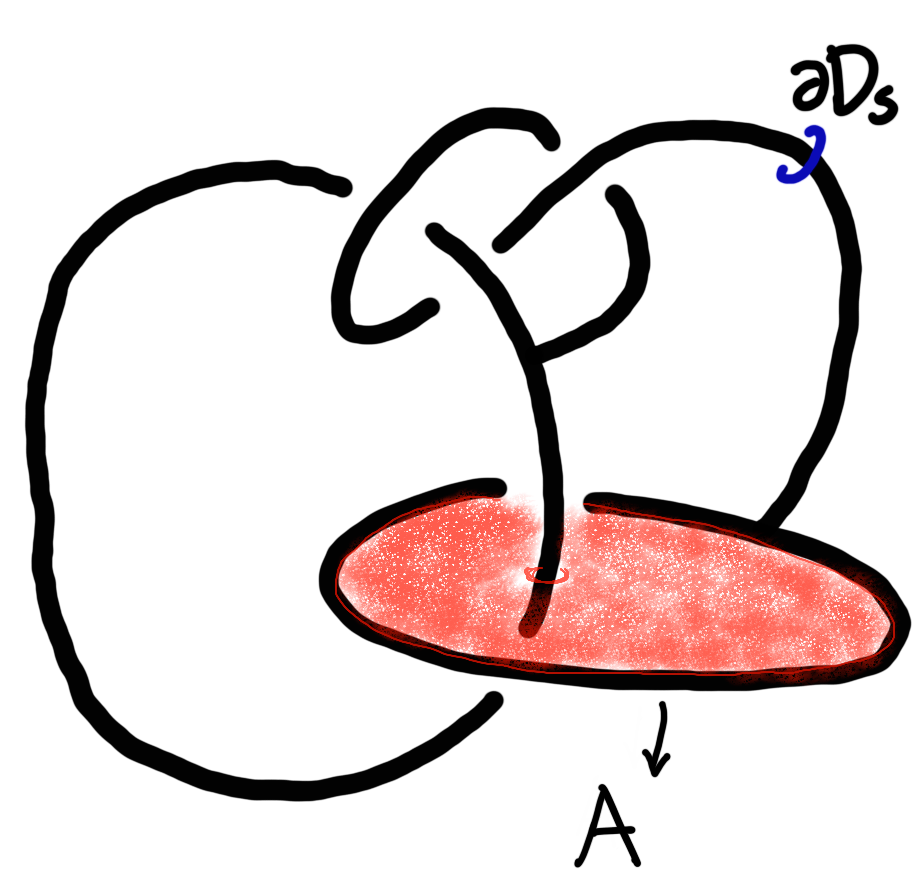}
\caption{Annulus in $\Compl {5_1}$.}
\label{fig:annuli_ext_fiveone}
\end{subfigure}
%
\begin{subfigure}{.45\linewidth}
\center
\includegraphics[scale=.1]{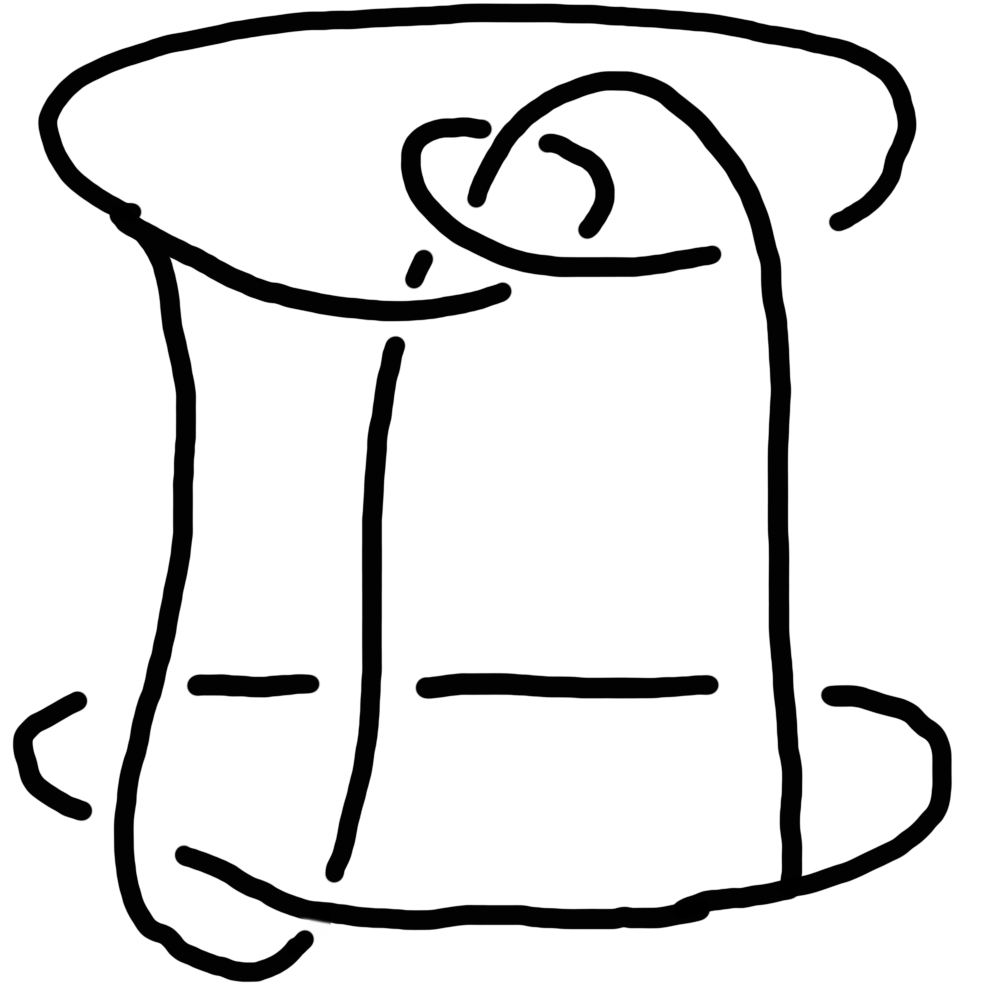}
\caption{$\pairsixone$.}
\label{fig:hksixone}
\end{subfigure}
\begin{subfigure}{.45\linewidth}
\center
\includegraphics[scale=.1]{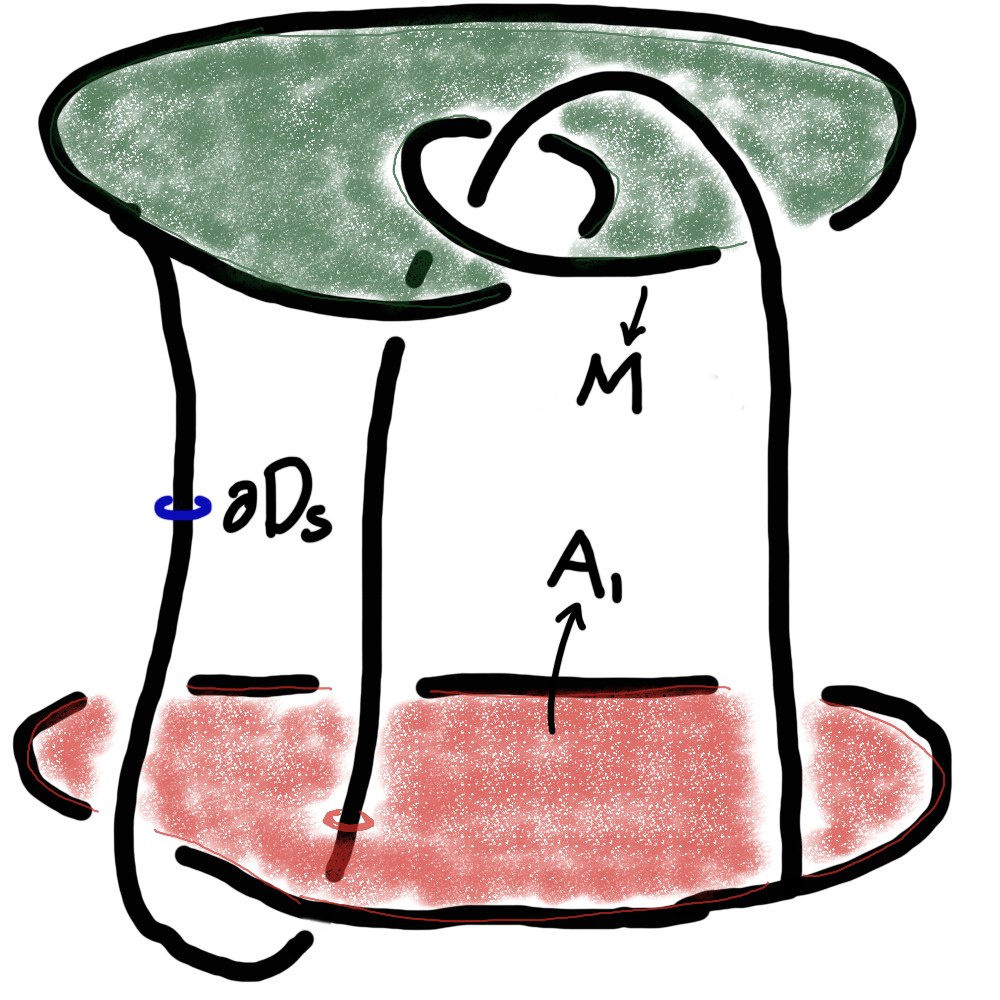}
\caption{Annuli in $\Compl {6_1}$.}
\label{fig:annuli_ext_sixone}
\end{subfigure}
\caption{Handlebody-knots in the knot table \cite{IshKisMorSuz:12}.}
\end{figure}

Consider first the handlebody-knot $\pairfiveone$ in the handlebody-knot table \cite[table $1$]{IshKisMorSuz:12} (Figs.\ \ref{fig:hkfiveone} and \ref{fig:equivalence_fiveone}). It admits a type $2$-$1$ annulus $\annulus$ (see Fig.\ \ref{fig:annuli_ext_fiveone}), which is the unique annulus in $\Compl {5_1}$ 
by \cite[Theorem $1.4$]{Wan:22p}, so its annulus diagram is  
\raisebox{-.2 cm}{\includegraphics[scale=.13]{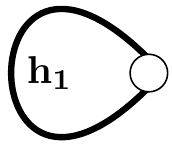}}.

\begin{figure}[t]
\begin{subfigure}{.3\linewidth}
\center
\includegraphics[scale=.13]{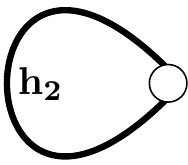}
\end{subfigure}
\begin{subfigure}{.35\linewidth}
\center
\includegraphics[scale=.13]{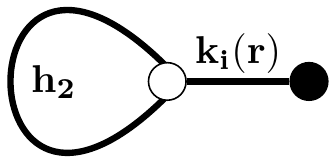}   
\hspace*{.1em}
\raisebox{2.5\height}{{\tiny, $r\in \mathbb{Q}$, $i=1,2$}}
\end{subfigure}
\begin{subfigure}{.3\linewidth}
\center
\includegraphics[scale=.13]{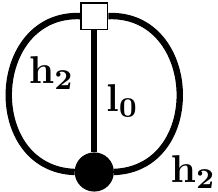}
\raisebox{2.5\height}{{\tiny, $\square=\hnode$ or $\snode$.}}
\end{subfigure}
\caption{Annulus diagrams.}
\label{fig:typetwotwo_ann_diags}
\end{figure}

Secondly, for the handlebody-knot $\pairsixone$ in \cite[Table $1$]{IshKisMorSuz:12} (see Fig.\ \ref{fig:hksixone} and \ref{fig:equivalence_sixone}), we observe that it 
admits a type $2$-$2$ annulus $\annulus_1$, 
so by \cite[Theorem $1.4$]{Wan:22p},
the annulus diagram is one of the diagrams in Fig.\ \ref{fig:typetwotwo_ann_diags}. 
\noindent
Note that $\Compl{6_1}$ also admits 
a M\"obius band $\mobius$ as shown in 
Fig.\ \ref{fig:annuli_ext_sixone}. 
Since there exists a separating disk $\disk_s$
disjoint from $\partial \mobius$, and 
since the core of the component $V$ 
of $\HK-\openrnbhd{\disk_s}$ 
containing $\partial \mobius$ is a trivial knot,
the frontier $\annulus_2$
of a regular neighborhood 
of $\mobius$ is a type $3$-$2$ii
annulus. Its annulus diagram is therefore the second one in Fig.\ \ref{fig:typetwotwo_ann_diags}. 
The boundary slope of $\mobius$ with respect to
$(\sphere,V)$ is $2$, so the annulus diagram of $\pairsixone$ is \raisebox{-.4\height}{\includegraphics[scale=.13]{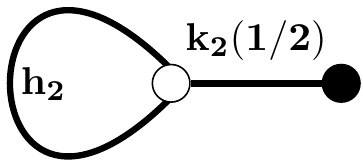}}.

\begin{figure}[t]
\begin{subfigure}{.25\linewidth}
\center
\includegraphics[scale=.08]{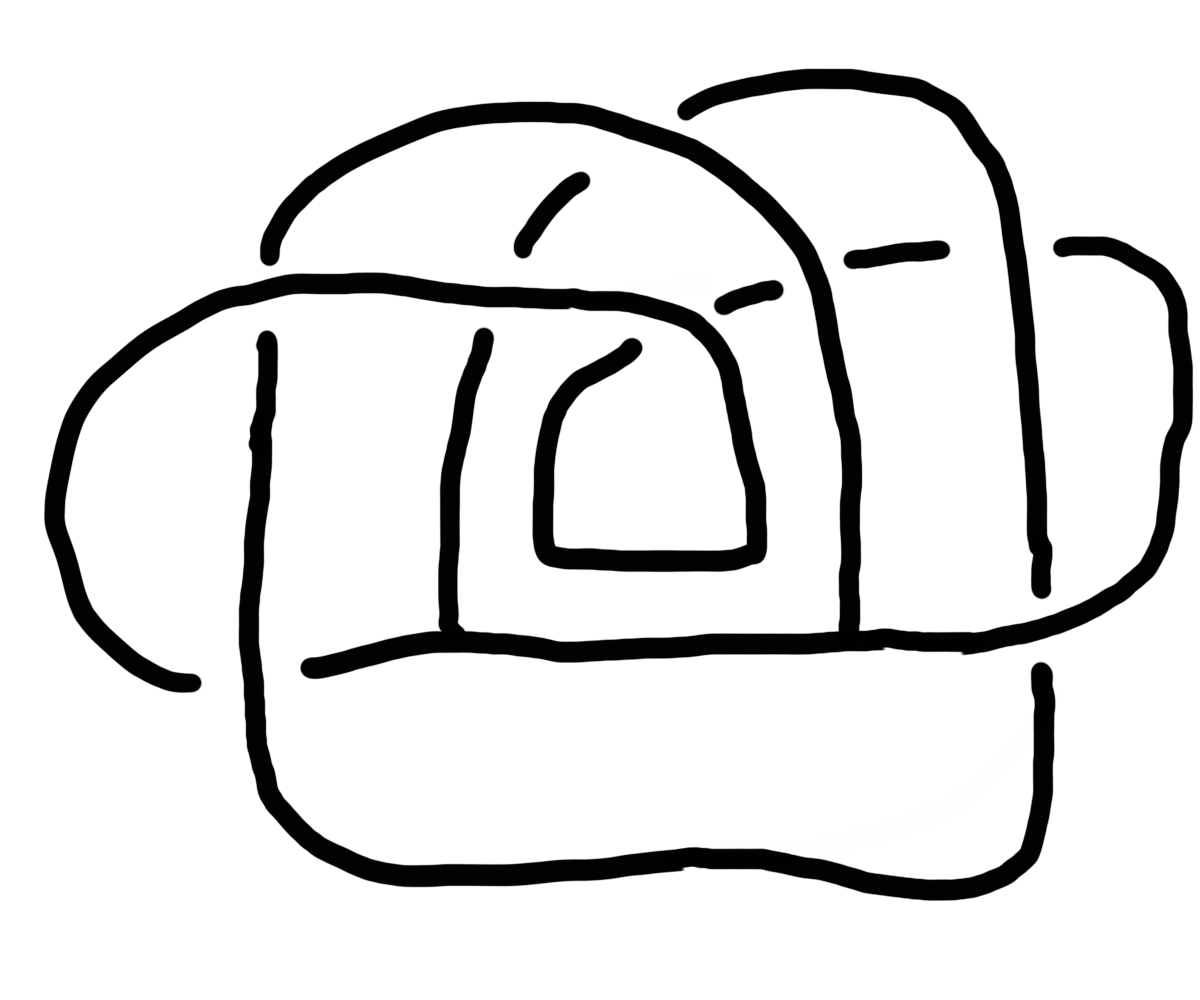}
\caption{$\pairfivetwo$.}
\label{fig:hkfivetwo}
\end{subfigure}
\begin{subfigure}{.35\linewidth}
\center
\includegraphics[scale=.08]{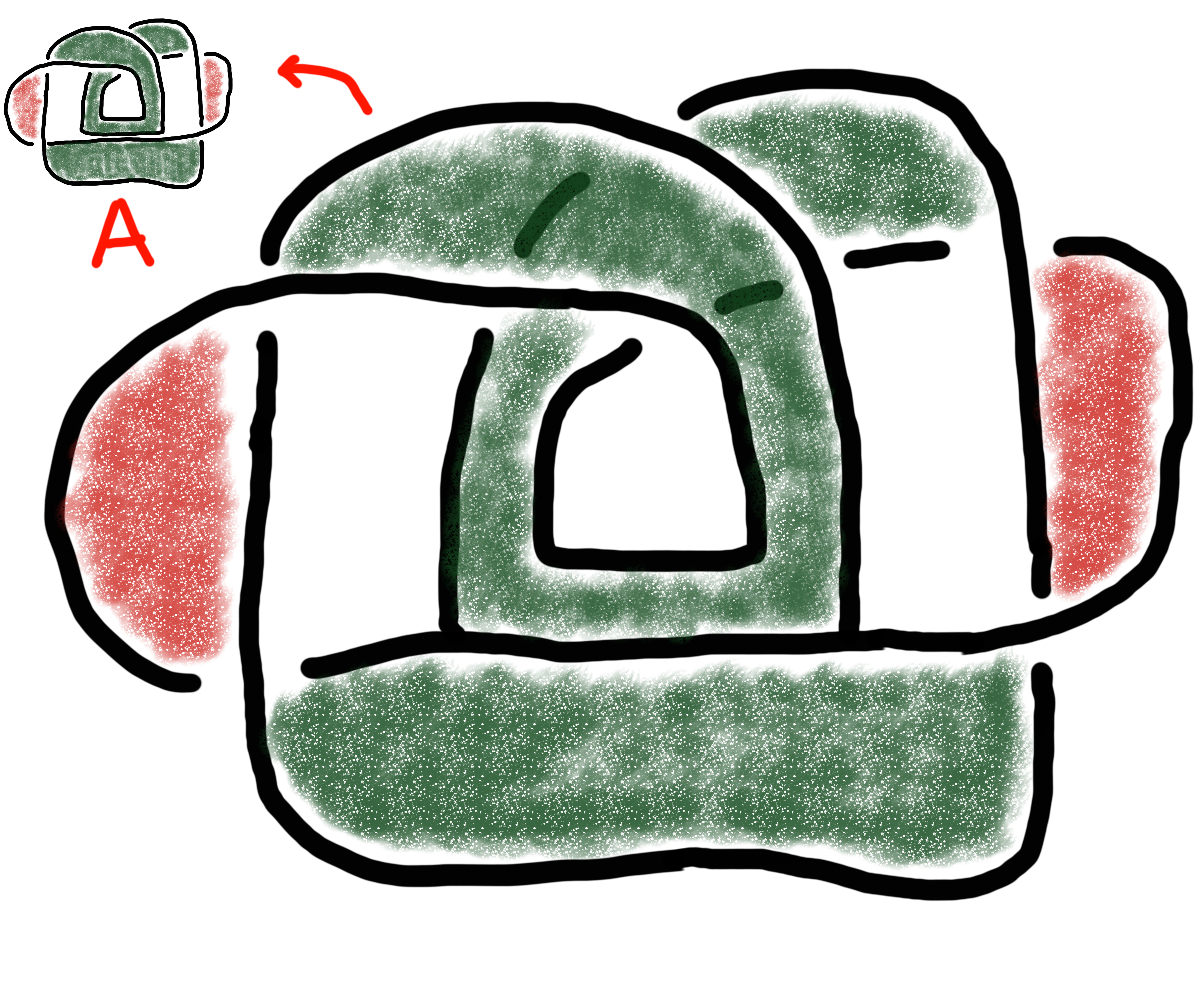}
\caption{Type $3$-$3$ annulus in $\Compl {5_2}$.}
\label{fig:annuli_ext_fivetwo}
\end{subfigure}
\begin{subfigure}{.35\linewidth}
\center
\includegraphics[scale=.08]{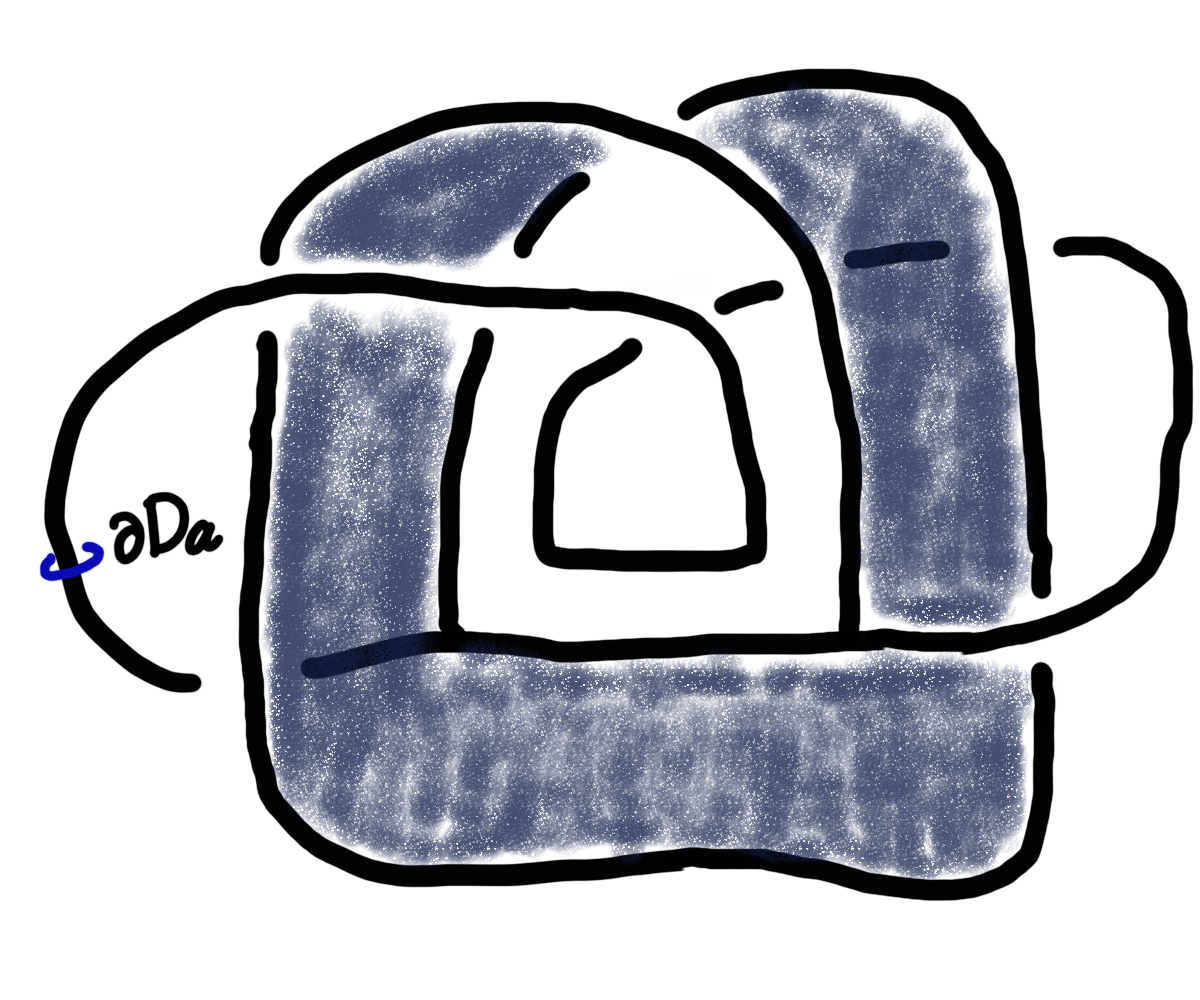}
\caption{M\"obius band in $\Compl {5_2}$.}
\label{fig:Mobius_ext_fivetwo}
\end{subfigure}

\begin{subfigure}{.45\linewidth}
\center
\includegraphics[scale=.1]{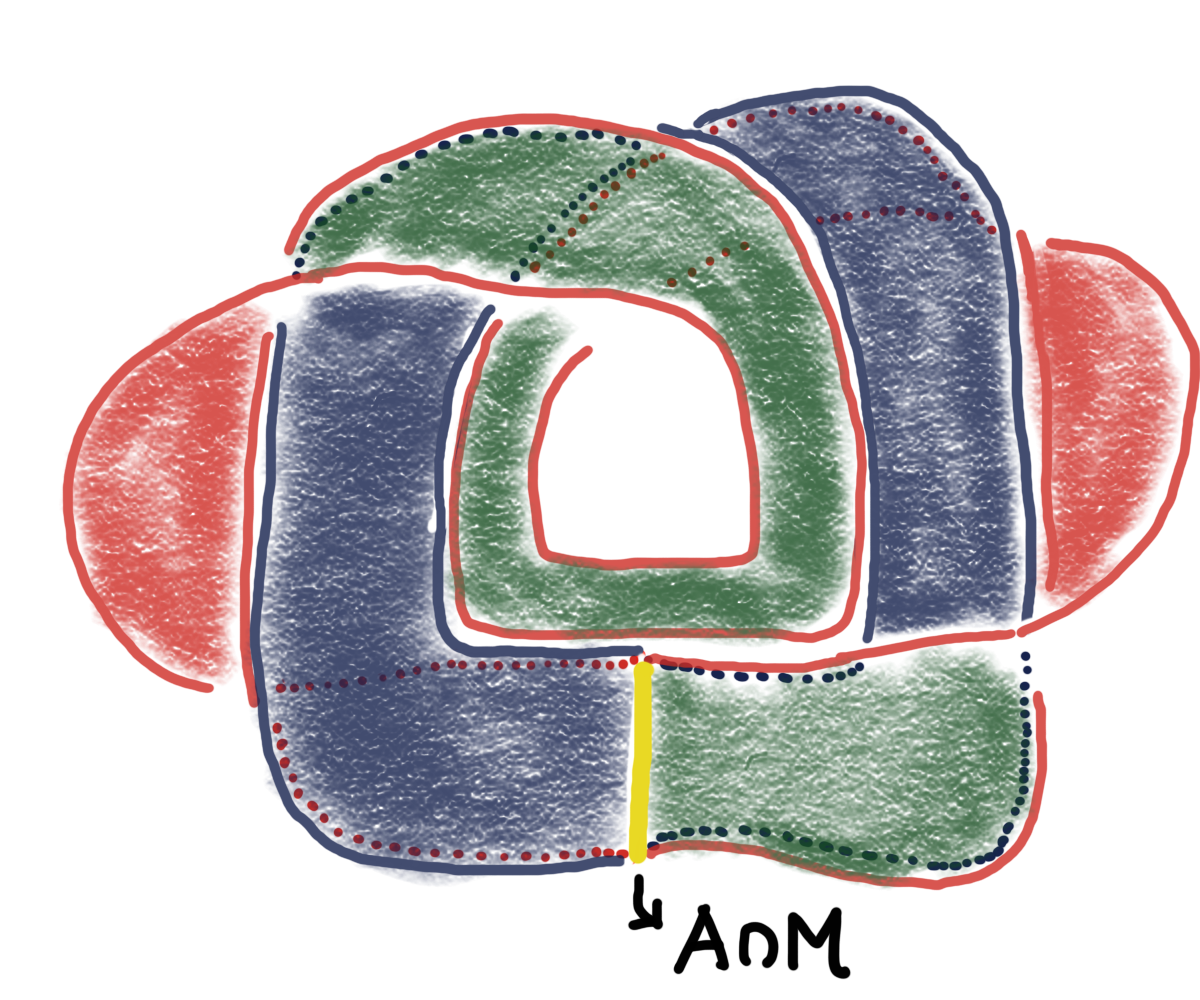}
\caption{$\annulus\cap \mobius$.}
\label{fig:intersection}
\end{subfigure} 
\begin{subfigure}{.45\linewidth}
\center
\includegraphics[scale=.1]{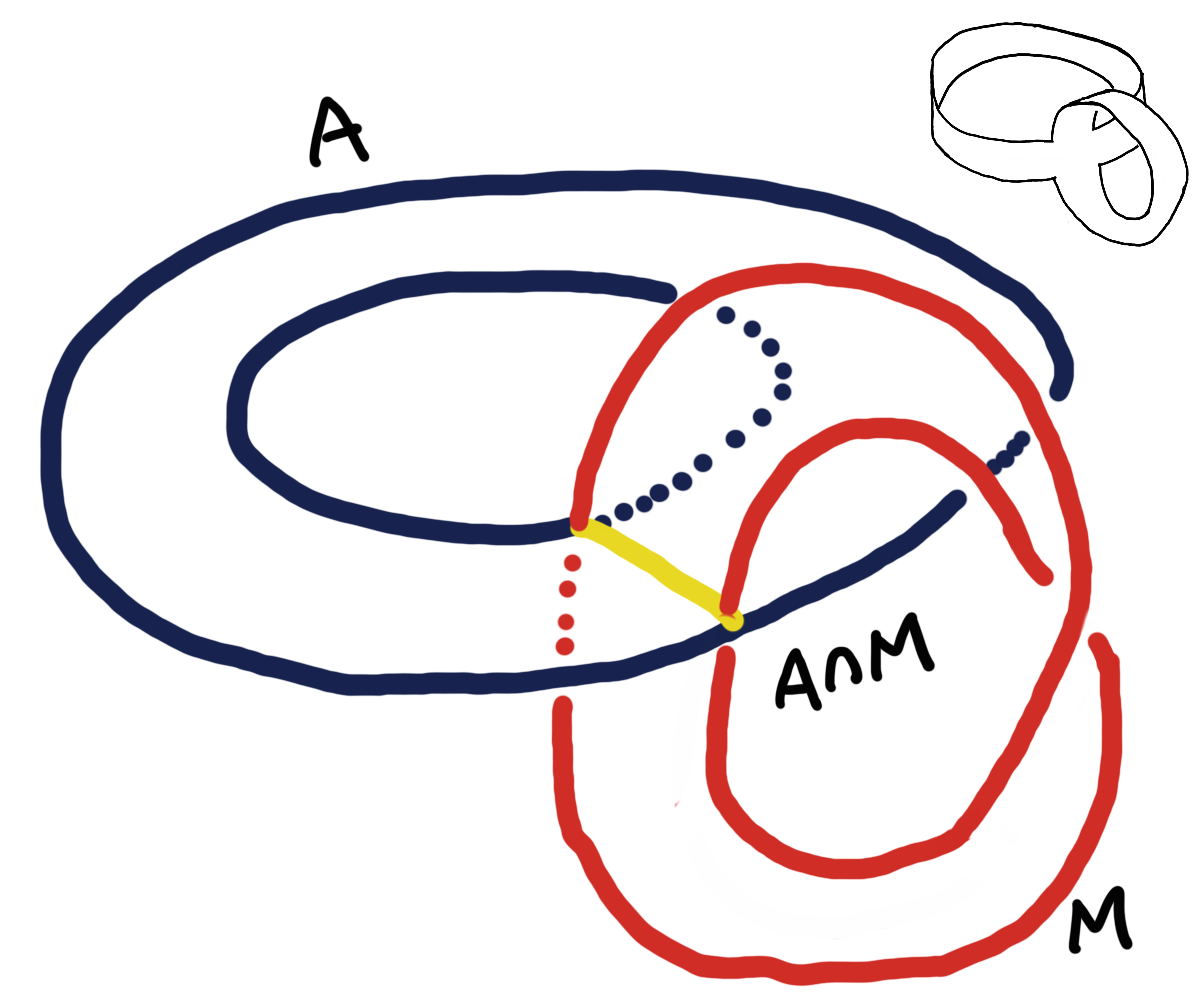}
\caption{$I$-bundle structure.}
\label{fig:bundle}
\end{subfigure} 
\caption{$\pairfivetwo$ and $\annulus,\mobius\subset \Compl{5_2}$.}
\end{figure}

\begin{figure}[H]
\begin{subfigure}{.39\linewidth}
\includegraphics[scale=.085]{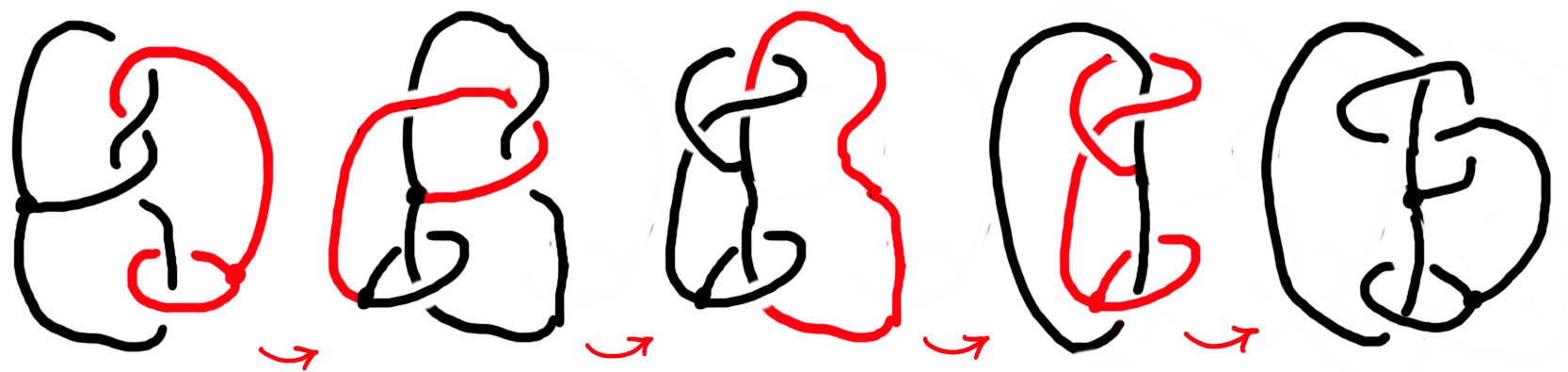}
\caption{\footnotesize $\pairfiveone$.}
\label{fig:equivalence_fiveone}
\end{subfigure}
\begin{subfigure}{.31\linewidth}
\centering
\includegraphics[scale=.08]{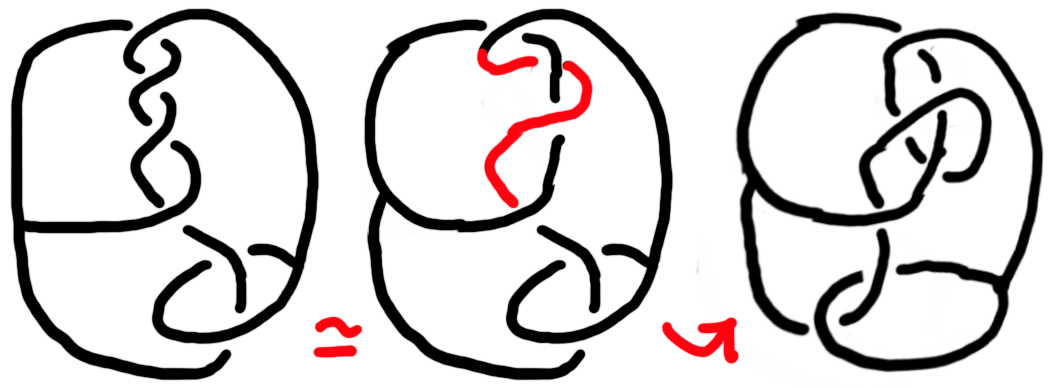}
\caption{\footnotesize $\pairsixone$.}
\label{fig:equivalence_sixone}
\end{subfigure}
\begin{subfigure}{.2\linewidth}
\centering
\includegraphics[scale=.076]{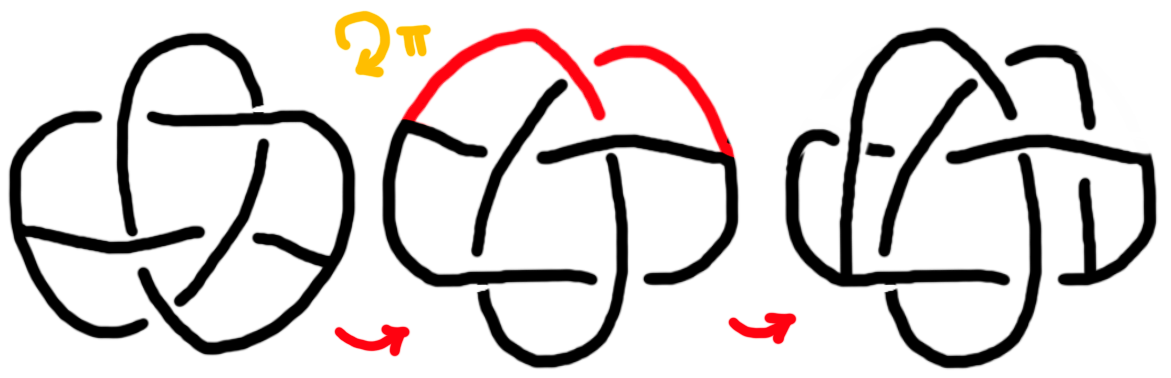}
\caption{\footnotesize $\pairfivetwo$.}
\label{fig:equivalence_fivetwo}
\end{subfigure}
\caption{Knot diagrams in Figs.\ \ref{fig:hkfiveone}, \ref{fig:hksixone}, \ref{fig:hkfivetwo} versus those in \cite{IshKisMorSuz:12}.}
\label{fig:equivalences}
\end{figure}

Now, we compute the annulus diagram of 
the handlebody-knot $\pairfivetwo$ in \cite{IshKisMorSuz:12} (Figs.\ \ref{fig:hkfivetwo} and \ref{fig:equivalence_fivetwo}). 
Observe that it admits a type $3$-$3$ annulus $\annulus$
and a M\"obius band $\mobius$ as shown in Figs.\ \ref{fig:annuli_ext_fivetwo} and \ref{fig:Mobius_ext_fivetwo}.
The frontier $\annulus_m$ 
of a regular neighborhood of $\mobius$ 
in $\Compl\HK$ is a type $3$-$2$ annulus since there exists a non-separating disk $\disk\subset 5_2$ 
disjoint from $\partial \mobius$. Because 
the core of the solid torus $V:=\HK-\openrnbhd{\disk}$ is 
a trivial knot in $\sphere$, $\annulus_m$ is of type $3$-$2$ii.

Note that $\annulus$ and $\mobius$ meets at an arc (Fig.\ \ref{fig:intersection}). 
Let $U$ be a regular neighborhood of 
$\annulus\cup \mobius$ in $\Compl{5_2}$. Then 
there is an admissible I-bundle structure $U\rightarrow B$, where $B$ is a Klein bottle
with one disk removed (Fig.\ \ref{fig:bundle}). In particular, this implies the characteristic diagram
of $\Compl{5_2}$ is 
\raisebox{-.07\height}{\includegraphics[scale=.17]{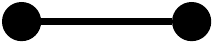}}, where the characteristic annulus $\annulus_c$ 
is the frontier of $U$ in $\Compl {5_2}$. 
To determine the type of $\annulus_c$, 
we need a few lemmas.

\begin{lemma}\label{lm:typethreetwoii_generator}
Let $\annulus$ be a type $3$-$2$ii essential annulus in the exterior of $\pair$, and $U\subset\Compl\HK$ is the $3$-manifold bounded by $\annulus$ and the non-annular component in $\partial\HK$ cut off by $\partial\annulus$.
Then the image of the induced homomorphism
$H_1(\annulus)\rightarrow H_1(U)$ is a generator.
\end{lemma}
\begin{proof}
Let $\disk\subset \HK$ be the unique non-separating disk disjoint from $\partial \annulus$; set 
$W:=\HK-\openrnbhd{\disk}$.   
Since $\annulus$ is $\partial$-compressible in 
$\Compl W$, it cobounds a solid torus $U'\subset \Compl W$ with an annulus in $\partial W$
such that 
the induced homomorphism
$H_1(\annulus_c)\rightarrow H_1(U')$ 
is an isomorphism. In particular, $U$ is obtained by 
digging a tunnel through $U'$, 
more precisely, $U'=U\cup\rnbhd{\disk}$. 
The claim then follows from the short exact sequence:
\[0\rightarrow H_2(U',U)\rightarrow H_1(U)\rightarrow H_1(U')\rightarrow 0.\]
\end{proof}

\begin{lemma}\label{lm:typefourone_exclusion}
If 
$\pair$ 
admits a type $4$-$1$ essential annulus $\annulus$, then
every essential annulus disjoint from $\annulus$ is separating.
\end{lemma}
\begin{proof}
Let $V$ be the solid torus cut off by $\annulus$ 
from $\Compl\HK$, 
and $W$ the closure of $\Compl\HK-V$. 
Denote by 
$\annulus'$ the annulus cut off by $\partial \annulus$ 
from $\partial \HK$. Then by the definition of a type $4$-$1$ annulus, the closure $\torus$ 
of $\partial\HK-\annulus'$, 
a torus with two open disks removed, 
is essential in $\Compl V$. 

Now, since no essential disk in $\HK$ disjoint from 
$\annulus$ exists, $\annulus'$ is essential in $\Compl W$.
By the essentiality of $\annulus'$, 
if there exists a non-trivial $\partial$-compressing disk of $\Compl W$, then there exists a non-trivial $\partial$-compressing disk of $\Compl W$ disjoint from $A'$, 
but this contradicts $\torus$ being essential in $\Compl V$. 
Therefore, $\Compl W$ is $\partial$-irreducible. 
This implies $W$ is $\partial$-reducible, so 
the frontier of the compression body of $W$ is empty, a torus or two tori.
The latter two cases are excluded since $\pair$ is atoroidal. $W$ is thus a handlebody of genus two. Applying \cite[Theorem $1.4$]{Prz:83}, 
$\torus$ induces an incompressible torus
$\hat{\torus}$ in the $3$-manifold $\hat{M}$ 
obtained by performing 
Dehn surgery on $V$ along the boundary of 
$\annulus$. Since $\annulus$ is essential in $\Compl\HK$, the boundary slope of $A$ with respect to $(\sphere,V)$ is non-integral. Hence,
by \cite[Lemma $3.14$]{KodOzaGor:15}, 
the core of $V$ is a hyperbolic knot, and therefore an Eudave-Mu\~noz knot by \cite{GorLue:04}. 
In particular, 
there exists an incompressible torus $\torus_0$
that separates $\hat M$ into two Seifert fiber spaces over the disk with two exceptional fibers \cite{Eud:97}, \cite{Eud:20}. 
Note that $\hat M$ 
is itself not Seifert-fibered by \cite{BoyZha:98},
being obtained by an non-integral Dehn surgery on $(\sphere,V)$. Isotope $\torus_0$ so that the number of 
components in 
$\hat \torus \cap \torus_0$ is minimized. If $\hat\torus\cap\torus_0\neq\emptyset$, then 
the closures $\annulus_1,\annulus_2$ of two neighboring
components in $\hat \torus- \torus_0$ are essential annuli
in $M_1,M_2$, respectively. By the vertical-horizontal theorem \cite[Corollary $5.7$]{Joh:79},
one can isotope the fibration on $M_1,M_2$ so that
$\annulus_1,\annulus_2$ are vertical, 
but this implies $\hat M$ 
is Seifert fibered, a contradiction. 
Thus $\hat \torus\cap \torus_0=\emptyset$, wherefrom
one deduces that $\hat \torus, \torus_0$ are isotopic. 
It may thus be assumed that $\hat\torus=\torus_0$; 
let $M_1$ be the component containing $W$.

Suppose $\Compl\HK$ admits a non-separating essential annulus $\bar \annulus$ disjoint from $\annulus$. 
Then $\bar \annulus$ is non-separating in $M_1$ as well,
and furthermore, it is essential in $M_1$ 
by the $\partial$-irreducibility of $M_1$, 
but this contradicts the fact that no essential non-separating annulus exists in a Seifert fiber space over a disk with two exceptional fibers.    
\end{proof}

\begin{lemma}\label{lm:cano_ann_stick_type}
If the characteristic diagram of 
the exterior $\Compl\HK$ of $\pair$
is \raisebox{-.07\height}{\includegraphics[scale=.17]{char_diag_stick}},
then its annulus diagram is \raisebox{-.07\height}{\includegraphics[scale=.17]{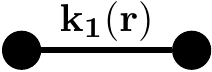}},
for some non-integral $r\in\mathbb{Q}$.
\end{lemma}
\begin{proof}
Let $\annulus$ be the characteristic annulus 
of $\Compl\HK$ and $V\subset\Compl\HK$ 
the solid torus cut off by $\annulus$. Set 
$U:=\Compl\HK-\mathring{V}$, and by \cite[Proposition $2.12v$]{Wan:22p}, there exists an I-bundle structure $\pi:U\rightarrow B$ with $B$ a Klein bottle $B$ with an open disk removed. 
By Lemma \ref{lm:typefourone_exclusion}, 
$\annulus$ cannot be of type $4$-$1$ 
since $U$ and hence $\Compl\HK$ admit a non-separating essential annulus.

Since $B$ is a Klein bottle with one open disk removed,
there exist two simple loops $l,l'\subset B$ with 
$l\cap l'$ a point, $\pi^{-1}(l)$ an annulus $\annulus_l\subset\Compl\HK$ and $\pi^{-1}(l')$ a M\"obius band $\mobius_l'\subset\Compl\HK$ such that $U$ is a regular neighborhood of $\annulus_l\cap \mobius_l'$. In particular, the homology classes $u_a,u_m$
of the cores of $\annulus_l, \mobius_l'$, respectively, generate $H_1(U)$.
$\annulus$ being the frontier of 
$U$ in $\Compl\HK$ implies that the image of 
a generator of $H_1(\annulus)$ under the homomorphism $H_1(\annulus)\rightarrow H_1(U)$ is $2u_a$, so 
$\annulus$ has to be of type $3$-$2$i by Lemma \ref{lm:typethreetwoii_generator}. 
\end{proof} 

\begin{figure}[h]
\begin{subfigure}{.45\linewidth}
\center
\includegraphics[scale=.1]{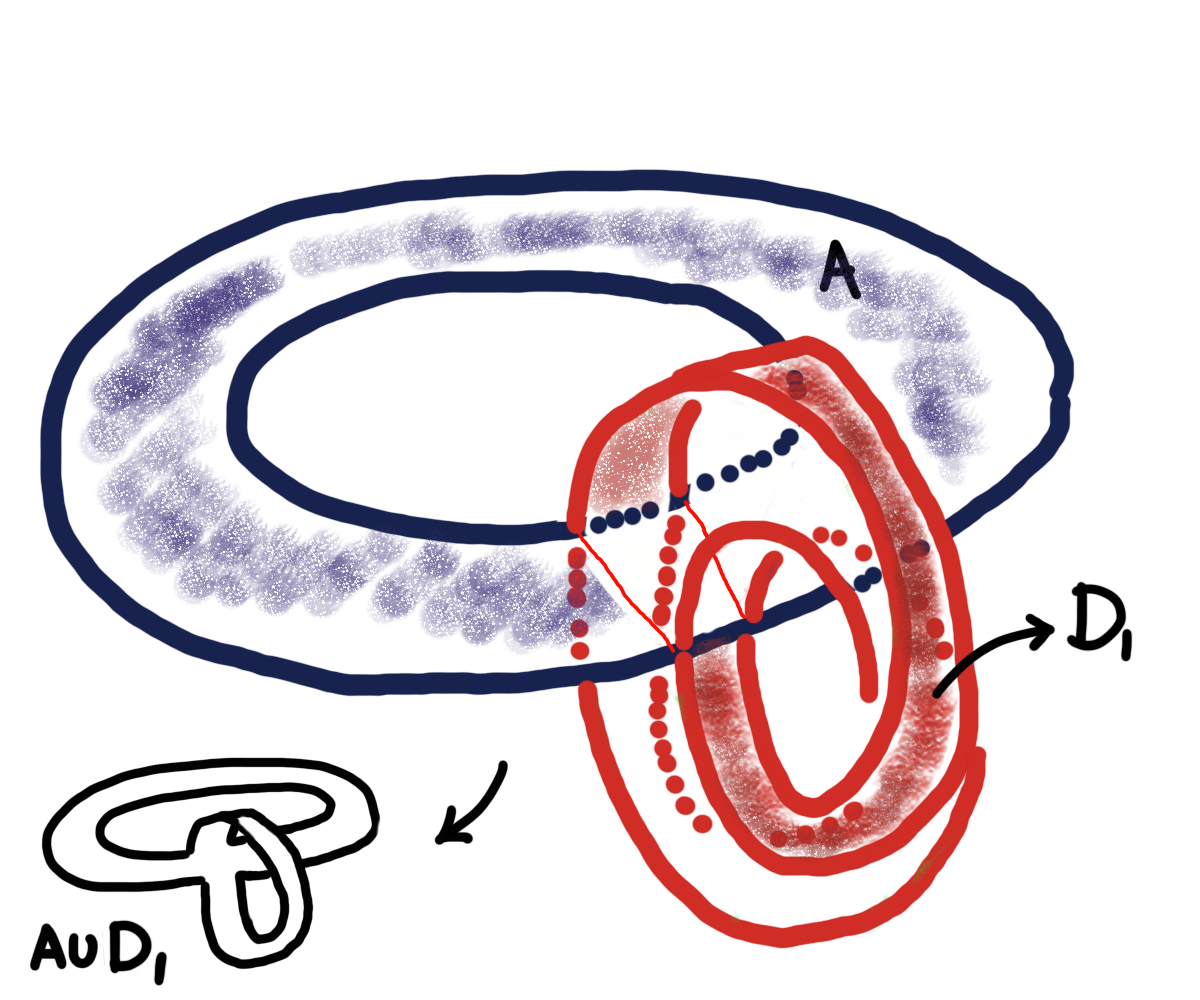}
\caption{$\annulus$ and $\disk_1\subset \annulus^\flat$.}
\label{fig:boundary}
\end{subfigure} 
\begin{subfigure}{.45\linewidth}
\center
\includegraphics[scale=.1]{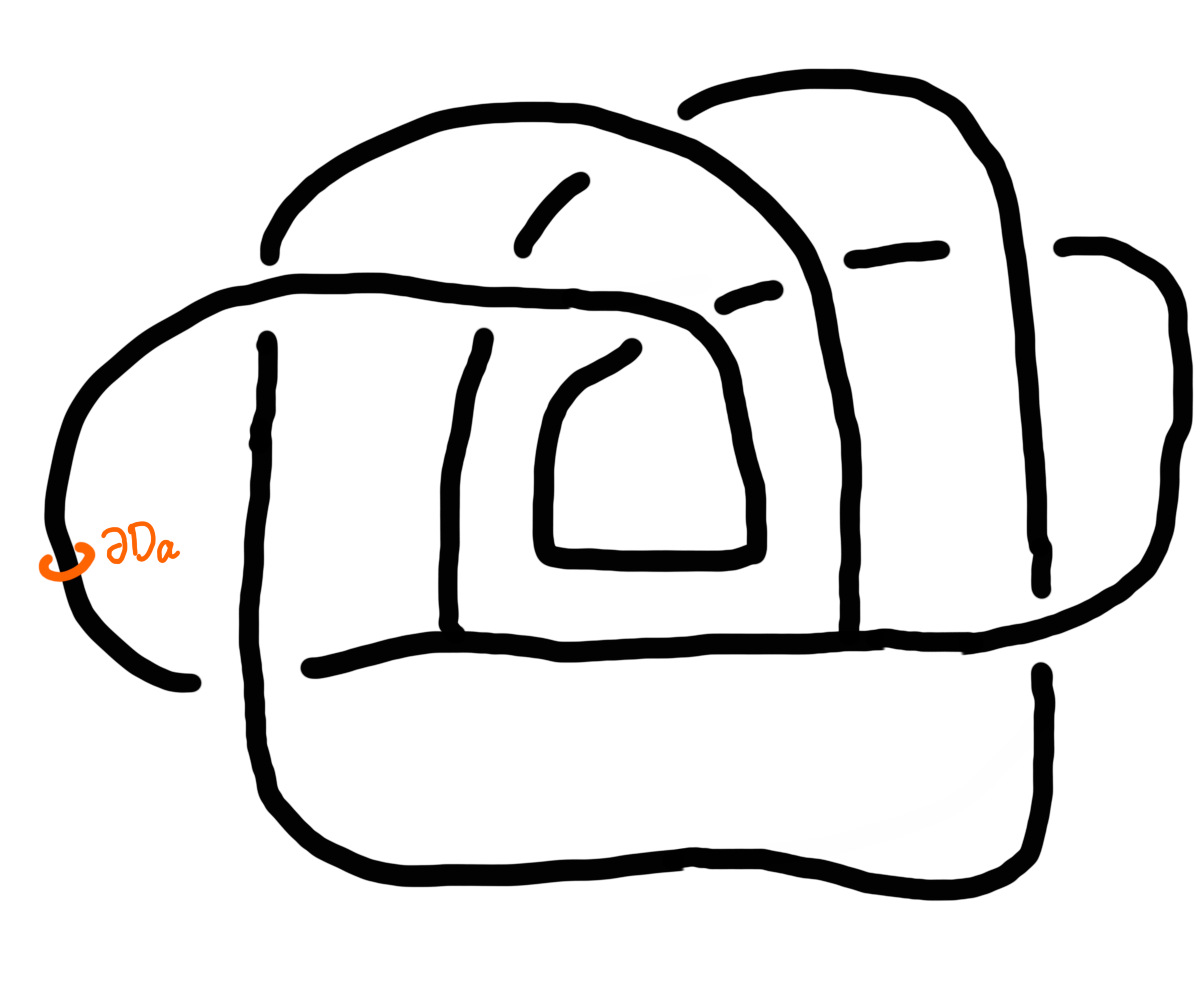}
\caption{Meridian disk $\disk_a$.}
\label{fig:meridian}
\end{subfigure}
\caption{$\partial \annulus_c\parallel \partial(\annulus\cup\disk_1)$ and $\disk_a$.} 
\end{figure}

Return to the example $\pairfivetwo$; by Lemma \ref{lm:cano_ann_stick_type}, the characteristic annulus
$\annulus_c\subset \Compl{5_2}$ is of type $3$-$2$i. 
To determine
the slope of $\annulus_c$, we recall 
that, if $\disk_c\subset 5_2$
is a non-separating disk disjoint from $\partial \annulus_c$, 
then $\annulus_c$ is the cabling annulus of 
the solid torus 
$W:=5_2-\openrnbhd{\disk_c}$ in $\sphere$. 
Thus to compute the slope of $\annulus_c$, 
it amounts to identifying the knot type of $(\sphere, W)$. 
To identify the knot type of $(\sphere, W)$, 
we search for an essential separating disk $\disk_s\subset 5_2$ disjoint from $\partial \annulus_c$ 
since such a $\disk_s$ cuts $5_2$ 
into two pieces with one being isotopic to $W$ in $\sphere$. 
Let $\annulus^\flat$ be 
the annulus $\rnbhd{\mobius}\cap 5_2$.
Note that
$\partial \annulus$ meets $\annulus^\flat$ at two arcs, which cut $\annulus^\flat$ 
into two disks $\disk_1,\disk_2$ (see Fig.\ \ref{fig:boundary}).  
Since $\annulus_c$ is 
the frontier of a regular neighborhood of 
$\annulus\cup \mobius$, the loop $l:=\partial(\annulus\cup \disk_1)$ is parallel to components of $\partial \annulus_c$. 

Consider now the disk
$\disk_a$ bounded by the meridian indicated in Fig.\ \ref{fig:meridian}. Then since 
$\disk_a\cap\partial \annulus$ is a point and
$\disk_a\cap \partial \mobius=\emptyset$, 
$\disk_a\cap l$ is a point. Therefore, 
the frontier $\disk_s$ of 
a regular neighborhood of $\disk_a\cup l\subset 5_2$
is a separating disk
disjoint from $\partial \annulus_c$.

\begin{figure}[h]
\begin{subfigure}{.45\linewidth}
\center
\includegraphics[scale=.1]{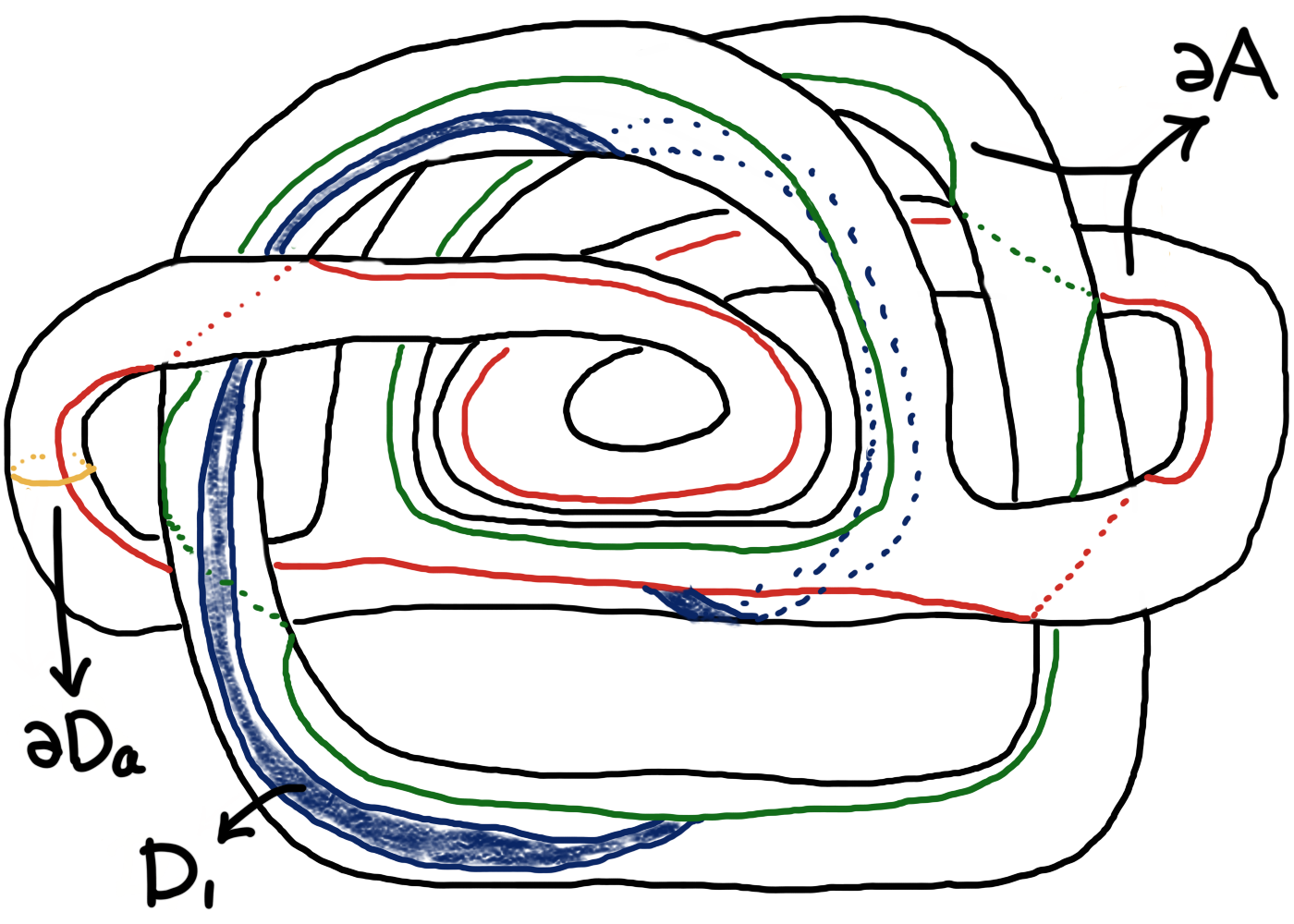}
\caption{$\partial \annulus, \disk_1$ in $5_2$.}
\label{fig:boundary_A_M1prime}
\end{subfigure} 
\begin{subfigure}{.45\linewidth}
\center
\includegraphics[scale=.1]{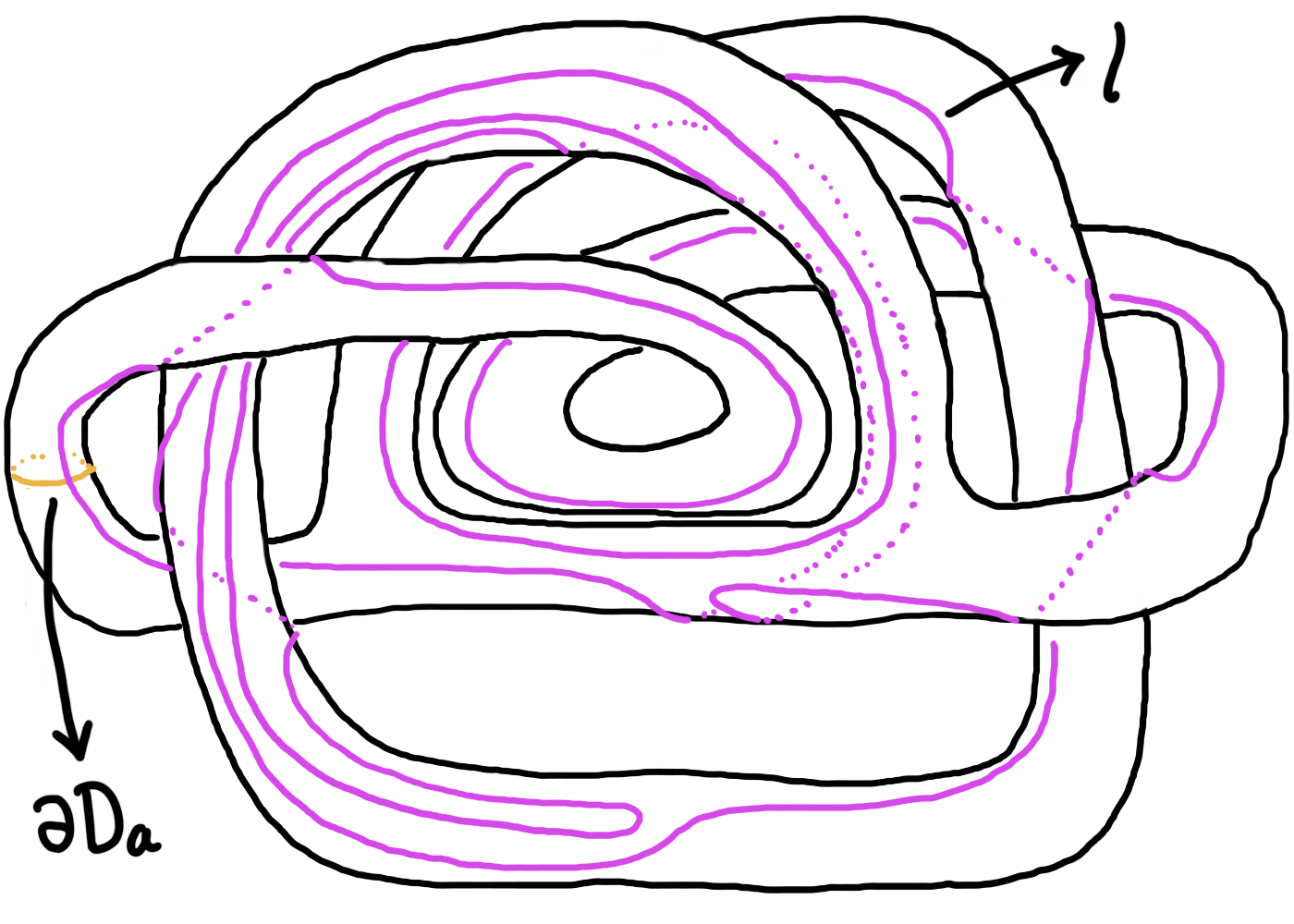}
\caption{$l$ in $\partial 5_2$.}
\label{fig:boundary_Ac}
\end{subfigure}
\begin{subfigure}{.45\linewidth}
\center
\includegraphics[scale=.1]{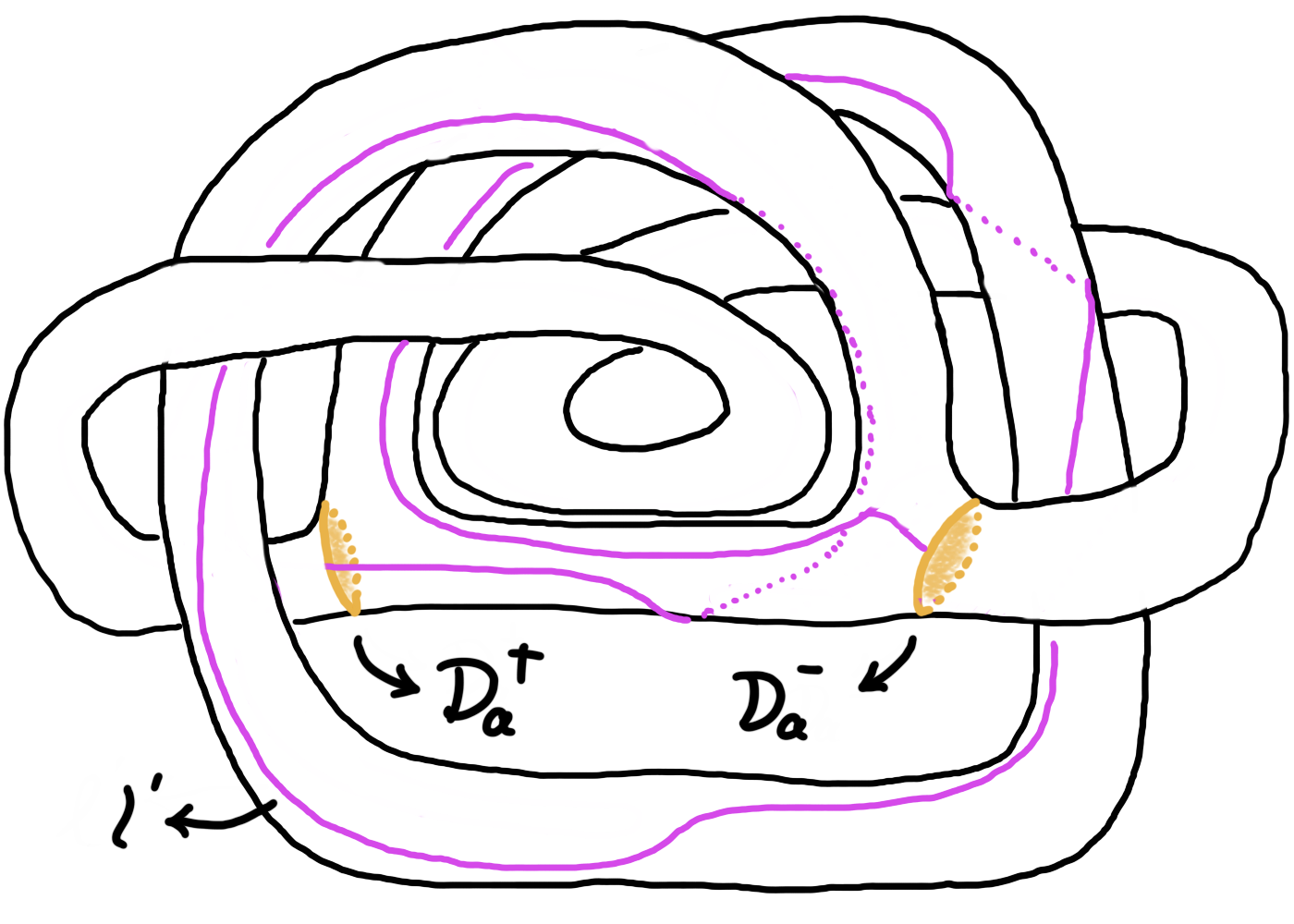}
\caption{$\disk_a^\pm$ and $l'$.}
\label{fig:DDa_lprime}
\end{subfigure} 
\begin{subfigure}{.45\linewidth}
\center
\includegraphics[scale=.1]{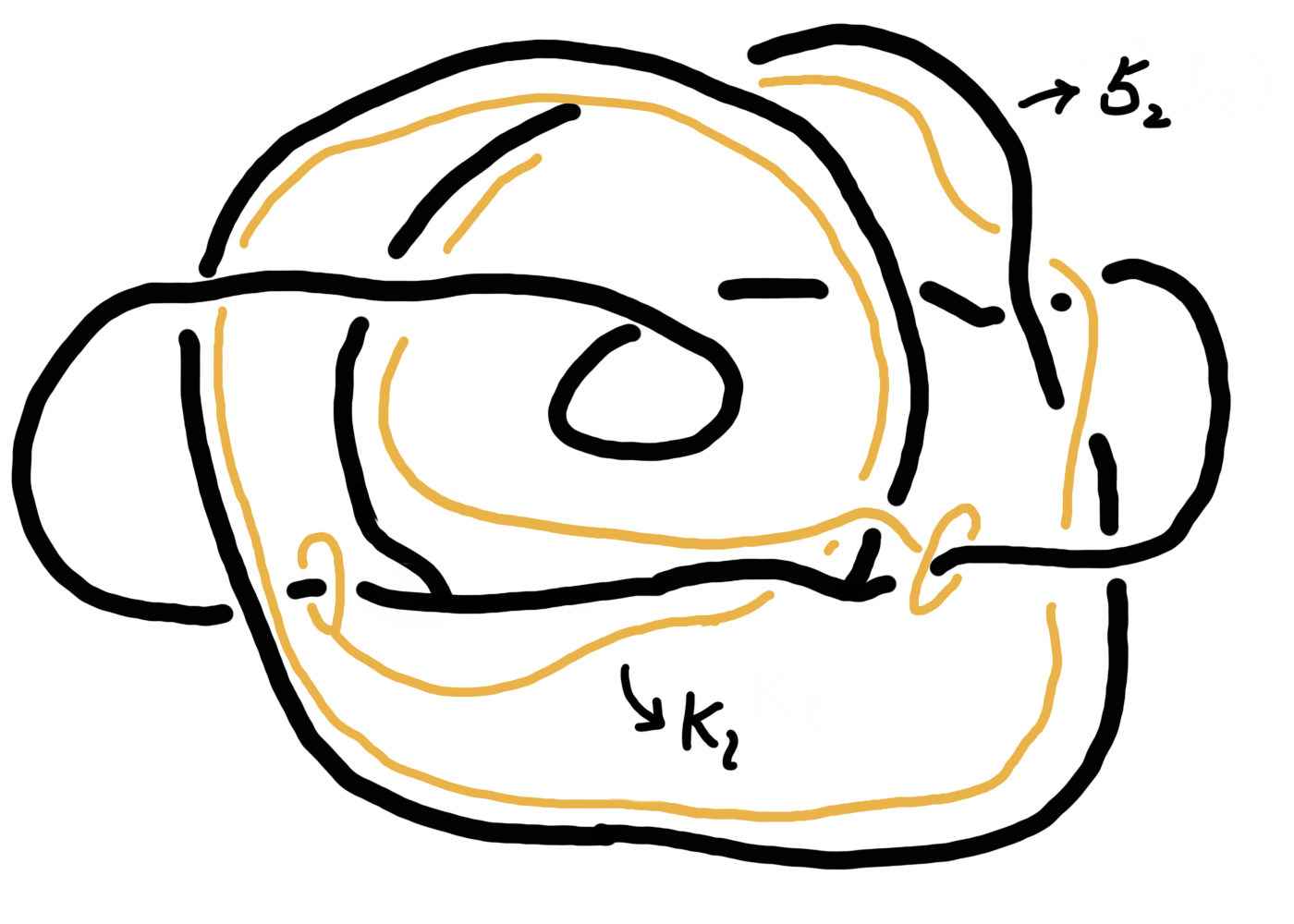}
\caption{$K_l$ in $\pairfivetwo$.}
\label{fig:Kl_fivetwo}
\end{subfigure} 
\caption{Identify $\partial \annulus,\disk_1, l$ and $K_l$ in $\partial 5_2$.}
\end{figure}

\begin{figure}[t]
\begin{subfigure}{.32\linewidth}
\center
\includegraphics[scale=.07]{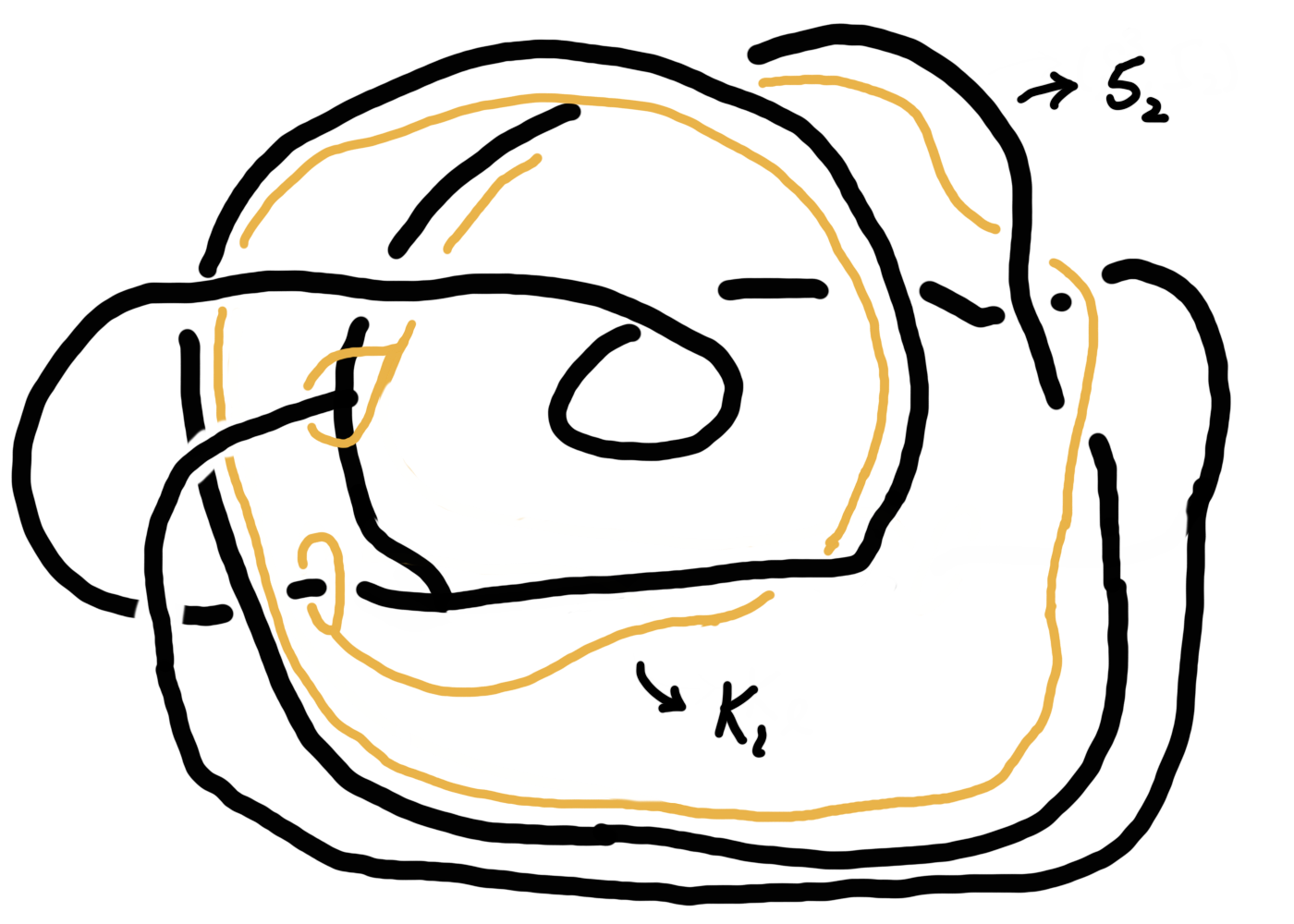}
\caption{Deform $\pairfivetwo$---I.}
\label{fig:deformation1}
\end{subfigure} 
\begin{subfigure}{.32\linewidth}
\center
\includegraphics[scale=.07]{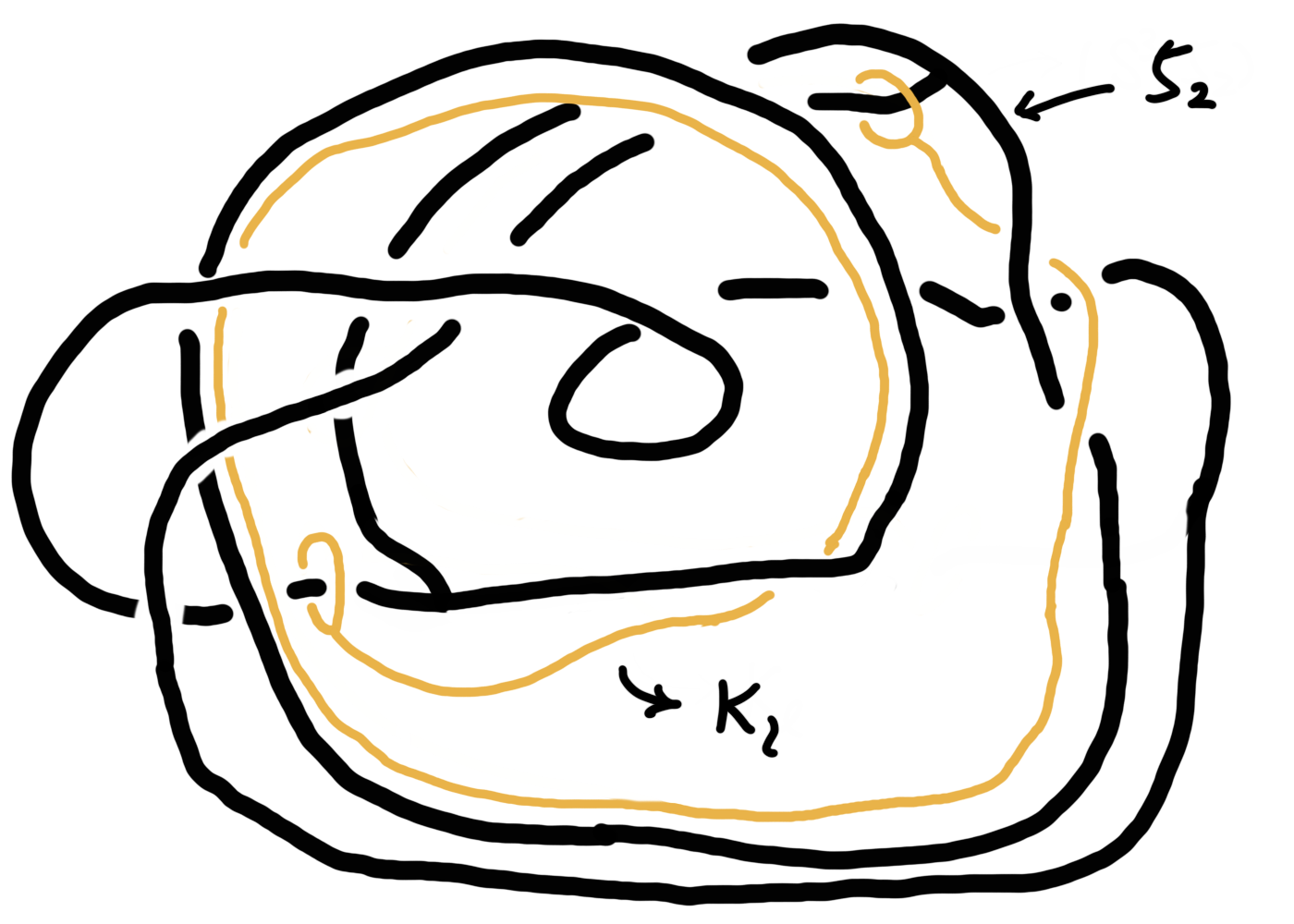}
\caption{Deform $\pairfivetwo$---II.}
\label{fig:deformation2}
\end{subfigure}
\begin{subfigure}{.32\linewidth}
\center
\includegraphics[scale=.07]{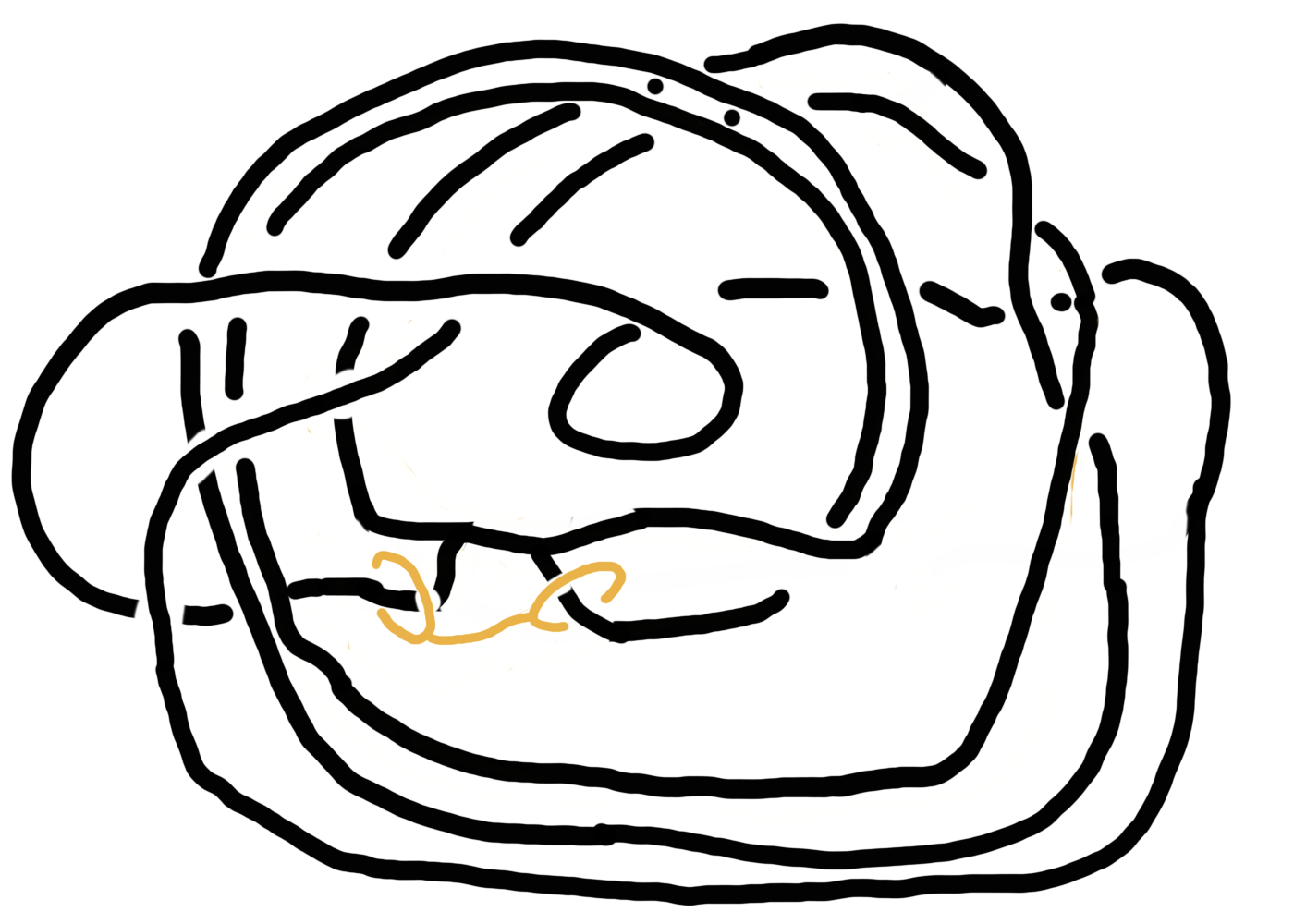}
\caption{Deform $\pairfivetwo$---III.}
\label{fig:deformation3}
\end{subfigure} 
\begin{subfigure}{.32\linewidth}
\center
\includegraphics[scale=.07]{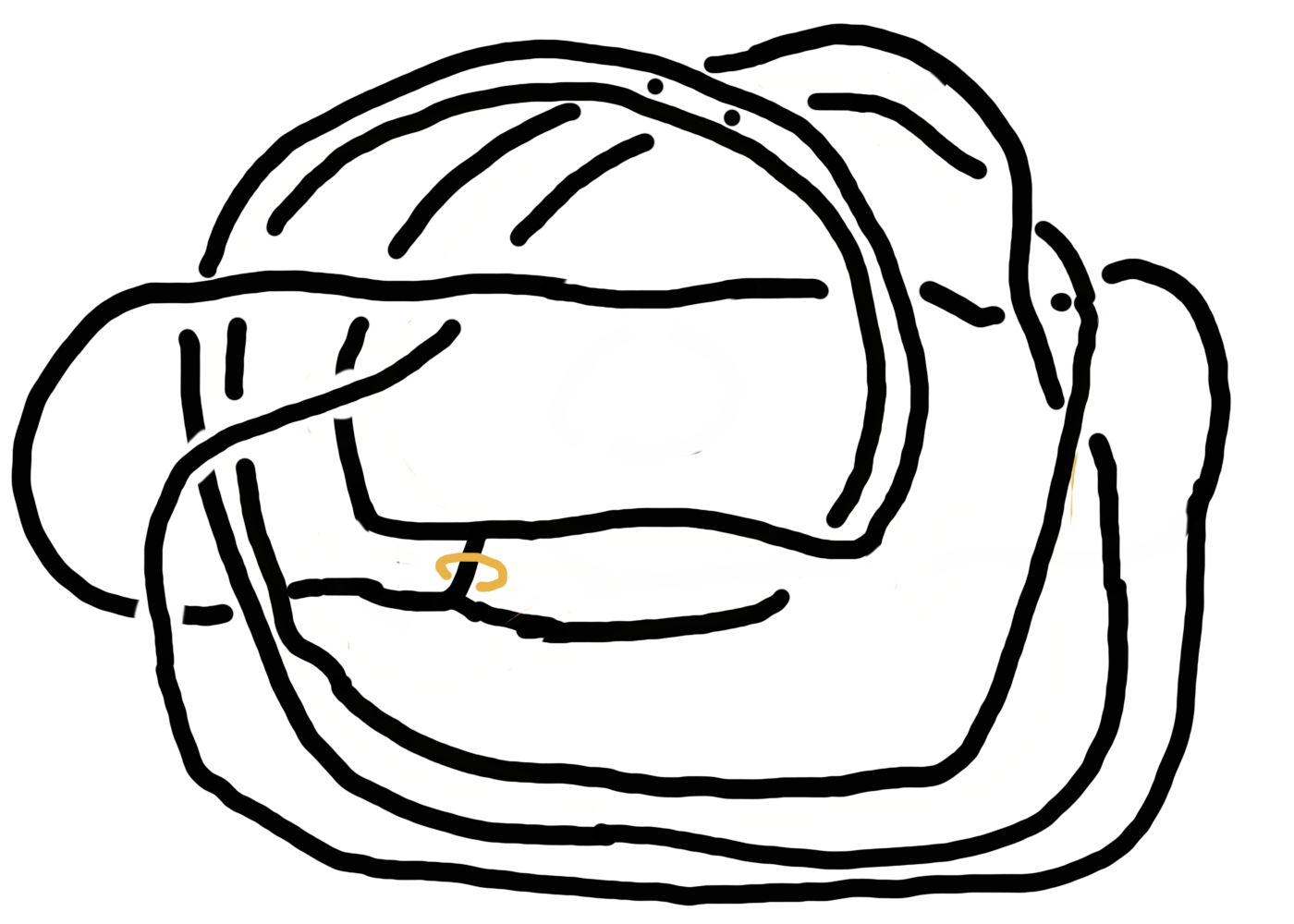}
\caption{Disk $\disk_s$.}
\label{fig:disk_Ds}
\end{subfigure} 
\begin{subfigure}{.32\linewidth}
\center
\includegraphics[scale=.07]{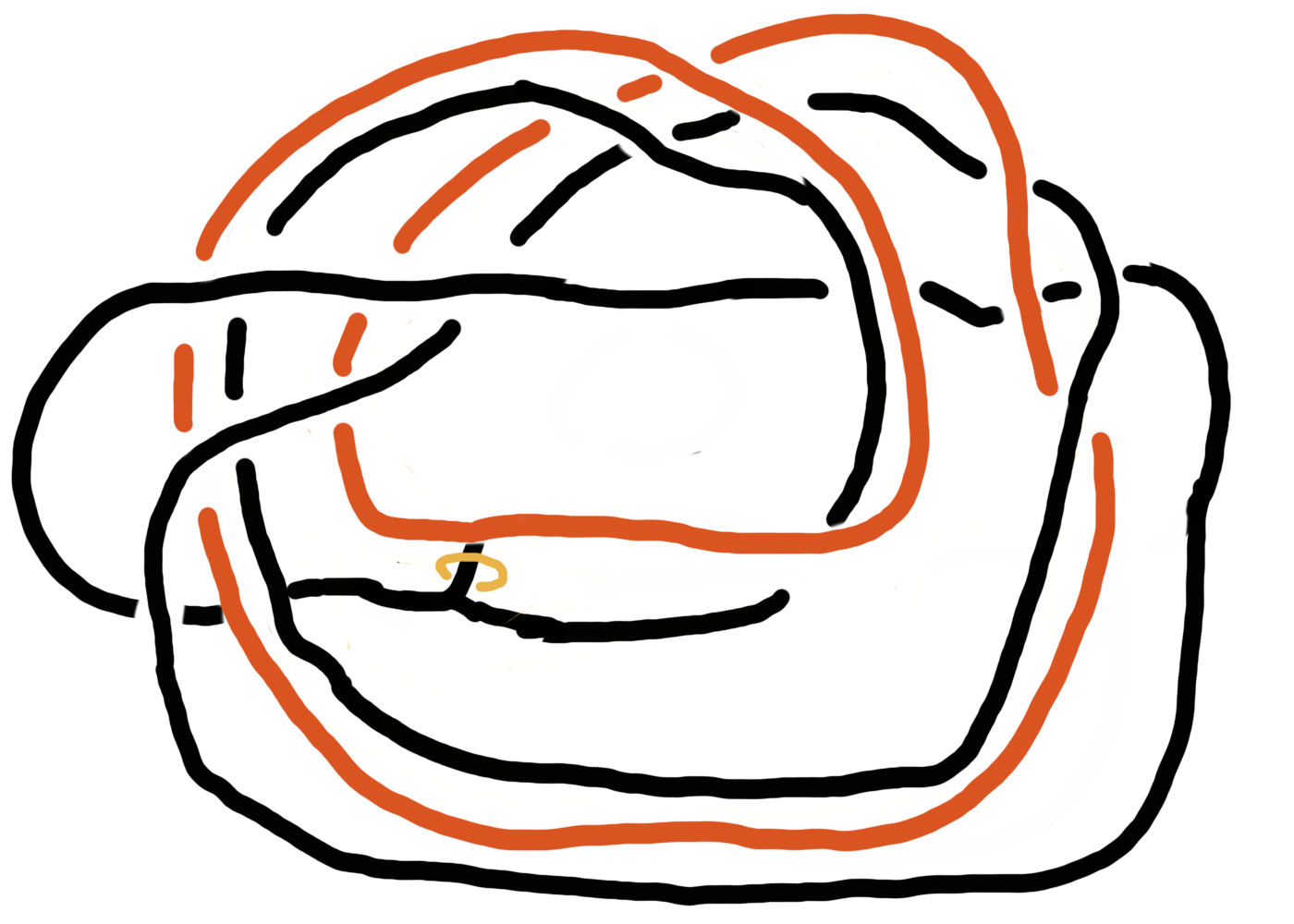}
\caption{Torus knot with an arc---I.}
\label{fig:torus_knot_form1}
\end{subfigure}
\begin{subfigure}{.32\linewidth}
\center
\includegraphics[scale=.07]{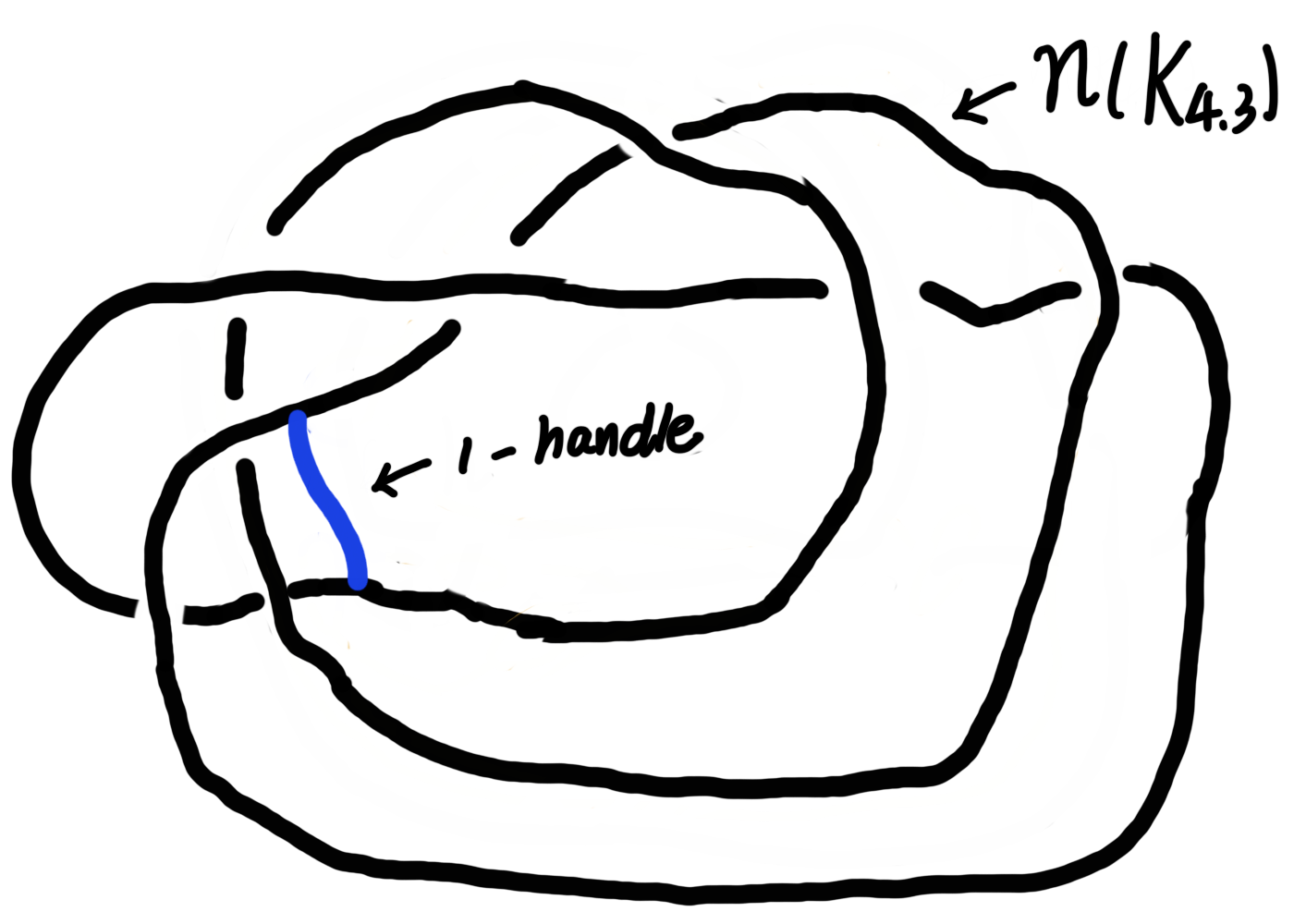}
\caption{Torus knot with an arc---II.}
\label{fig:torus_knot_form2}
\end{subfigure}
\caption{Deforming $K_l,\disk_s$.}
\end{figure}

To see how $\disk_s$ is embedded in $5_2$, 
we first observe that the boundary components
of $\annulus$ and $\disk_1$ are embedded in $\partial 5_2$ as depicted in Fig.\ \ref{fig:boundary_A_M1prime}; $l$ is hence the loop in $\partial 5_2$ shown in Fig.\ \ref{fig:boundary_Ac}. 
Let $\disk_a^+, \disk_a^-$ be components of 
the frontier of $\rnbhd{\disk_a}$ in $5_2$.
Then it may be assumed that 
the disk $\disk_s$ is
the component of 
the frontier of a regular neighborhood of $K_l:=\disk_a^+\cup \disk_a^-\cup l'$ not parallel $\disk_a^\pm$, where $l':=l-\openrnbhd{\disk_a}$. Isotope $K_l$ so that it is as shown in Fig.\ \ref{fig:DDa_lprime}.
The relation between $K_l$
and $5_2$ is shown in Fig.\ \ref{fig:Kl_fivetwo}. 
Isotope $5_2$ in $\sphere$ based on how $l'$ is twisted around $5_2$. After a series of isotopies shown 
in Figs.\ \ref{fig:deformation1}, \ref{fig:deformation2} and \ref{fig:deformation3}, 
we end up with  
the handlebody-knot in Fig.\ \ref{fig:disk_Ds} with 
a simpler expression of $\disk_s$.
%
Therefore $\pairfivetwo$ can be thought of 
as a regular neighborhood of 
$(4,3)$-torus $K_{4,3}$ 
with a $1$-handle attached as highlighted in Figs.\ \ref{fig:torus_knot_form1} and \ref{fig:torus_knot_form2}, so 
the slope of $\annulus_c$ is $\frac{4}{3}$, and hence 
the following.
\begin{theorem}\label{teo:hkfivetwo_ann_diag}
The annulus diagram of $\pairfivetwo$ is 
\raisebox{-.07\height}{\includegraphics[scale=.17]{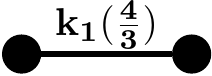}}.
\end{theorem}

Note that by \cite[Theorem $1.1$]{Wan:22p}
the annulus diagram of $\pair$ 
is \raisebox{-.07\height}{\includegraphics[scale=.17]{char_diag_stick}} if and only if 
its exterior $\Compl\HK$ admits infinitely many
essential annuli. In this case, all essential annuli 
but one are separating with the unique non-separating one being necessarily 
of type $3$-$3$ by \cite[Theorem $1.4$]{Wan:22p}, and all essential annuli but one 
is non-characteristic with the characteristic one being of 
type $3$-$2$i by Lemma \ref{lm:cano_ann_stick_type}. 
For such a $\pair$, it is then interesting to consider the following questions.
\begin{question}
What are the types of the non-separating, non-characteristic annuli in $\Compl\HK$?
\end{question}
\begin{question}
Are there any constraints on the slope pair of 
the non-separating annulus or on the slope
of the characteristic annulus in $\Compl\HK$?
\end{question}
 



\section{Annulus diagrams of Handlebody-knot families}\label{sec:families}
\subsection{Handlebody-knot families}\label{subsec:hk_family}
To construct an infinite family of handlebody-knots
with homeomorphic exteriors, we consider the following 
twisting operation.

\begin{definition}
A \textit{slicing} surface of a handlebody-knot $\pair$
is a pair $(\slice, c)$, where $\slice$ a planar surface
in $\Compl\HK$ with at most two components of $\partial \slice$ not bounding a disk in $\HK$ and
$c$ being one such component. 

\end{definition}

\begin{figure}
\begin{subfigure}{.47\linewidth}
\center
\includegraphics[scale=.12]{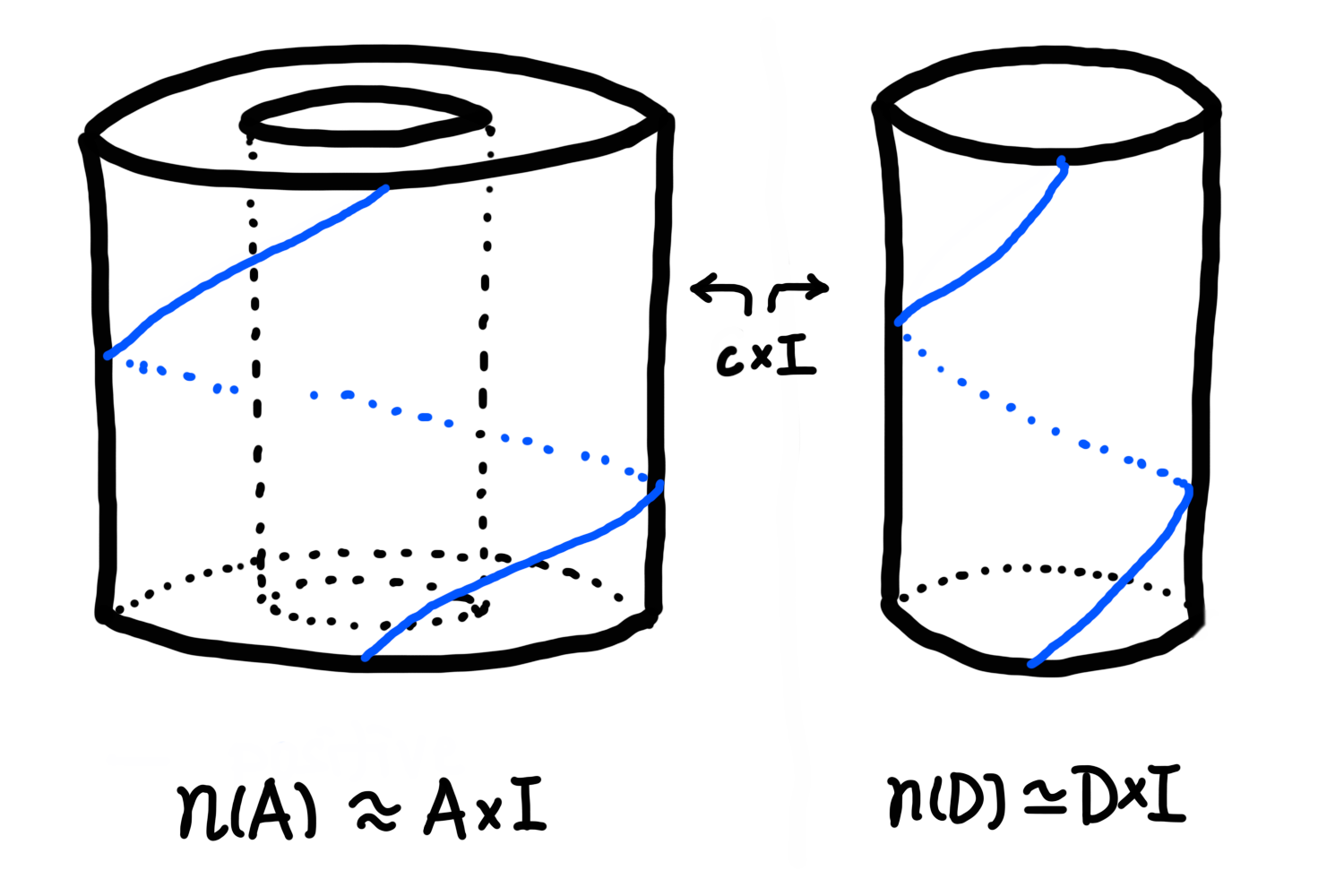}
\caption{$(+)$-twisting.}
\label{fig:p_twisting}
\end{subfigure}
\begin{subfigure}{.47\linewidth}
\center
\includegraphics[scale=.12]{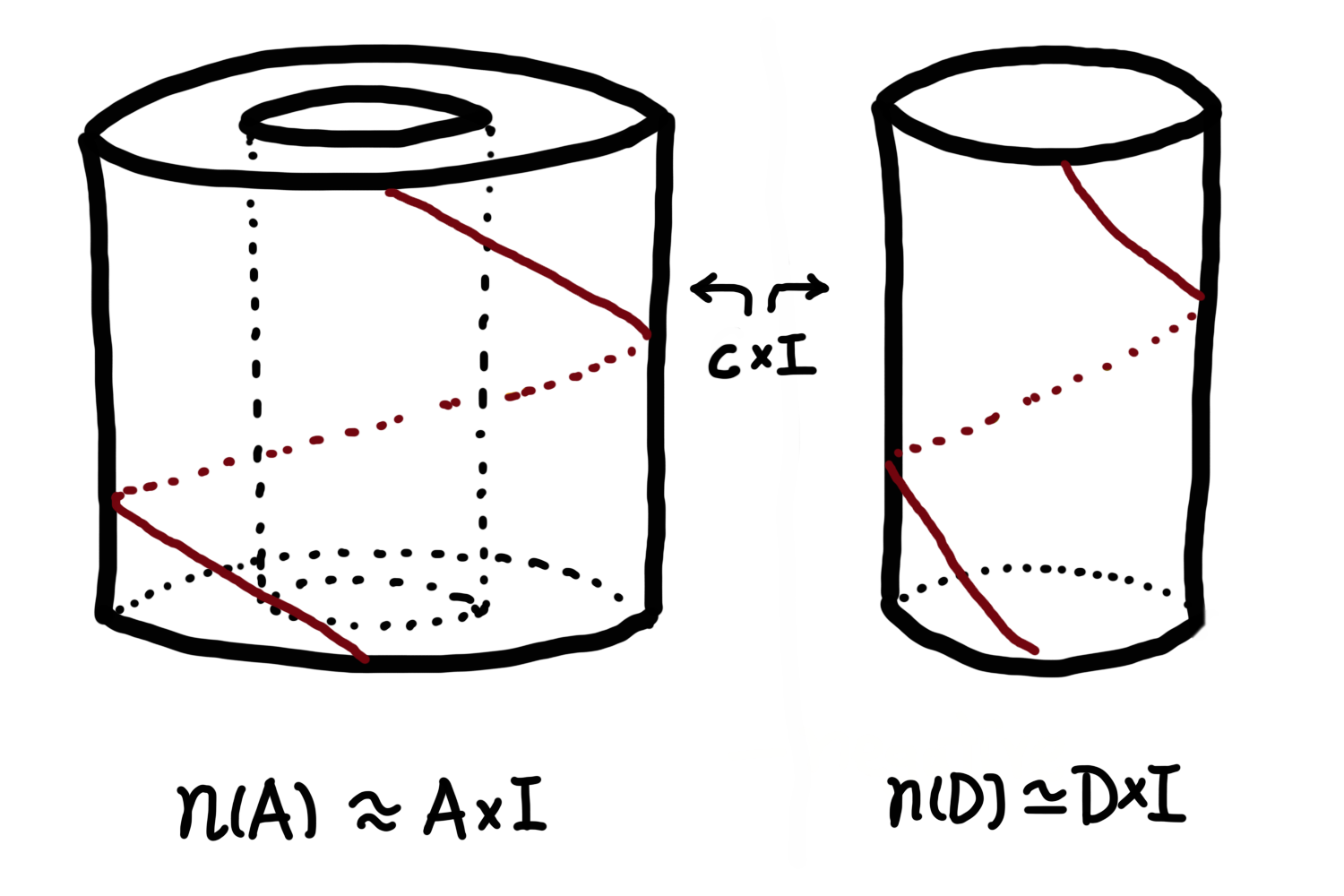}
\caption{$(-)$-twisting.}
\label{fig:n_twisting}
\end{subfigure}
\caption{Sign of a twisting.}
\label{fig:pm_twisting}
\end{figure}

Let $c_1,\dots, c_n$ be the components of $\partial \slice$ that bound a disk in $\HK$ 
and $\disk_i\subset \HK$, $i=1,\dots, n$, be disks 
bounded by them. 
Then $\tF:=\slice\bigcup_{i=1}^n \disk_i$ is either an
annulus $\tA$ or a disk $\tD$, called 
a \emph{twisting} annulus or disk of $\pair$.
%
%
Consider now the union $M'$ of $\Compl\HK$ 
and a regular neighborhood of $\mathbf{D}:=\bigcup_{i=1}^n \disk_i\subset \HK$, and choose 
a regular neighborhood $\rnbhd{\tF}$ of $\tF\subset M'$ 
so that $\tF\cap \HK$ is a regular neighborhood of $\mathbf{D}\subset\HK$. Let $M$ be 
the union of $\Compl\HK$ and $\rnbhd{\tF}$. 
Then the twisting operation is to reembeds
the handlebody-knot exterior $\Compl\HK$ via the following composition
\begin{equation}\label{eq:twisting_map}
T_n:\Compl\HK\subset M\xrightarrow{t_n} M\subset \sphere,
\end{equation}
where $t_n$ is given by twisting $n$ times along $\tF$ with the sign convention
in Fig.\ \ref{fig:pm_twisting}; note that the convention depends on the selected component $c\subset \partial \slice$. 
Set $M_n:=T_n(\Compl\HK)\subset\sphere$. 
\begin{lemma}
$M_n$ is the exterior of some handlebody-knot. 
\end{lemma}
\begin{proof}
Since $t_n$ is a self-homeomorphism of $M$, $t_n(M)=M\subset \sphere$, and hence $t_n(M)$ is the exterior of 
the union of some $3$-balls $B$ and 
a $2$-component link or knot $L$. On the other hand, 
$\Compl\HK$ is the exterior of some arcs $\gamma$ in $M$,
so $M_n$ is the exterior of the handlebody
$B \cup \rnbhd{L}\cup\rnbhd{t_n(\gamma)}$, whose boundary 
is necessarily of genus $2$. 
\end{proof}

The resulting handlebody-knot   
is said to be 
obtained by twisting $\pair$ along $\tF$ $n$ times.
Its exterior is homeomorphic to $\Compl\HK$ as $M_n$ is homeomorphic to $\Compl\HK$.
%
\begin{figure}[t]
\begin{subfigure}{.32\linewidth}
\center
\includegraphics[scale=.1]{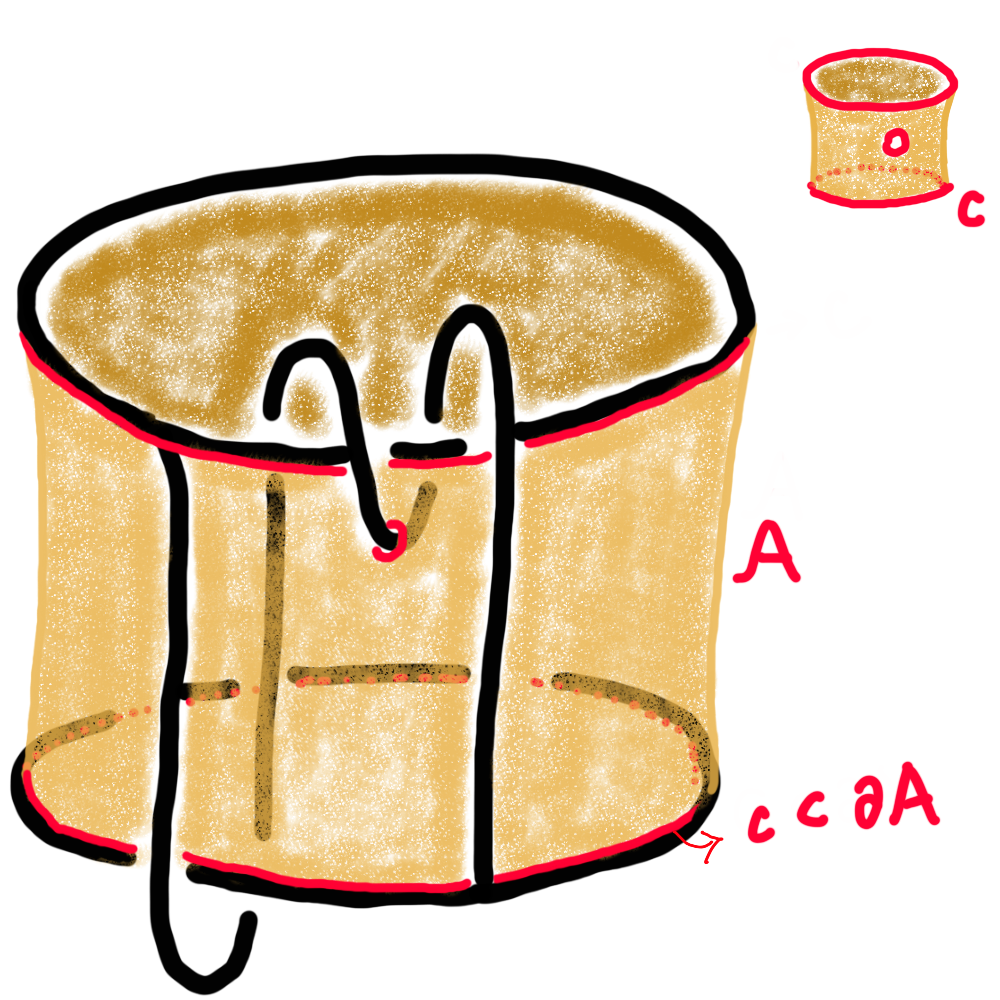}
\caption{Twisting annulus $\tA$.}
\label{fig:hksixone_A_flat}
\end{subfigure}
\begin{subfigure}{.33\linewidth}
\center
\includegraphics[scale=.2]{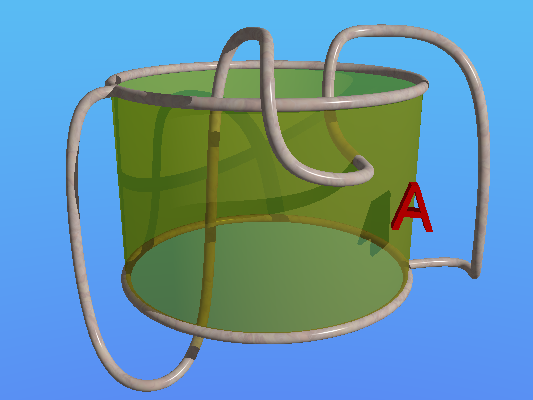}
\caption{Twisting annulus $\tA$ in 3D.}
\label{fig:hksixone_A}
\end{subfigure}
\begin{subfigure}{.33\linewidth}
\center
\includegraphics[scale=.1]{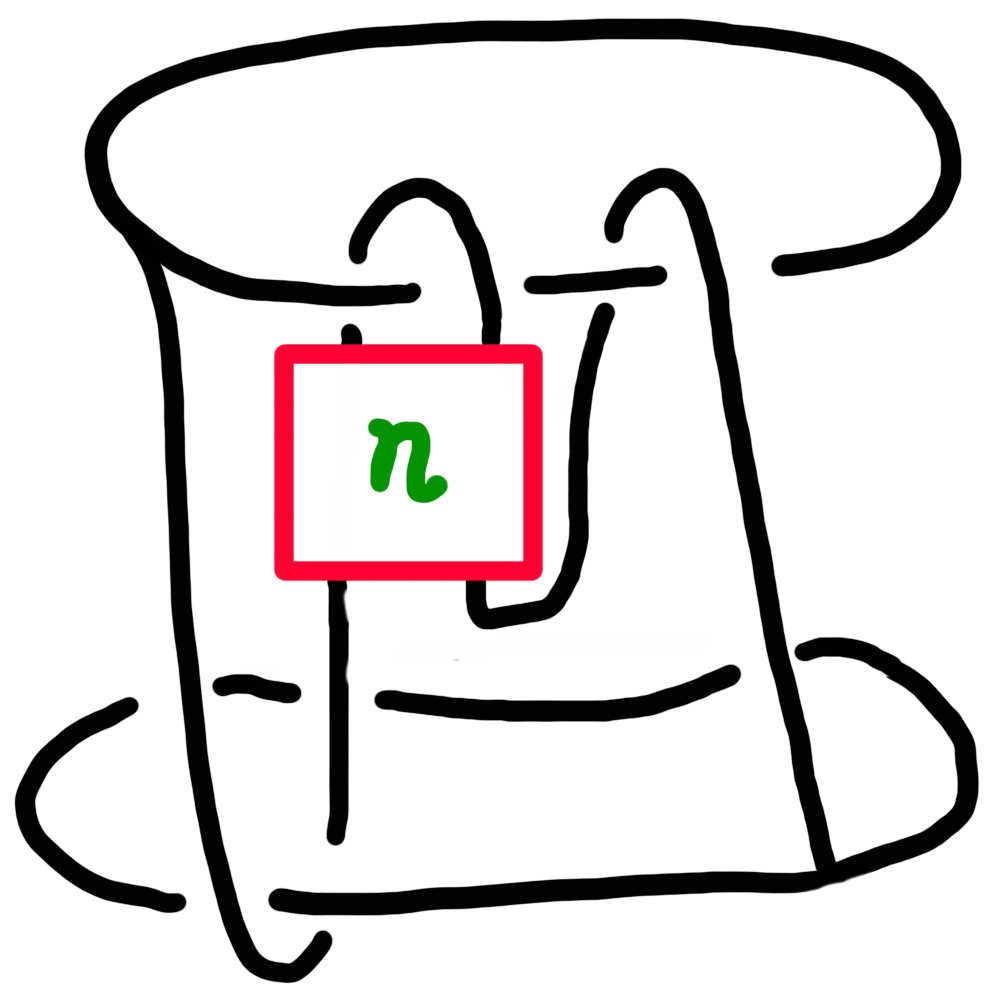}
\caption{$\pairMo n$.}
\label{fig:motto_hk_family}
\end{subfigure}
\caption{Twisting annulus $\tA$ of $\pairsixone$.}
\label{fig:twisting_ann_hksixone}
\end{figure}

\subsubsection*{Motto's handlebody-knot family} 
Consider the handlebody-knot $\pairsixone$ 
and the twisting annulus $\tA$ depicted in Figs.\ \ref{fig:hksixone_A_flat} and \ref{fig:hksixone_A}.
Twisting $\pairsixone$ $n$ times along $\tA$, we 
obtain Motto's handlebody-knots family $\familyMo$ in 
\cite{Mott:90} (see Fig.\ \ref{fig:motto_hk_family}).
  
   
\begin{figure}[t]
\center
\includegraphics[scale=.1]{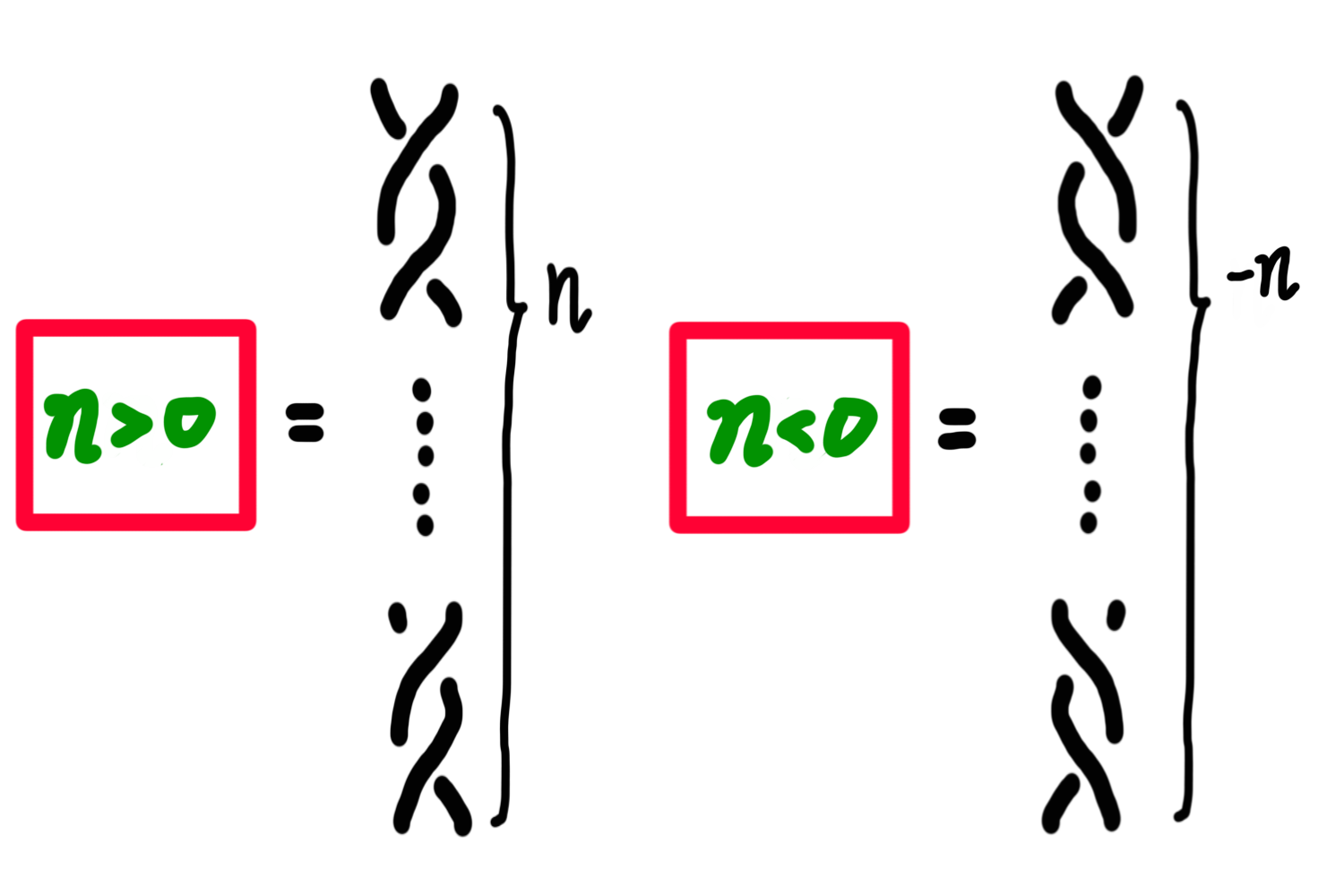}
\caption{Sign convention.}
\label{fig:pm_crossings.}
\end{figure}

\subsubsection*{LeeLee's handlebody-knot family I}
\begin{figure}[t]
\begin{subfigure}{.45\linewidth}
\center
\includegraphics[scale=.2]{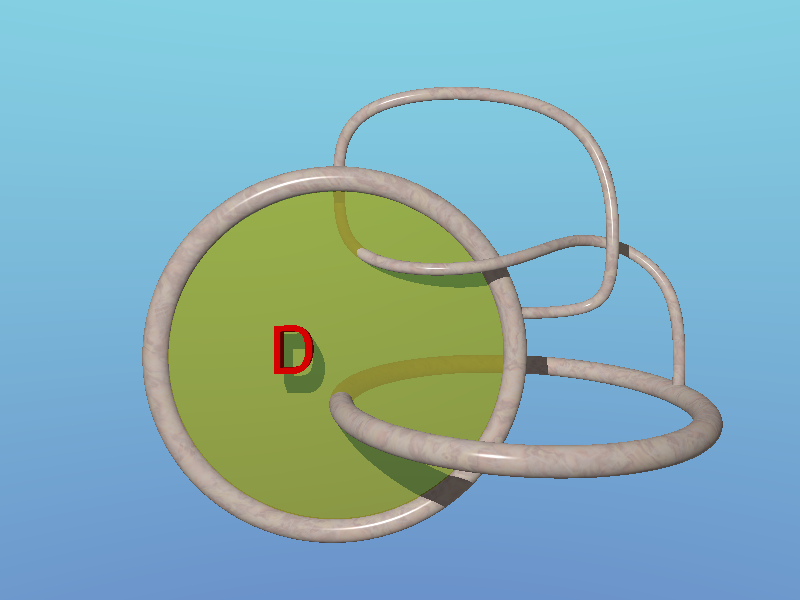}
\caption{Twisting disk $\tD$ of $\pairfiveone$.}
\label{fig:hkfiveone_D}
\end{subfigure}
\begin{subfigure}{.45\linewidth}
\center
\includegraphics[scale=.1]{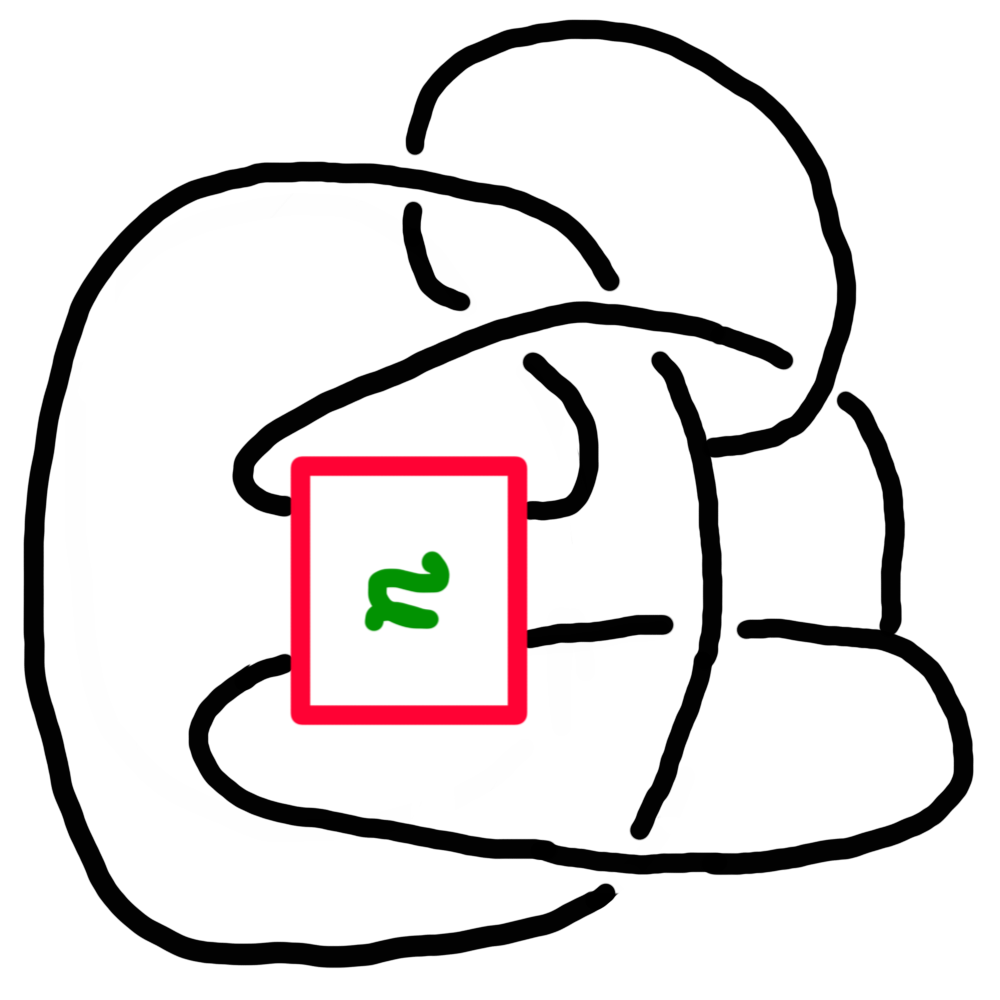}
\caption{$\pairLLone n$.}
\label{fig:LLone_hk_family}
\end{subfigure}
\caption{Twisting disk $\tD$ and $\pairLLone n$.}
\end{figure}

Similarly, 
the first Lee-Lee handlebody-knot family $\familyLLone$ (see Fig.\ \ref{fig:LLone_hk_family}) in \cite{LeeLee:12}
can be constructed by twisting 
the handlebody-knot $\pairfiveone$ 
along the twisting disk
$\tD$ shown in Fig.\ \ref{fig:hkfiveone_D} $n$ times.

\subsubsection*{A variant}
\begin{figure}
\begin{subfigure}{.33\linewidth}
\center
\includegraphics[scale=.18]{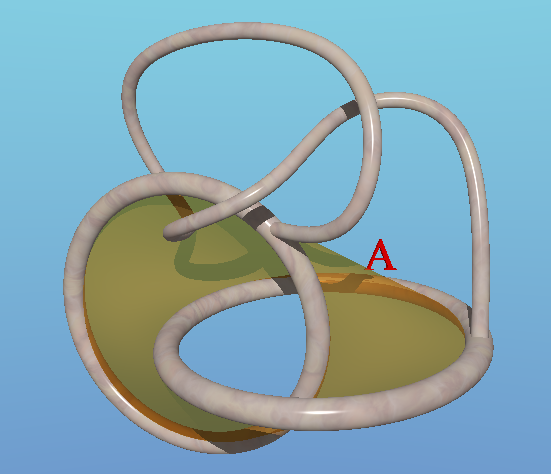}
\caption{Twisting annulus $\tA$.}
\label{fig:hkfiveone_A}
\end{subfigure}
\begin{subfigure}{.32\linewidth}
\center
\includegraphics[scale=.07]{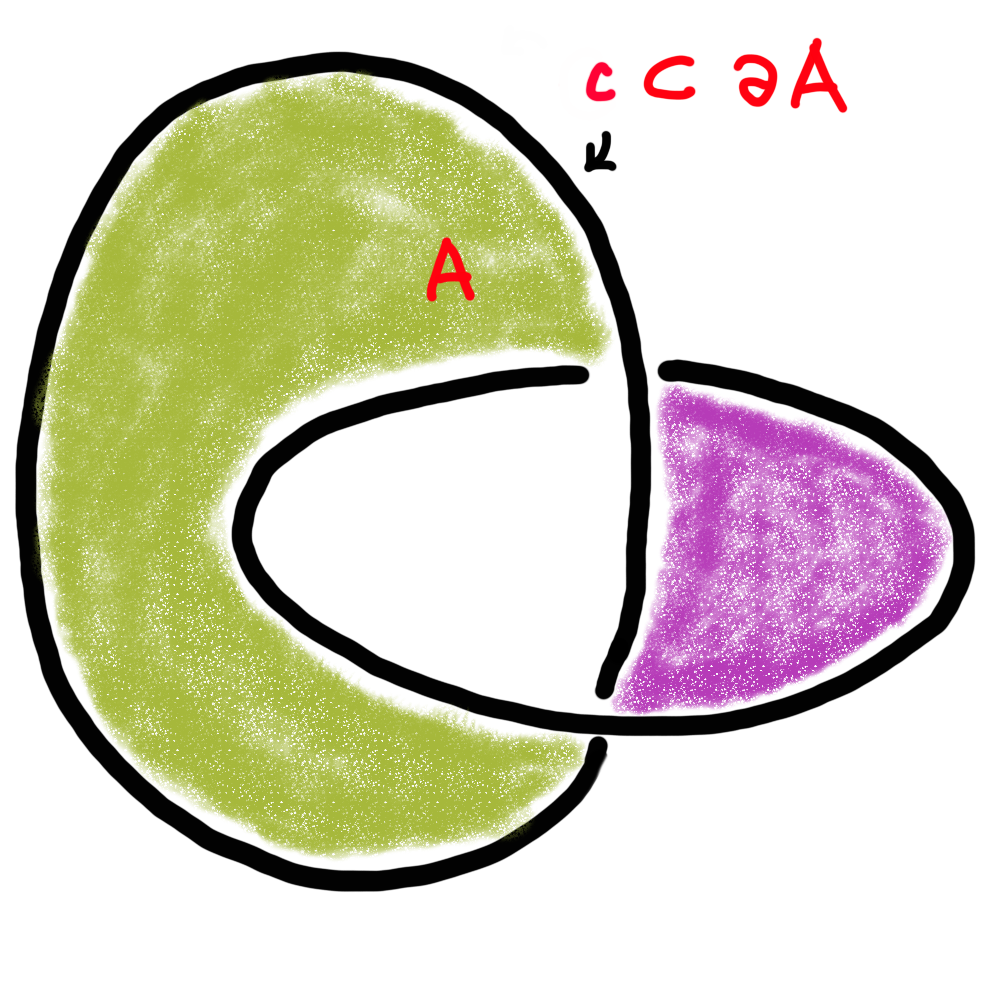}
\caption{Marked component $C\subset \partial \tA$.}
\label{fig:hkfiveone_A_highlight}
\end{subfigure}
\begin{subfigure}{.31\linewidth}
\center
\includegraphics[scale=.08]{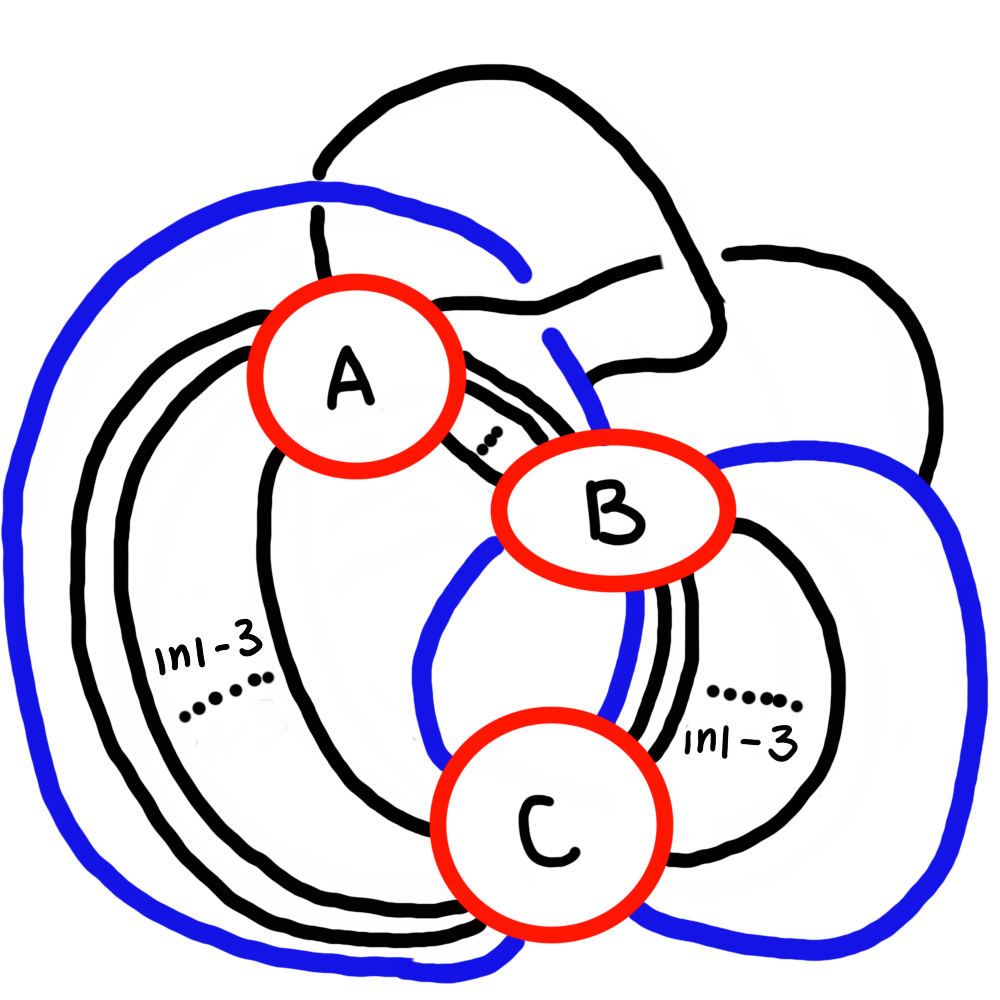}
\caption{$\pairLLvariant n$.}
\label{fig:LL1_variant}
\end{subfigure}
\begin{subfigure}{.24\linewidth}
\center
\includegraphics[scale=.07]{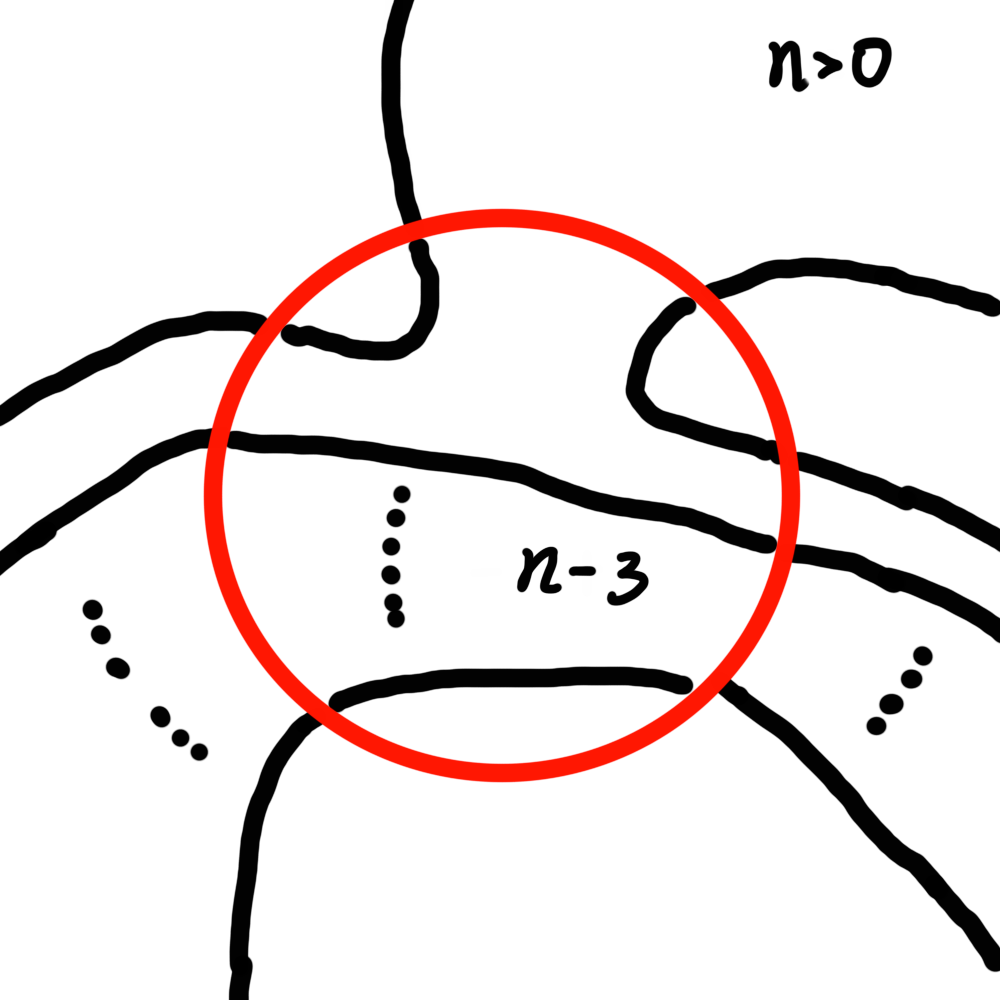}
\caption{Tangle $A$, $n>0$.}
\label{fig:spider_pA}
\end{subfigure}
\begin{subfigure}{.24\linewidth}
\center
\includegraphics[scale=.07]{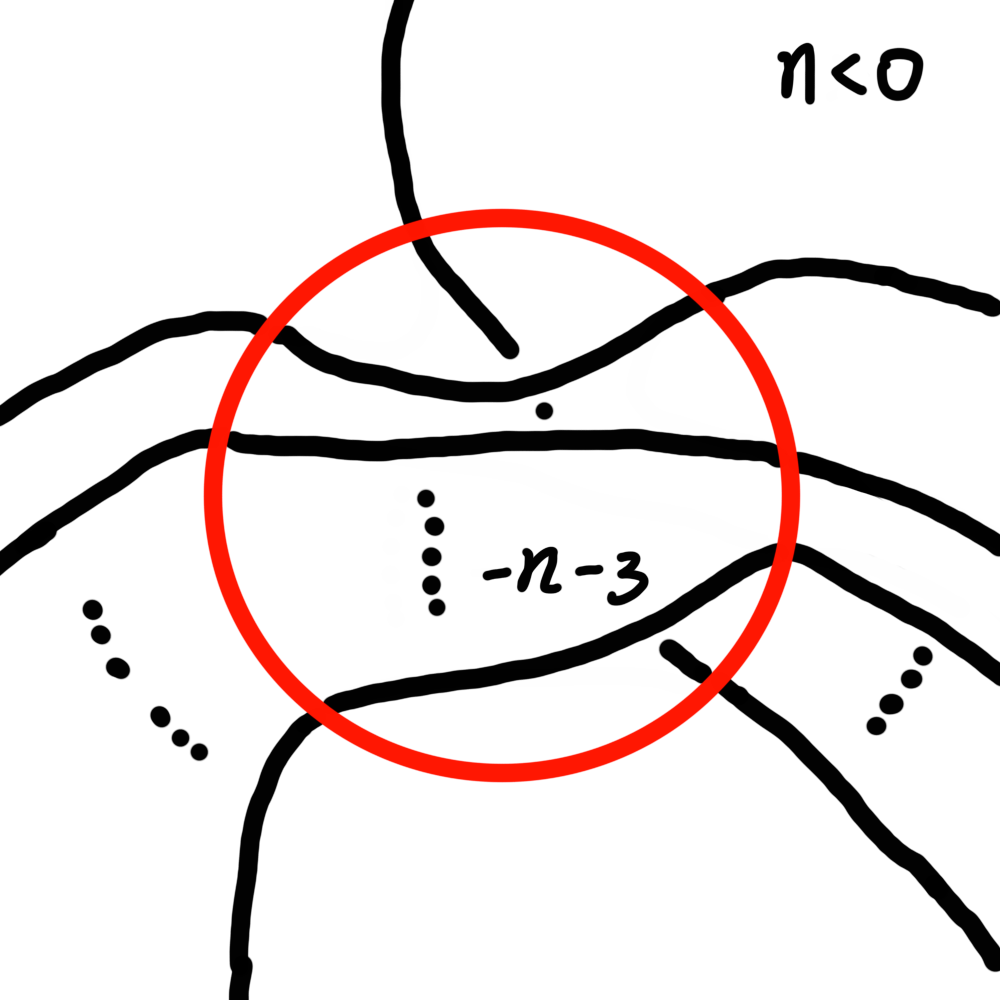}
\caption{Tangle $A$, $n<0$.}
\label{fig:spider_mA}
\end{subfigure}
\begin{subfigure}{.24\linewidth}
\center
\includegraphics[scale=.07]{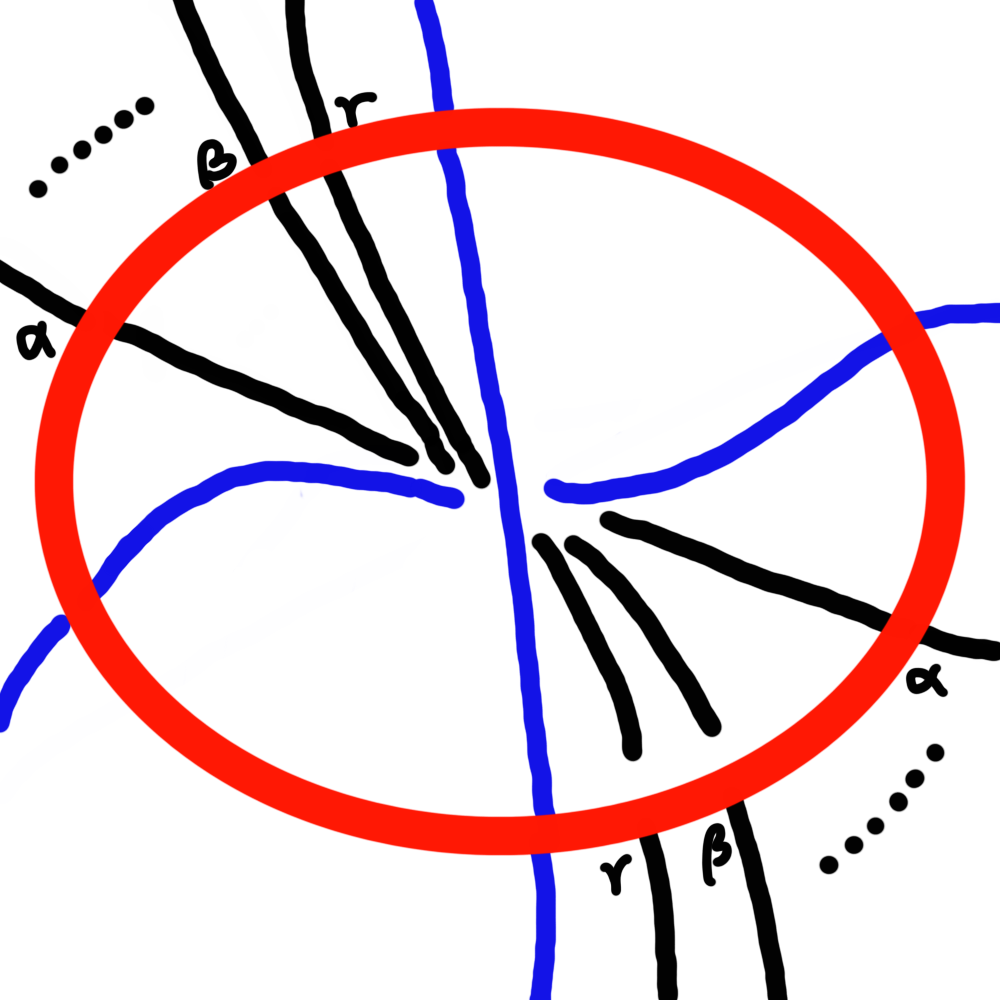}
\caption{Tangle $B$.}
\label{fig:spider_B}
\end{subfigure}
\begin{subfigure}{.24\linewidth}
\center
\includegraphics[scale=.07]{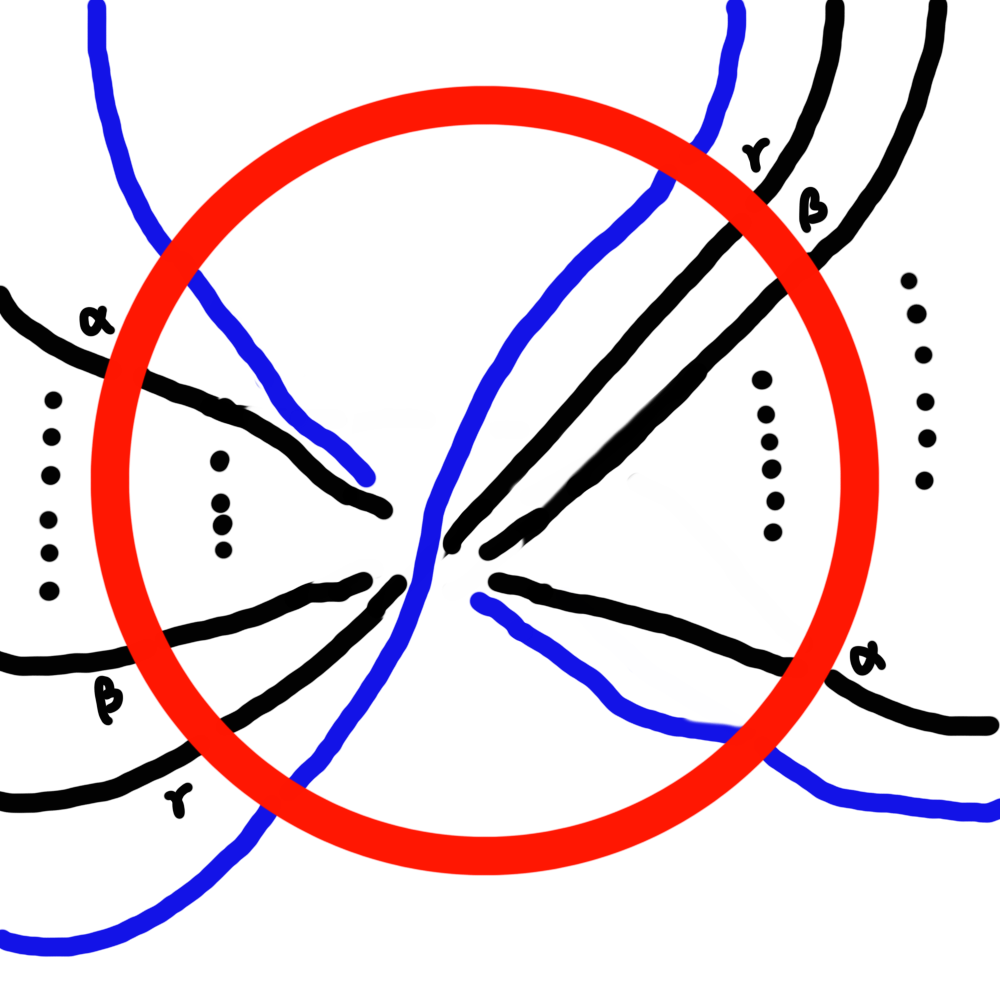}
\caption{Tangle $C$.}
\label{fig:spider_C}
\end{subfigure}
\caption{Twisting annulus $\tA$ and $\pairLLvariant n$.}
\end{figure}

Note that $\pairfiveone$ also admits a 
twisting annulus $\tA$ as depicted in Figs.\ \ref{fig:hkfiveone_A}, \ref{fig:hkfiveone_A_highlight}, and
we denote by $\familyLLvariant$, 
the family of handlebody-knots 
obtained by twisting $\pairfiveone$ along $\tA$ $n$ times  
(see Figs.\ \ref{fig:LL1_variant}--\ref{fig:spider_C}).
In particular, 
the exteriors of $\pairLLvariant n$ and $\pairLLone n$ are homeomorphic to $\Compl {5_1}$, for every $n$. 

\subsubsection*{LeeLee's handlebody-knot family II} 
\begin{figure}[t]
\begin{subfigure}{.32\linewidth}
\center
\includegraphics[scale=.08]{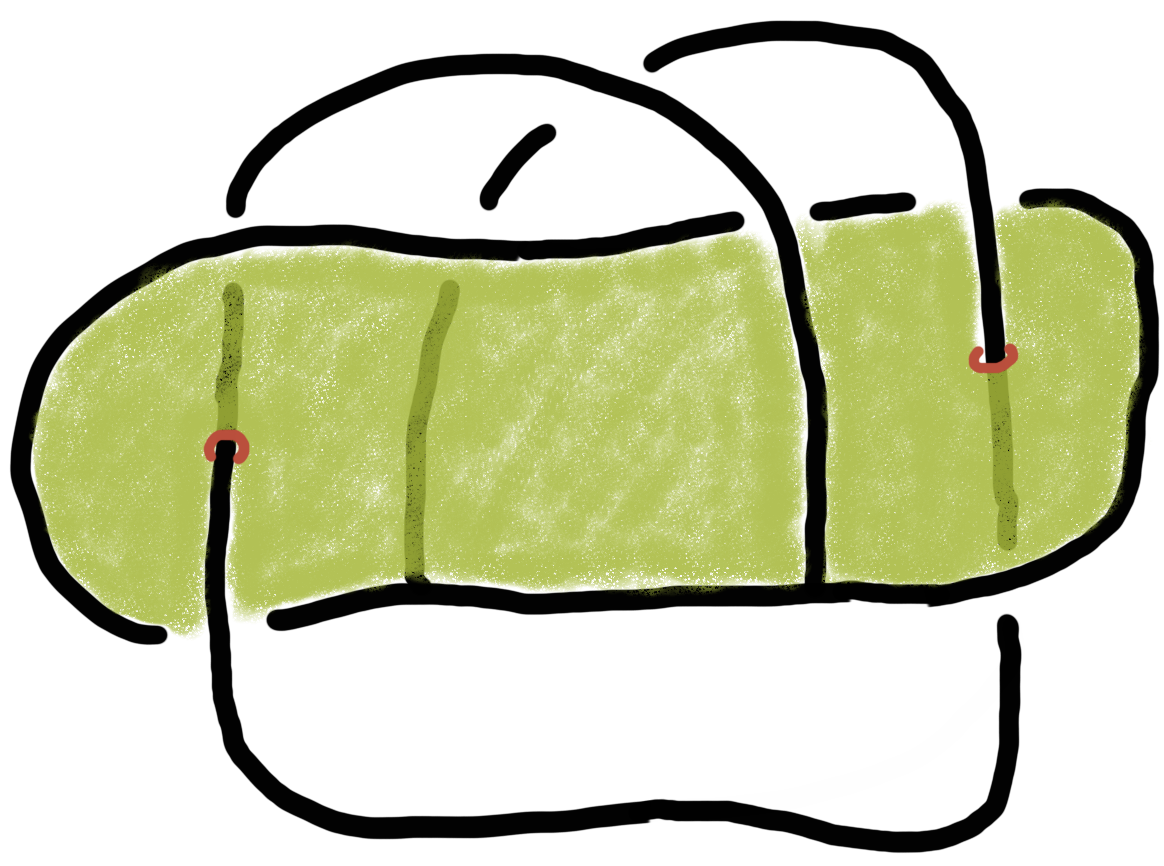}
\caption{Twisting disk $\tD$.}
\label{fig:hkfivetwo_D_flat}
\end{subfigure}
\begin{subfigure}{.32\linewidth}
\center
\includegraphics[scale=.2]{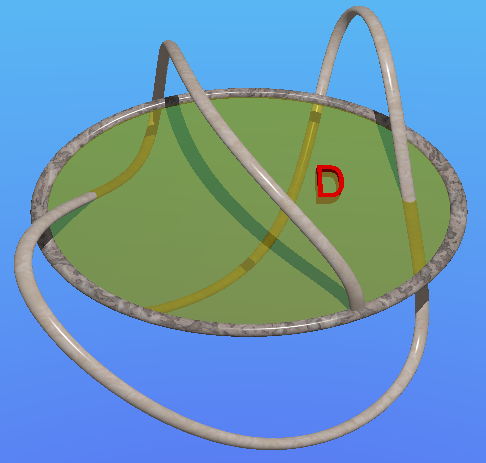}
\caption{Twisting disk $\tD$ in 3D.}
\label{fig:hkfivetwo_D_3d}
\end{subfigure}
\begin{subfigure}{.32\linewidth}
\center
\includegraphics[scale=.1]{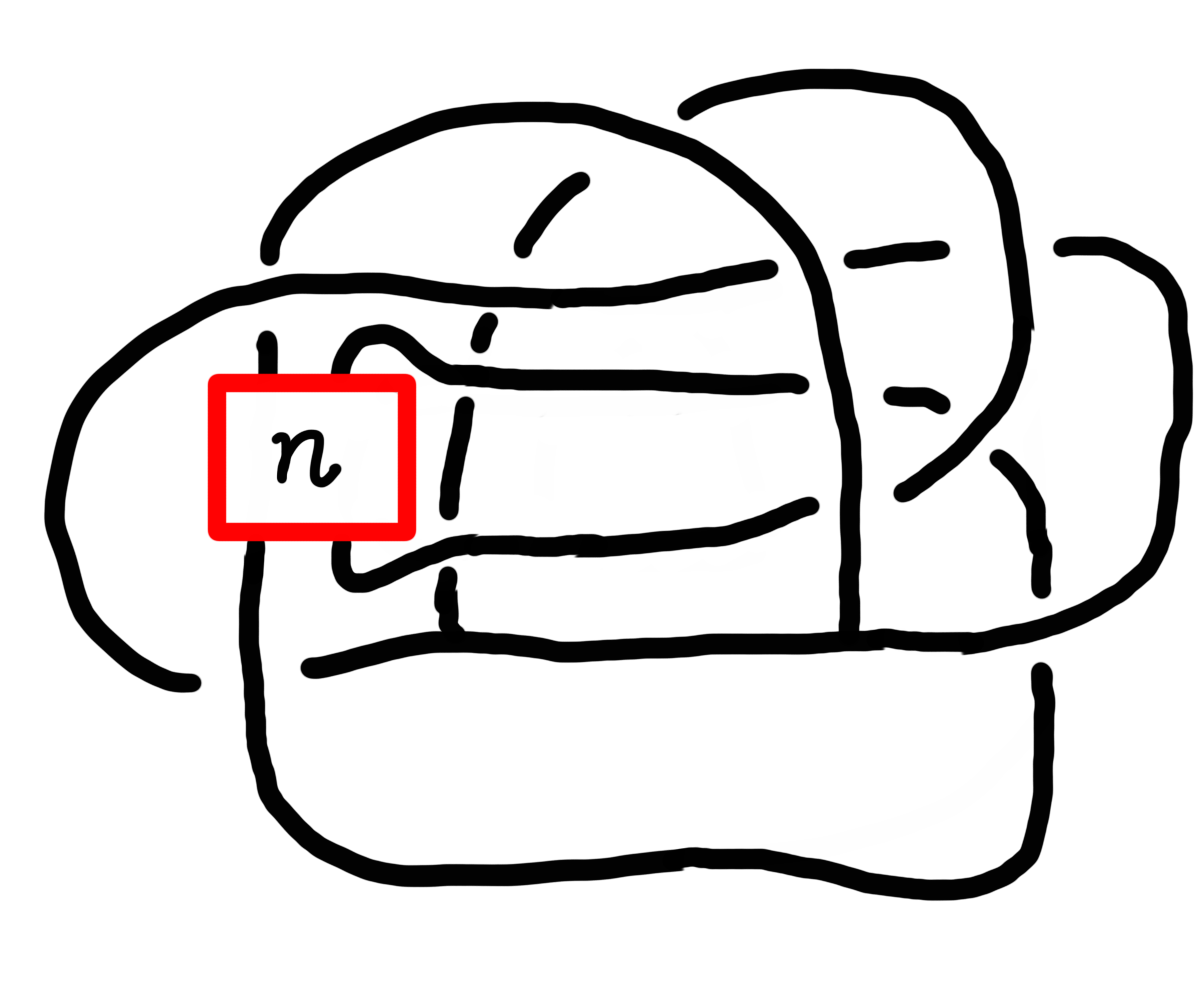}
\caption{$\pairLLtwo n$.}
\label{fig:LLtwo_hk_family}
\end{subfigure}
\caption{Twisting disk $\tD$ and $\pairLLtwo n$.}
\label{fig:hkfivetwo_D}
\end{figure}

Lastly, we consider the twisting disk $\tD$ 
of $\pairfivetwo$ shown in Figs.\ \ref{fig:hkfivetwo_D_flat}--\ref{fig:hkfivetwo_D_3d}.
Twisting $\pairfivetwo$ along 
$\tD$ $n$ times yields the second handlebody-knot family 
$\familyLLtwo$ in \cite{LeeLee:12}  
(see Fig.\ \ref{fig:LLtwo_hk_family}).   

\subsection{Annulus diagram}\label{subsec:ann_diag}

Here we compute the annulus diagram 
for handlebody-knot families in Section \ref{subsec:hk_family}.

\begin{figure}[t]
\begin{subfigure}{.32\linewidth}
\center
\includegraphics[scale=.07]{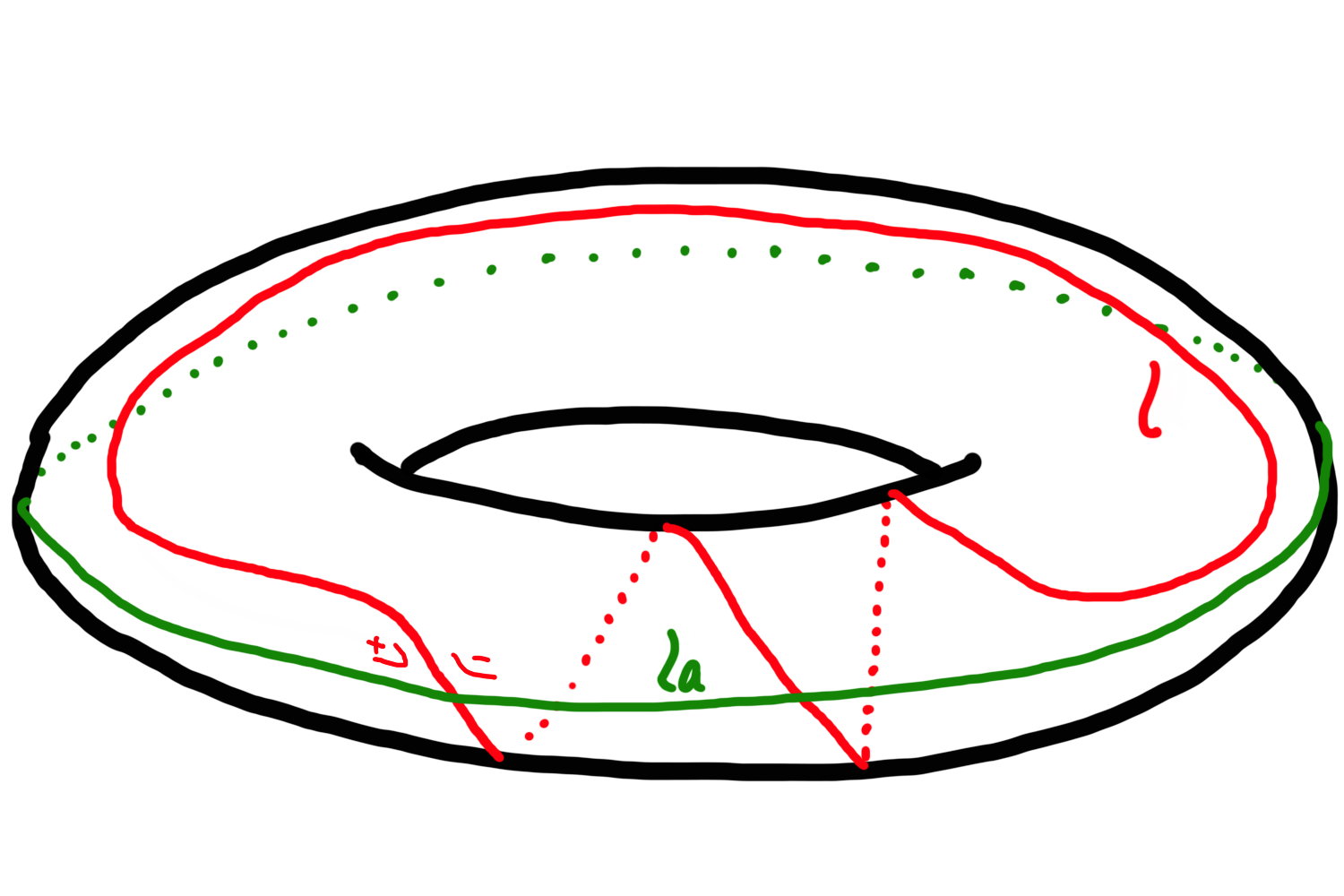}
\caption{Intersection $l\cap l_a$.}
\label{fig:mo_intersection}
\end{subfigure}
\begin{subfigure}{.32\linewidth}
\center
\includegraphics[scale=.1]{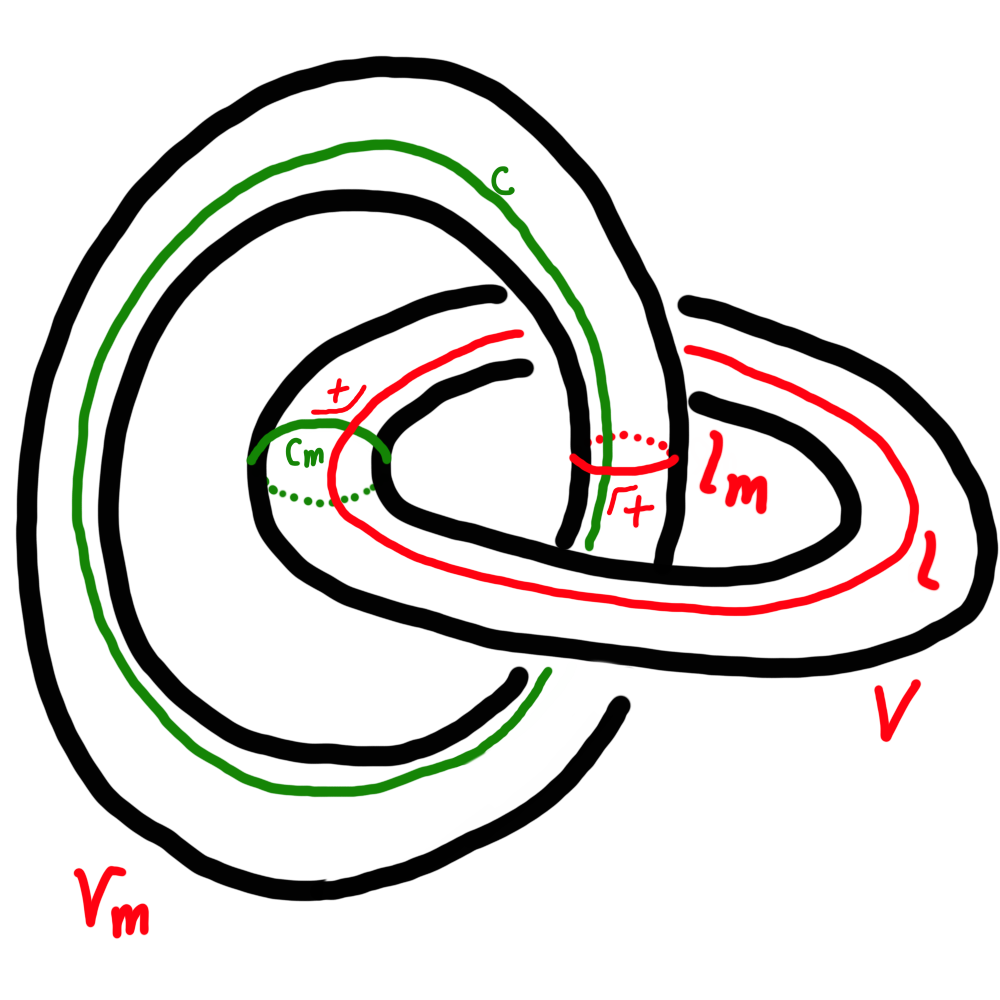}
\caption{Intersection $\partial \annulus$ and $c\cup c_m$.}
\label{fig:LLone_intersection}
\end{subfigure}
\begin{subfigure}{.32\linewidth}
\center
\includegraphics[scale=.1]{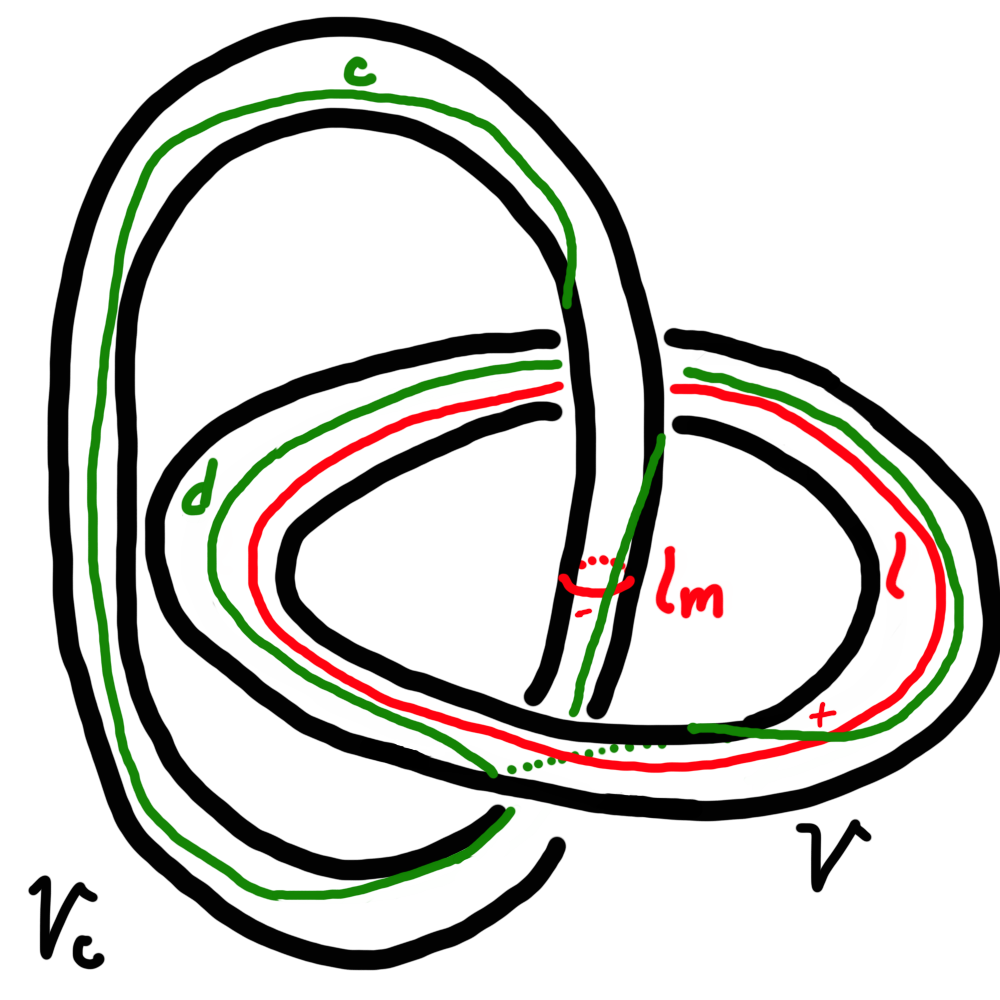}
\caption{Intersection $\partial \annulus\cap \partial\tA$.}
\label{fig:LLvariant_intersection}
\end{subfigure}
\begin{subfigure}{.32\linewidth}
\center
\includegraphics[scale=.1]{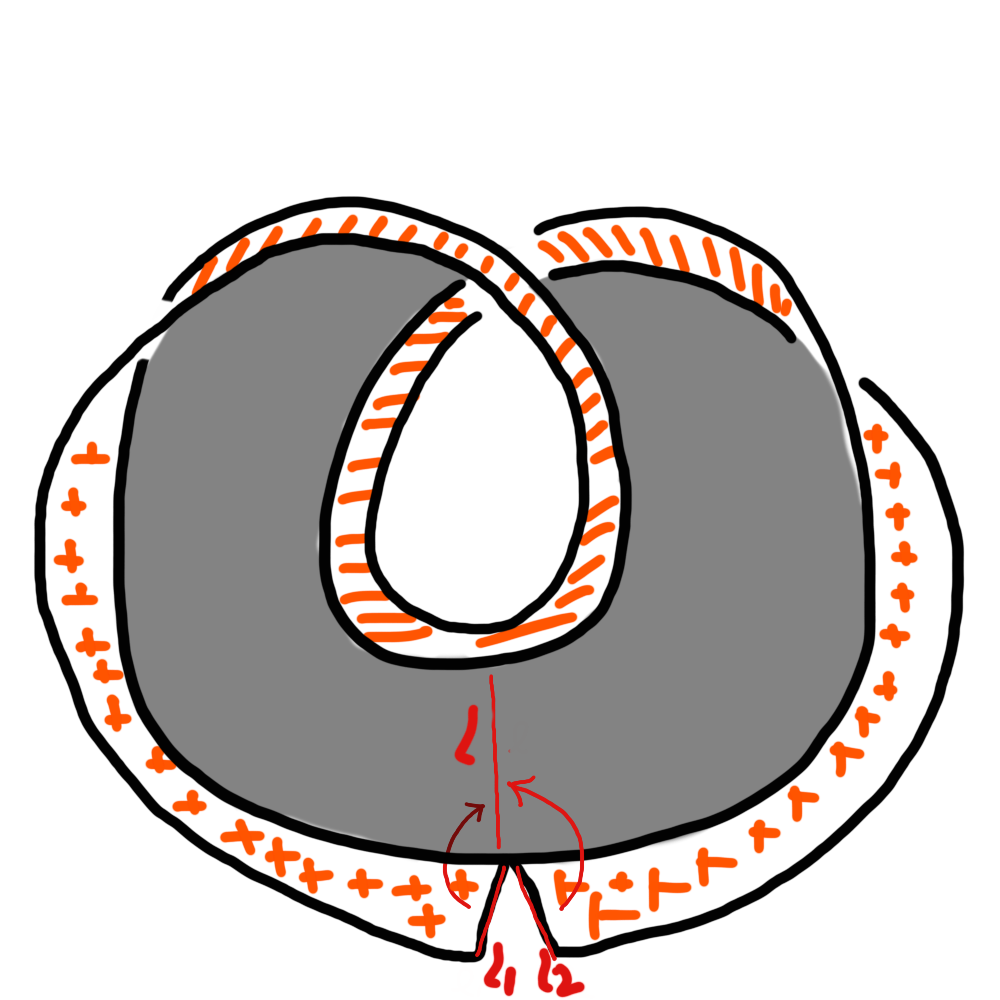}
\caption{Union $5_2\cup \annulus\cup \mobius$.}
\label{fig:union}
\end{subfigure}
\begin{subfigure}{.32\linewidth}
\center
\includegraphics[scale=.1]{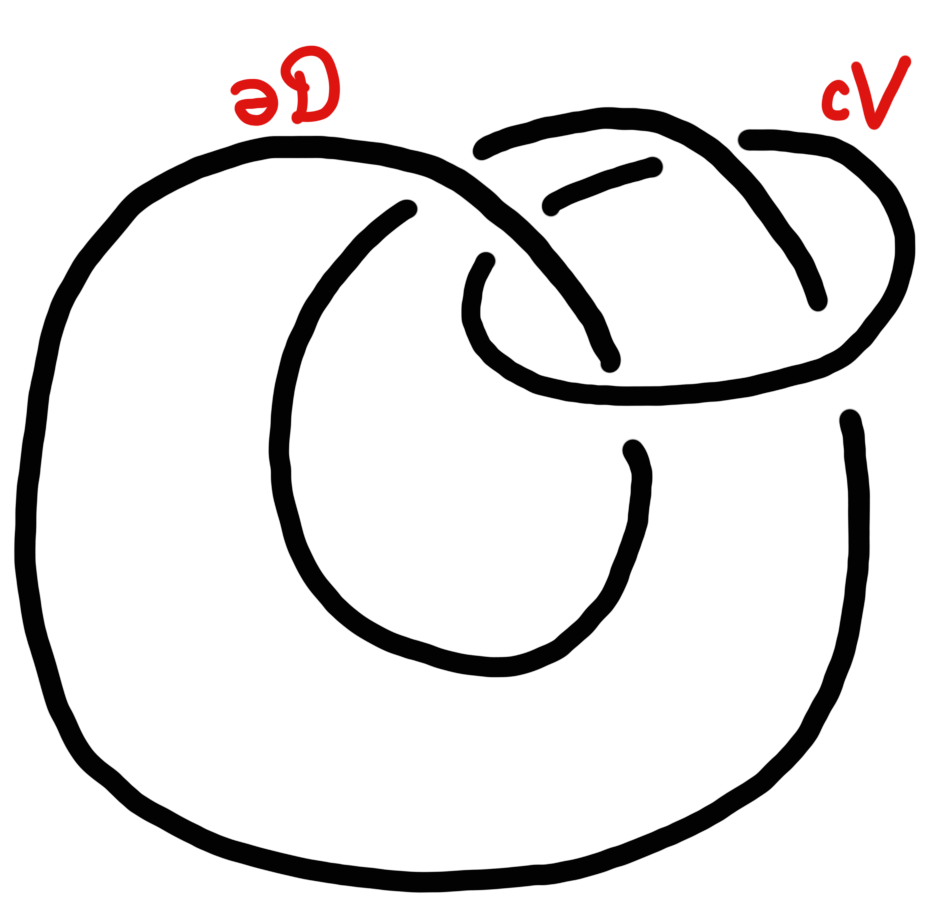}
\caption{$\partial \tD$ and $c V$.}
\label{fig:btDandcA}
\end{subfigure}
\begin{subfigure}{.32\linewidth}
\center
\includegraphics[scale=.1]{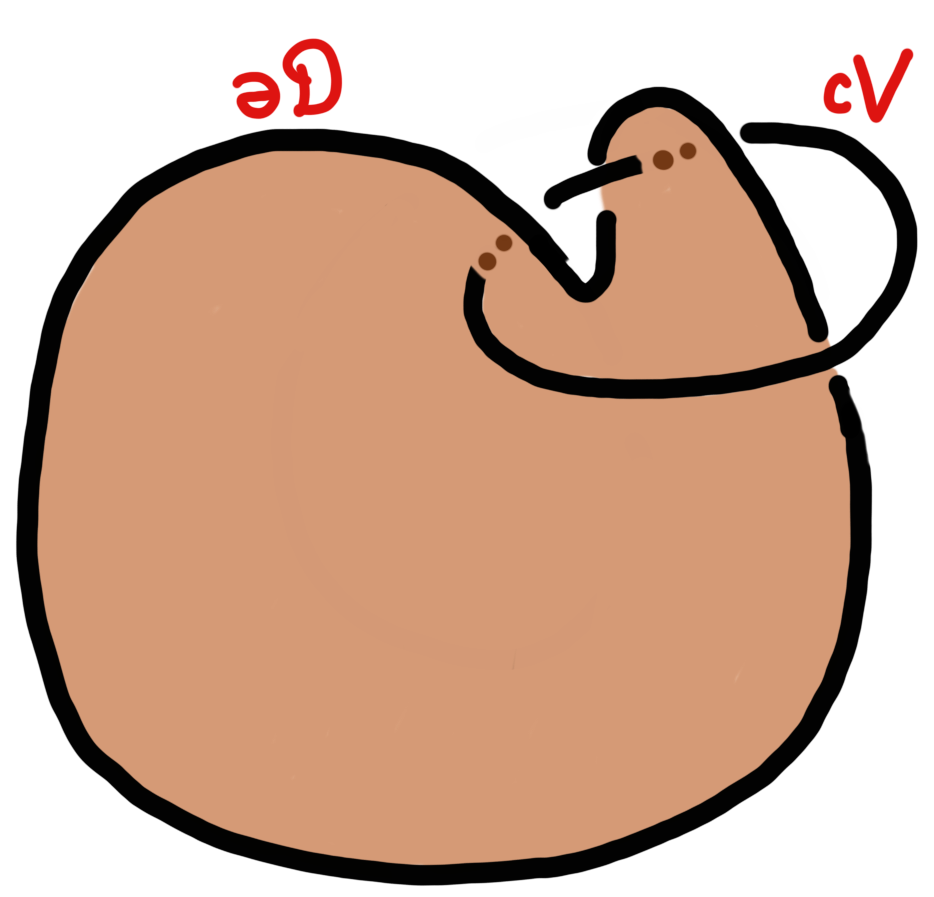}
\caption{$\tD$ and $c V$.}
\label{fig:tDandcA}
\end{subfigure}
\caption{Annulus diagram computation.}
\end{figure}

\begin{theorem}\label{teo:ann_diag_motto}
The annulus diagram of $\pairMo n$ is   
\raisebox{-.4\height}{\includegraphics[scale=.2]{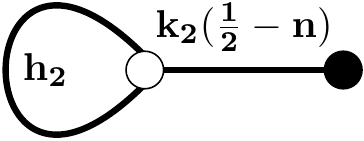}} . 
\end{theorem}
\begin{proof}
Note first that 
the interior of the twisting annulus $\tA$ of $\pairsixone$ 
meets $\HK$ at a separating disk $\disk_s\subset\HK$ 
and $\HK-\openrnbhd{\disk_s}$ consists of 
two solid tori $V,V'$,
so what the twisting map $t_n$ does is to 
reembed the $1$-handle $\rnbhd{\disk_s}$ connecting $\partial V,\partial V'$
in the exterior $M:=\Compl{V\cup V'}$. 
Denote by $\annulus, \annulus'\subset\Compl{6_1}$ 
the characteristic essential annuli of type $3$-$2$
and of type $2$-$2$, respectively, and by 
$\annulus_n, \annulus_n'$ the annuli 
$T_n(\annulus),T_n(\annulus')\subset\Compl {\Mo n}$,
respectively. We may assume that $\partial \annulus\subset V$; this implies $\partial \annulus'\cap V'\neq\emptyset$. 
Now, since $\tA\cap \annulus'=\emptyset$, 
$\annulus_n'\subset \Compl{\Mo n}$ is of type $2$-$2$, whereas, $(\sphere,V)$ being 
a trivial solid torus knot in $\sphere$, $\annulus_n\subset\Compl{\Mo n}$ 
is of type $3$-$2$ii.
On the other hand, $\partial\tA\cap \partial \annulus\neq\emptyset$ and the intersection of the component $l_a=V\cap \tA$ and a component $l$ of $\annulus$ is shown in Fig.\ \ref{fig:mo_intersection}, wherefrom
we deduce $T_n(l)\subset V$ 
has a slope of $\frac{2}{1-2n}$ with respect to $(\sphere, V)$ and hence the theorem.
\end{proof}

\begin{theorem}\label{teo:ann_diag_leeleeone}
The annulus diagram of $\pairLLone n$, $n\neq 0$, is \raisebox{-.4\height}{\includegraphics[scale=.2]{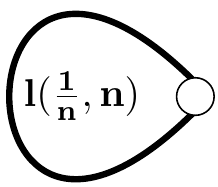}} .
\end{theorem}
\begin{proof}
Let $\annulus$ be the type $2$-$1$ annulus in $\Compl{5_1}$,
and denote by $l_m$ the component of $\partial \annulus$
that bounds a disk in $5_1$ and by $l$ the other component.
Let $\disk_s$ be the separating disk of $\HK$ in 
$\HK\cap\tD$. Then $\disk_s$ 
separates $\HK$ into two solid tori
$V,V_m$ with $l_m\subset \partial V_m, l\subset \partial V$.
Denote by $c_m$ the intersection loop $\tD\cap \partial V$
and $c=\partial \tD$; the intersection of $l\cup l_m$ and 
$c\cup c_m$ is shown in Fig.\ \ref{fig:LLone_intersection})
from where we see that
the slope of $T_n(l_m)$ with respect to $(\sphere, V_m)$
is $\frac{1}{n}$, and the slope of $T_n(l)$ with respect to $(\sphere, V)$ is $n$, so $\annulus_n$ is of type $3$-$3$ with a slope pair $(\frac{1}{n},n)$.
\end{proof}

\begin{remark}
The annulus $\annulus$ is precisely the annulus $\annulus_0$ in \cite{LeeLee:12}, and \cite{LeeLee:12} uses the fact
that $T_n(\partial\annulus_0)$ is a $(2,2n)$-torus link 
to distinguish members in $\familyLLone$.  
\end{remark}

\begin{theorem}\label{teo:ann_diag_leeleevariant}
The annulus diagram of $\pairLLvariant n$, $n\neq 0$, is \raisebox{-.4\height}{\includegraphics[scale=.2]{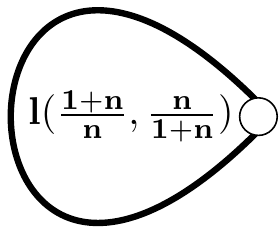}} .
\end{theorem}
\begin{proof}
Let $\annulus$, $\partial \annulus=l\cup l_m$ be 
as in the preceding proof. 
The interior of the twisting annulus $\tA$ 
in Fig.\ \ref{fig:hkfiveone_A}
meets $5_1$ at a separating disk $\disk_s$.    
Let $V,V_c\subset 5_1$ be the solid tori cut off by 
the disk $\disk_s$ with the selected component $c\subset \partial \tA$ in $\partial V_c$. 
It follows from the intersection 
$\partial \tA\cap \partial\annulus$ drawn in 
Fig.\ \ref{fig:LLvariant_intersection}, where $d$
is the other component of $\partial \tA$, 
that the slope of 
$T_n(l)$ with respect to 
$(\sphere, V)$ is $\frac{n}{n+1}$, and the slope 
of $T_n(l_m)$ with respect to $(\sphere, V_c)$
is $\frac{n+1}{n}$. The theorem thence follows. 
\end{proof}
\begin{theorem}\label{teo:ann_diag_leeleetwo}
The annulus diagram of $\pairLLtwo n$ is 
\raisebox{-.1\height}{\includegraphics[scale=.25]{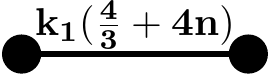}} .
\end{theorem}
\begin{proof}
Let $\annulus_c$ be the characteristic annulus of $\Compl{5_2}$
and $\annulus, \mobius$ be the type $3$-$3$ annulus and M\"obius band 
in Fig.\ \ref{fig:intersection}, respectively.
Denote by $V\subset\Compl{5_2}$ the solid torus 
cut off by $\annulus_c$. Then $\Compl{V}$ is the union 
of $5_2$ and a regular neighborhood of $\annulus\cup \mobius$ in $\Compl{5_2}$. Since 
$\annulus\cup \mobius\cup 5_2$ can be obtained by gluing $l_1,l_2$ to $l\subset \mobius$ in Fig.\ \ref{fig:union} with one from below and one from above, $\Compl V$ is a solid torus; particularly, 
the core of $V$ is a trivial knot in $\sphere$. 
%
Now observe that 
the boundary of the twisting disk $\tD$ and $cV$ 
in $\sphere$ is the link $L4a1$ in Rolfsen's link table
(see Fig.\ \ref{fig:btDandcA}), so their linking number is $\pm 2$. This, along with Theorem \ref{teo:hkfivetwo_ann_diag}, implies the core of $\annulus_{c,n}:=T_n(\annulus_c)$ 
has a slope of $\frac{4}{3}+4n$
with respect to $(\sphere V_n)$, where 
$V_n$ is the solid torus cut off by $\annulus_{c,n}$ from $\Compl{\LLtwo n}$.   

\end{proof}

\begin{remark}
It follows from Fig.\ \ref{fig:tDandcA} that the 
core of $V_n$ is a $(2n+1,2)$-torus knot.
\end{remark}

\begin{remark}
While \cite{LeeLee:12} uses the unique \emph{non-separating} annulus in $\Compl{\mathbf{L}^-_{n}}$ 
to differentiate the handlebody-knots, 
we employ the \emph{characteristic} annulus. 
\end{remark}


\section{Classification problems}\label{sec:classification}
Here we discuss to what extent \emph{the annulus diagram} and \emph{the handlebody-knot exterior} determine the knot type of a handlebody-knot. 
Let $\anndiag,\anndiagprime$ be the annulus diagrams 
of the handlebody-knots
$\pair,\pairprime$, respectively.
\subsection{Gordon-Luecke type theorems}
\begin{theorem}\label{teo:circle_stick}
If both $\anndiag$, $\anndiagprime$ are   
\raisebox{-.4\height}{\includegraphics[scale=.12]{typetwotwo_diagram2}}\ ,
where $r\in\mathbb{Q}$ and $i=1$ or $2$, and  
$\Compl\HK,\Complprime \HK$ are homeomorphic,
then $\pair, \pairprime$ are equivalent.
\end{theorem}
\begin{proof}
By the assumption, there exists a homeomorphism
$f:\Compl\HK\rightarrow \Complprime \HK$.  
Let $\annulus_1$ (resp.\ $\annulus_1'$) 
be the type $2$-$2$ annulus, and 
$\annulus_2$ (resp.\ $\annulus_2'$) 
the type $3$-$2$ annulus 
corresponding to the edges of the annulus diagram. 
Since they are the unique type $2$-$2$ and $3$-$2$ 
annuli in $\Compl\HK,\Complprime \HK$. 
It may be assumed that $f(\annulus_i)=\annulus_i'$, $i=1,2$.

Let $\disk_s\subset\HK$ be a disk bounded 
by a component of $\partial \annulus_1$;
isotope $\disk_s$ so that it is 
disjoint from $\annulus_1\cup \annulus_2$. 
Then $\disk_s$ cuts $\HK$ into two solid tori $V_1,V_2$ with $\partial \annulus_2\subset V_2$ and $\partial\annulus_1\subset V_1$. 
Denote by $\disk_s'$ a disk bounded by the image $f(\partial \disk_s)$ and disjoint from $\annulus_1'\cup\annulus_2'$. 

Observe that $f$ sends a preferred longitude of $(\sphere, V_1)$ 
to a preferred longitude of $(\sphere, V_1')$ 
since $f(\annulus_1)=\annulus_1'$, 
and sends a meridian of $V_1$
to a curve in $\partial V_1'$ of slope $\frac{1}{n}$ with respect to $(\sphere,V_1')$, $n\in\mathbb{Z}$. 
Let $t:\Complprime\HK\rightarrow \Complprime\HK$ be the homeomorphism given by twisting along $\annulus_1'$ once. 
Then the composition $f_1:=t^{-n}\circ f:\Compl\HK\rightarrow \Complprime\HK$ sends a meridian of $V_1$ to a meridian of $V_1'$, so $f_1$ can be extended to a homeomorphism
$f_2$ from $\Compl {V_2}$ to $\Compl {V_2'}$. 

For the homological reason, 
$f_2$ sends a preferred longitude of 
$(\sphere,V_2)$ to a preferred longitude of $(\sphere,V_2')$, and in terms of meridians and preferred longitudes, the induced homomorphism ${f_2}_\ast$ 
on the first homology is represented by the matrix
\begin{equation}\label{eq:repre_matrix_f_2}
\begin{pmatrix}
1&0\\
k&1
\end{pmatrix},\quad k\in\mathbb{Z}. 
\end{equation} 

Let $r=\frac{p}{q}$, $p,q\in\mathbb{Z}$.
Note that $pq\neq 0$ 
by the essentiality of $\annulus_2$, $\annulus_2'$, and components of $\partial\annulus_2$ (resp.\ $\partial\annulus_2'$) 
have a slope of $pq$ with respect to $(\sphere, V_2)$ (resp.\ $(\sphere,V_2')$) if $\annulus_2$ (resp.\ $\annulus_2'$) is of type $3$-$2$i; 
otherwise, they have a slope of $\frac{q}{p}$. 
Since $f(\partial \annulus_2)=f(\partial \annulus_2')$,
by \eqref{eq:repre_matrix_f_2},
either $pq+k=pq$ or $q+kp=q$. 
This implies $k=0$, and thus $f_2$ can be 
extended to a homeomorphism from $\pair$ to $\pairprime$.
\end{proof}

\begin{theorem}\label{teo:theta}
If both $\anndiag$, $\anndiagprime$ are   
\raisebox{-.4\height}{\includegraphics[scale=.12]{typetwotwo_diagram3}}\ ,
where $\square=\snode$ or $\hnode$, and  
$\Compl\HK,\Complprime \HK$ are homeomorphic,
then $\pair, \pairprime$ are equivalent. 
\end{theorem}
\begin{proof}
Let $f$ be a homeomorphism
from $\Compl\HK$ to $\Compl{\HK'}$, and 
$\annulus, \annulus_1, \annulus_2$ (resp.\ $\annulus', \annulus_1', \annulus_2'$) 
be the type $3$-$3$ annulus and two type $2$-$2$ annuli
in $\Compl\HK$ (resp.\ $\Compl{\HK'}$), respectively.
Since $\annulus\cup \annulus_1\cup \annulus_2$ is a characteristic
surface of $\Compl\HK$, it may be assumed that
$f(\annulus\cup \annulus_1\cup \annulus_2)=\annulus'\cup \annulus_1'\cup \annulus_2'$.
In addition, one boundary component of $\annulus_i$ (resp.\ $\annulus_i'$)
is separating in $\partial\Compl\HK$, $i=1,2$, 
while no boundary component of $\annulus$ (resp.\ $\annulus'$)
is separating, so we may further assume
$f(\annulus)=\annulus', f(\annulus_i)=\annulus_i'$, $i=1,2$.
Let $\disk_1, \disk_2$ (resp.\ $\disk_1', \disk_2'$) 
be disjoint disks bounded by boundary components of $\annulus_1, \annulus_2$ (resp.\ $\annulus_1', \annulus_2'$), respectively. Then
$\disk_1, \disk_2$ (resp.\ $\disk_1', \disk_2'$)
are parallel and hence cobound a $3$-ball $B$ 
in $\HK$ (resp.\ $B'$ in $\HK'$), which
cuts $\HK$ (resp.\ $\HK'$) 
into two solid tori $V,W$ (resp.\ $V',W'$).
Since $f(\annulus_i)= \annulus_i'$, $i=1,2$, one can extend 
$f$ to a homeomorphism
$f_1:\big(\Compl{V\cup W},B\big)\rightarrow 
\big(\Compl{V'\cup W'},B'\big)$.  

Note that $f_1$ sends a preferred longitude of 
$V$ (resp.\ of $W$) 
to a preferred longitude of $V'$ (resp.\ of $W'$), 
and sends a meridian of $V$ (resp.\ of $W$) 
to a curve of slope $\frac{1}{k_v}$ in $V'$
(resp.\ of slope $\frac{1}{k_w}$ in $W'$). 
Let $t_v,t_w$ be the homeomorphisms: 
$\big(\Compl {V'\cup W'}, B'\big)
\rightarrow \big(\Compl {V'\cup W'},B'\big)$
given by twisting along the disks $\annulus_1'\cup \disk_1'$,
$\annulus_2'\cup \disk_2'$, respectively.
Then the composition 
$t_v^{-k_v}\circ t_w^{-k_w}\circ f_1$ 
sends a meridian of $V$ (resp.\ of $W$)
to a meridian of $V'$ (resp.\ of $W'$) and can 
therefore be extended to a homeomorphism
between $\pair$ and $\pairprime$. 
\end{proof}

\begin{corollary}\label{cor:theta_gordon_luecke}
If the exterior of $\pair$ admits three non-isotopic, non-separating annuli, then the exterior determines the knot type of $\pair$.
\end{corollary}
\begin{proof}
By \cite[Theorem $1.5$]{Wan:22p}, its annulus diagram 
is 
\raisebox{-.4\height}{\includegraphics[scale=.12]{typetwotwo_diagram3}}\ with $\square=\snode$ or $\hnode$.
\end{proof}
 
\subsection{Non-completeness}
\begin{figure}[b]
\begin{subfigure}{.47\linewidth}
\center
\includegraphics[scale=.1]{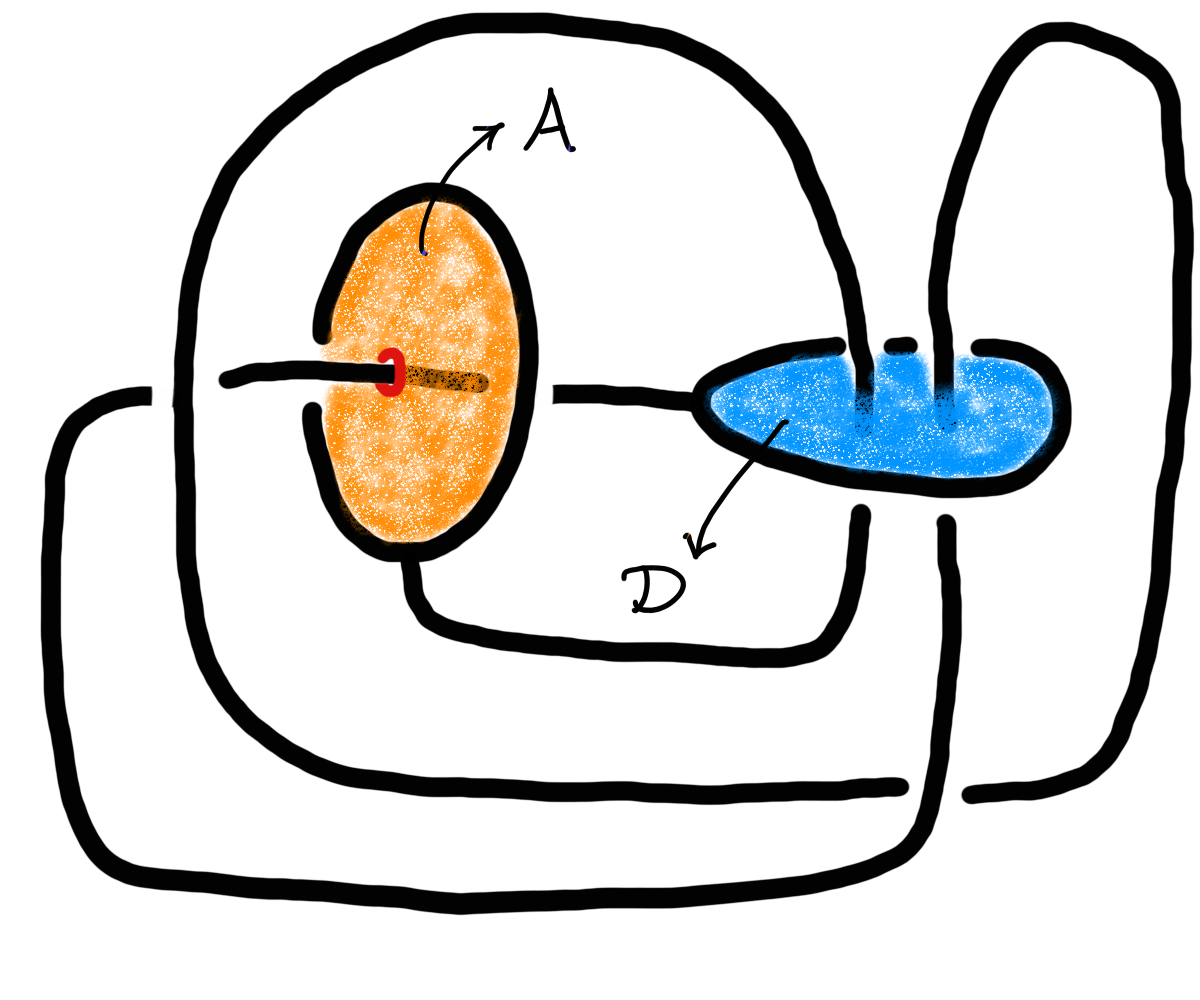}
\caption{$\pair$, $\annulus$ and $\tD$.}
\label{fig:E0}
\end{subfigure}
\begin{subfigure}{.47\linewidth}
\center
\includegraphics[scale=.1]{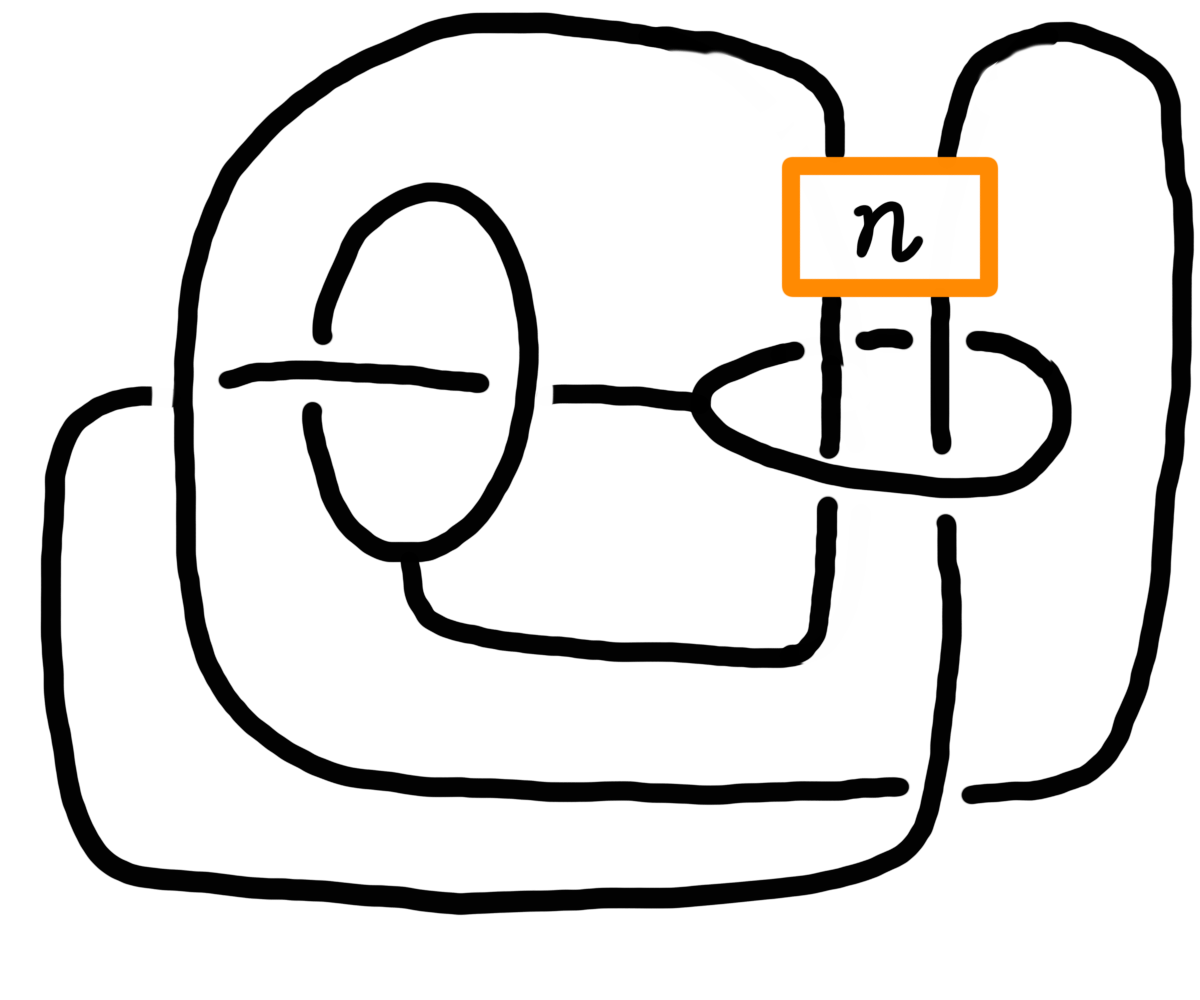}
\caption{$\pairE n$.}
\label{fig:En_family}
\end{subfigure}
\begin{subfigure}{.47\linewidth}
\center
\includegraphics[scale=.1]{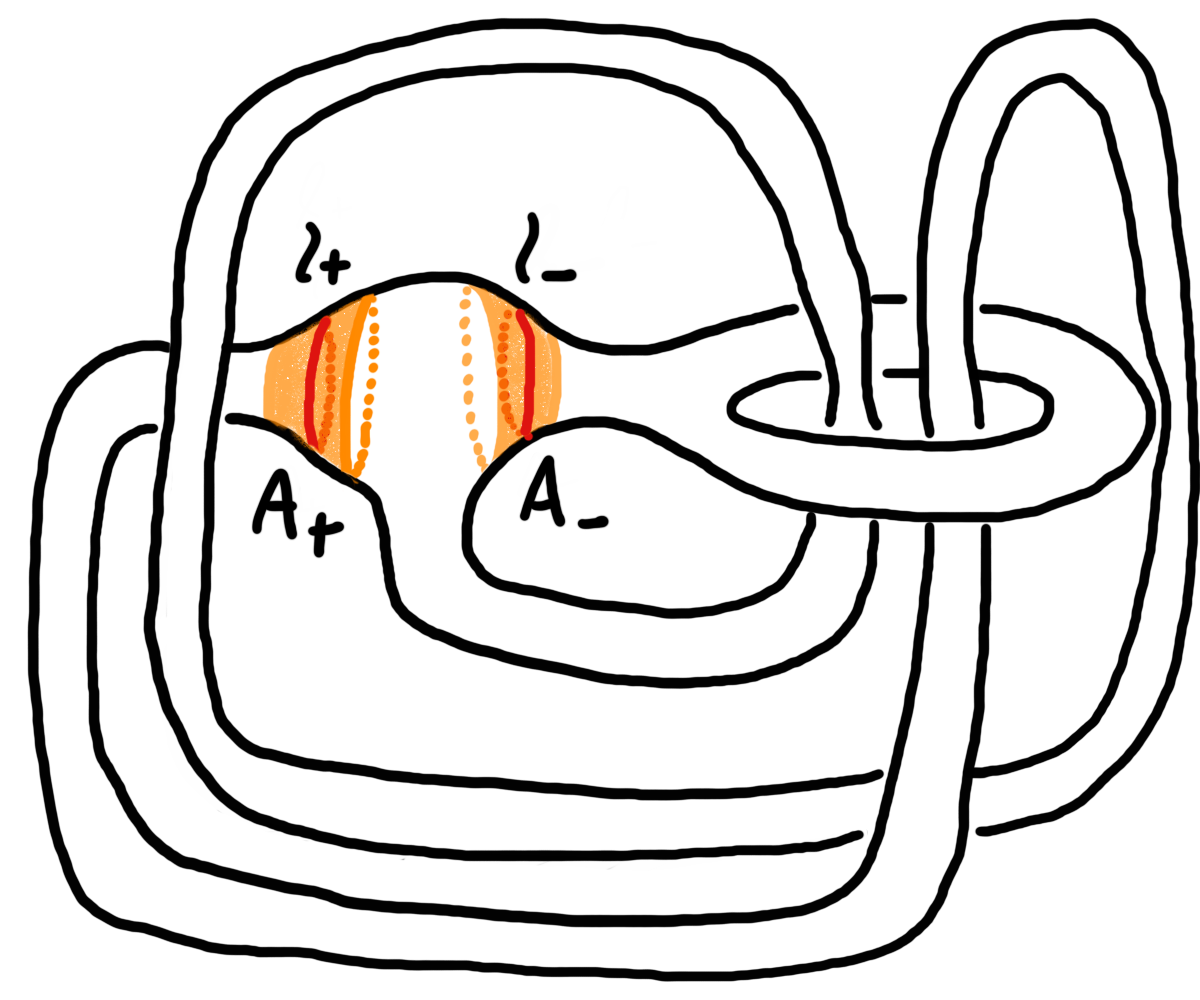}
\caption{$\pairA$ and $l_+,l_-$.}
\label{fig:unknotting_E0}
\end{subfigure}
\begin{subfigure}{.47\linewidth}
\center
\includegraphics[scale=.1]{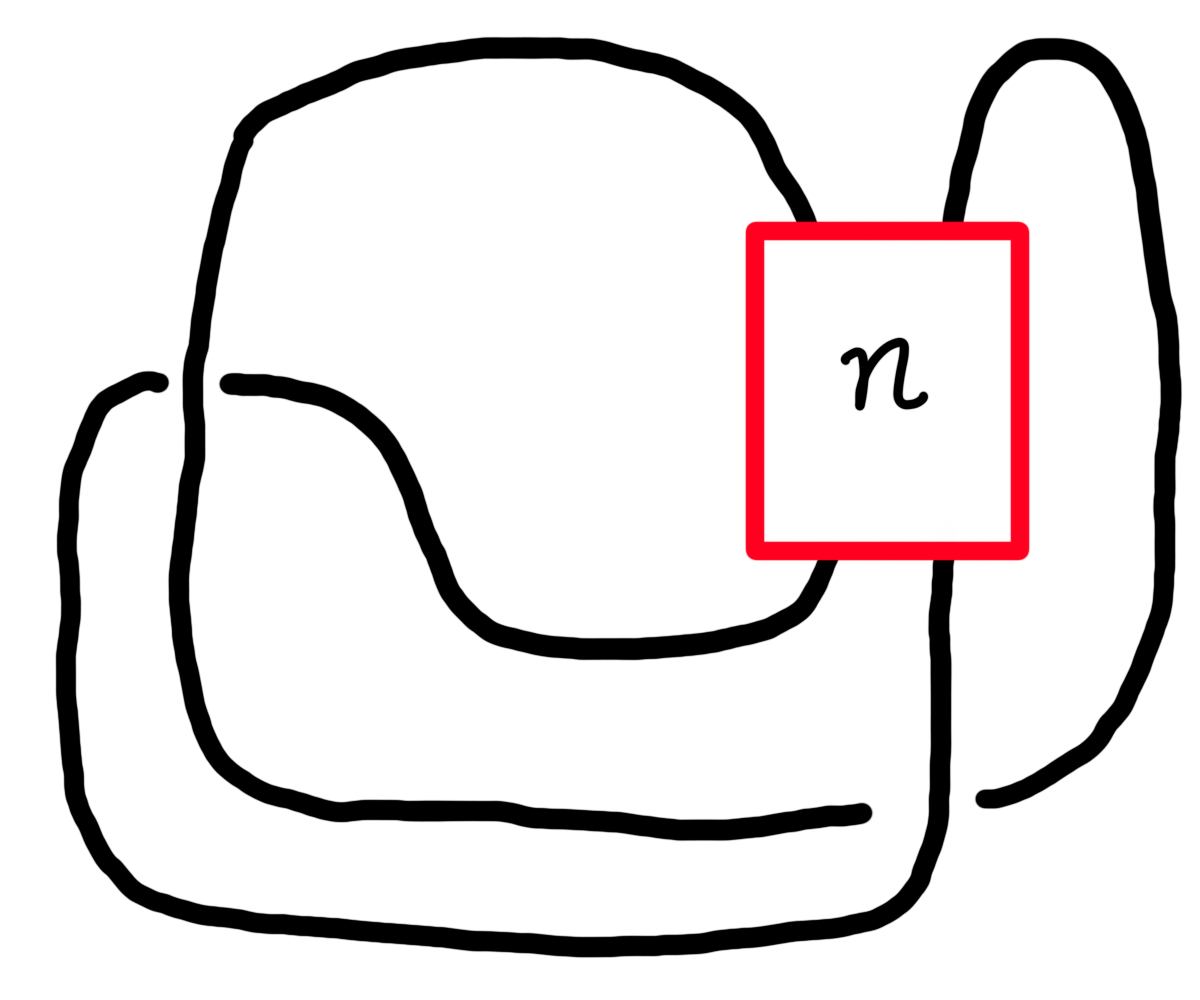}
\caption{Constituent knot $(\sphere, K_n)$.}
\label{fig:constituent_knot}
\end{subfigure}
\caption{Handlebody-knot family $\familyE$.}
\end{figure}
The annulus diagram and the handlebody-knot exterior 
are not a complete invariant, and especially so
when the exterior admits a unique essential annulus.
 
\begin{theorem}\label{teo:circle}
There exist infinitely many inequivalent atoroidal 
handlebody-knots with homeomorphic exteriors 
and the annulus diagram
\raisebox{-.4\height}{\includegraphics[scale=.12]{typetwotwo_diagram1}}. 
\end{theorem}  
\begin{proof}
Consider the handlebody-knot $\pair$ in Fig.\ \ref{fig:E0}, 
which admits a type $2$-$2$ annulus $\annulus$ and a twisting disk $\tD$. Twisting $\pair$ along $\tD$ $n$ times yields an infinite family of handlebody-knots $\mathcal{E}:=\{\pairE n\}_{n\in\mathbb{Z}}$ with $\pairE 0=\pair$ (see Fig.\ \ref{fig:En_family}). 
Denote by $\annulus_n\subset \Compl {\mathbf{E}_n}$ the image of $\annulus$ under the twisting map $T_n$ in \eqref{eq:twisting_map}; 
note that $\annulus_n$ is of type $2$-$2$, for every $n$.

To see members in $\mathcal{E}$ are atoroidal, 
we observe that
$\annulus$ is an unknotting annulus, namely $\pairA$ is a trivial handlebody-knot, where $\HKA$ 
is the union of $\HK$
and a regular neighborhood $\rnbhd \annulus$ of $\annulus$ in $\Compl\HK$. The frontier of $\rnbhd{\annulus}$ 
consists of two annuli $\annulus_+, \annulus_-$, 
parallel to $\annulus$, one of which, say $\annulus_-$, separates $\partial \HKA$ 
(Fig.\ \ref{fig:unknotting_E0}). 
Since the core $l_-$ of $\annulus_-$ 
does not bound a disk in $\Compl\HKA$---it determines 
a non-trivial conjugate class in $\pi_1(\Compl\HKA)$. 
By \cite[Proposition $5.10$]{Wan:22p}, 
$\pair$ is atoroidal and 
$\annulus$ is essential. As (essentiality) atoroidality is a property of (surfaces in) a handlebody-knot exterior, 
$\pairE n$ is atoroidal with $\annulus_n$ essential,
for every $n\in\mathbb{Z}$. 

\textbf{Claim: $\annulus\subset \Compl\HK$ is the unique type $2$-$2$ annulus.}
%
Suppose there exists another type $2$-$2$
annulus $\annulus'\subset \Compl\HK$, 
which is necessarily essential by \cite[Corollary $3.18$]{KodOzaGor:15}, \cite[Lemma $5.4$]{Wan:22p}. 
Let $\slice$ be the slicing surface, namely the closure of 
$\tD-\HK\subset\Compl\HK$. By \cite[Lemma $3.10$]{Wan:22p}, 
the component $l$ of $\partial \annulus'$
not bounding a disk in $\HK$ is isotopic to $ \partial \tD$, while the other component $l_m\subset\partial \annulus'$ is isotopic to components of $\partial \slice-\partial \tD$. 
In particular, one can isotope $\annulus'$ so 
that $\partial \annulus'\cap  \partial \slice=\emptyset$. 
Choose $\annulus', \slice$   
that minimizes 
$\# \{\annulus'\cap \slice\mid\partial \annulus'\cap \partial \slice=\emptyset\}$ in their isotopy classes.

Suppose $\annulus'\cap \slice \neq \emptyset$. 
Then there exists an annulus $\annulus''\subset \annulus'$ containing 
$l_m$ such that $c:=\annulus''\cap \slice=\partial \annulus''-l_m$. 
Since $\slice$ is a disk with two open disks removed,
$c$ cuts off an annulus $\annulus_s$ from $\slice$. 
Either $\partial \tD\subset \annulus_s$ or $\partial \tD\cap \annulus_s=\emptyset$. 
In the former, $\annulus\cup \annulus''$ 
induces a type $2$-$2$ annulus
$\hat \annulus$ having less intersection with $\slice$ 
than $\annulus'$ does. Since $\partial \hat \annulus$ is parallel to $\partial A'$, by \cite[Lemma $3.10$]{Wan:22p}, they are isotopic, contradicting the minimality. 
If $\partial \tD\cap \annulus_s=\emptyset$,
then $\annulus_s\cup \annulus''$ 
and $\partial \HK$ cobound a solid torus
with the core of $\annulus''$ its a longitude; thus one can isotope $\slice$ to decrease the number of components in $\annulus'\cap \slice$, contradicting the minimality.
Therefore $\annulus'\cap \slice=\emptyset$.
Let $\annulus^\flat$ be the annulus cut off by 
$l\subset\partial \annulus'$ 
and $\partial \tD$.
Then $\annulus^\flat\cup \annulus'\cup \slice$ is a 
pair of pants
$P\subset \Compl\HK$ that separates $\Compl\HK$,
an impossibility as components in 
$\partial P$ are parallel in $\partial \HK$. This proves the claim.

Now observe that $l_-$ bounds a separating disk in $\HKA$
and hence induces a handcuff spine of $\HKA$
whose constituent link we denote by $(\sphere,L_0)$. 
Let $K_0$ be
the component of $L_0$ dual to a 
meridian disk
bounded by the core $l_+$ of $\annulus_+$. 
Likewise the type $2$-$2$ annulus $\annulus_n\subset\Compl {\mathbf{E}_n}$
induces a trivial handlebody-knot $\big(\sphere,\HK_{\annulus_n}\big)$, where 
$\HK_{\annulus_n}:=
\mathbf{E}_n\cup \rnbhd{\annulus_n}$, 
and the non-separating component
$\annulus_{n+}\subset \partial\HK_{\annulus_n}$ of the frontier
of $\rnbhd{\annulus_n}$ in $\Compl {\mathbf{E}_n}$
induces a knot $(\sphere,K_n)$
(see Fig.\ \ref{fig:constituent_knot}). 

Now if $\pairE n,\pairE m$ are equivalent, 
then the uniqueness of $\annulus_n, \annulus_m$ implies 
there exists a homeomorphism
$f:\pairE n\rightarrow \pairE m$ 
sending $\big(\rnbhd{\annulus_n}, \annulus_{n+}\big)$ to $\big(\rnbhd{\annulus_m}, \annulus_{m+}\big)$, and 
hence $f$ induces an equivalence
between $(\sphere,K_n)$ and $(\sphere,K_m)$. On the other hand, when $n>0$, the diagram in Fig.\ \ref{fig:constituent_knot} is reduced and alternating, 
so the crossing number of $(\sphere,K_n)$ is
$n+2$ by the Tait conjecture 
(see \cite[Chap.\ $5$]{Lic:97}). 
This implies members in $\{\pairE n\}_{n\in\mathbb{N}\cup \{0\}}$ are mutually inequivalent. 
It then follows from \cite[Theorem $1.4$]{Wan:22p}
and Theorem \ref{teo:circle_stick} 
that the annulus diagram 
of $\pairE n$ is 
\raisebox{-.4\height}{\includegraphics[scale=.12]{typetwotwo_diagram1}}, for every $n$.
\end{proof}     




\section{Appendix: equivalences}
\label{sec:equivalences}


\end{document}